\documentclass[a4paper,11pt]{article}

\usepackage{amsmath}
  \usepackage{fullpage}

\usepackage{amsthm,amssymb}

\usepackage{tikz}
\usepackage{tikz-cd}
\usepackage{xcolor}

\definecolor{dgreen}{HTML}{008F00}
\definecolor{dgray}{HTML}{555555}
\definecolor{dmagenta}{HTML}{8B008B}
\definecolor{dblue}{rgb}{0.2,0.2,0.7}
\definecolor{Cerulean}{rgb}{0.0, 0.48, 0.65}

\usepackage{amssymb}
\usepackage{latexsym}
\usepackage{enumerate}

\newtheorem{Theorem}{Theorem}[section]

\newtheorem{Definition}[Theorem]{Definition}
\newtheorem{Proposition}[Theorem]{Proposition}
\newtheorem{Lemma}[Theorem]{Lemma}
\newtheorem{Corollary}[Theorem]{Corollary}
\newtheorem{Remark}[Theorem]{Remark}

\begin{document}

\title{Representations of the modular group 
into the isometries of $\mathrm{SL}_3(\mathbb{R})/\mathrm{SO}(3)$}

\author{Joan Porti\footnote{The author acknowledges support by grant
PID2024-157757NB-I00 funded by MICIU/AEI/10.13039/501100011033 and by ERDF/EU,  and by the
Mar\'\i a de Maeztu Program CEX2020-001084-M
}} 

\date{\today}
\maketitle

\begin{abstract}
We describe the topology of the variety of characters  of the modular group
$\mathrm{PSL}_2(\mathbb{Z})$ into the isometry group of the 
 symmetric space
$\mathrm{SL}_3(\mathbb{R})/\mathrm{SO}(3)$. 
This topology is determined 
in terms of the fixed point sets of the 
generators of the modular group. 
We also find Anosov representations on each nontrivial 
component.

% This connected component contains the family of representations constructed by 
% Schwartz via Pappus' theorem, 
% as well as their Anosov deformations studied by Barbot,
% Lee, and Val\'erio.
% I show that certain  representations in this component  are Anosov.
% The main results of this paper were previously proved by Schwartz \cite{Schwartz2025}, though with different arguments. I was unaware of \cite{Schwartz2025} when I posted the first version of this paper.
\end{abstract}

\section{Introduction}
\label{Section:introduction}

This article investigates the character variety of representations of the modular group
\begin{equation*}
\Gamma =\mathrm{PSL}_2(\mathbb Z)\cong \langle a,b\mid a^2=b^3=1\rangle    
\end{equation*}
into the isometry group of the  symmetric space
$
X=\mathrm{SL}_3(\mathbb{R})/\mathrm{SO}(3),
$
that we denote by
$$
\mathcal{X}(\Gamma, \mathrm{Isom}(X))
.
$$
Our purpose is to give a description of this character variety, as a concrete rank-two example, in which the topology of the components, the Fuchsian and elliptic loci, and the occurrence of Anosov representations can all be analyzed geometrically.

We first determine the topology of $\mathcal{X}(\Gamma, \mathrm{Isom}(X))
$.
Since $\Gamma\cong \mathbb{Z}_2 * \mathbb{Z}/3\mathbb{Z}$, the connected components of
\(
\mathcal{X}(\Gamma, \mathrm{Isom}(X))
\)
are naturally in bijection with the Cartesian product of the sets of conjugacy classes of elements of order $2$ and of elements of order $3$ in $\mathrm{Isom}(X)$.

% In $\mathrm{Isom}(X)$ there are precisely two conjugacy 
% classes of elements of order $3$, represented by the 
% identity and by a rotation of order $3$ in $\mathrm{SO}
% (3)$. Furthermore, there are four conjugacy classes of 
% elements of order $2$: two orientation-preserving classes, 
% the identity and a rotation of order $2$ in 
% $\mathrm{SO}(3)$, 
% and two orientation-reversing classes, namely 
% inversions (with a unique fixed point) and generalized 
% reflections with respect to a totally geodesic hyperbolic 
% plane $\mathbb{H}^2$ congruent to
% \(
% \mathrm{SO}(1,2)\big/ \bigl(\mathrm{SO}(1,2) \cap 
% \mathrm{SO}(3)\bigr).
% \)
% Consequently, $\mathcal{X}(\Gamma, \mathrm{Isom}(X))$ has 
% a total of $8$ connected components.

Within $\mathrm{Isom}(X)$ there exists exactly one conjugacy 
class of nontrivial elements of order $3$, represented by 
a rotation of angle  $\frac{2\pi}{3}$ in $\mathrm{SO}
(3)$. Moreover, there are precisely three nontrivial conjugacy classes of 
elements of order $2$: one consisting of rotations of order $2$ in 
$\mathrm{SO}(3)$, 
and two consisting of orientation-reversing isometries, namely 
inversions (each having a unique fixed point) and generalized 
reflections with respect to a totally geodesic hyperbolic 
plane $\mathbb{H}^2$ congruent to
\[
\mathrm{SO}(1,2)\big/ \bigl(\mathrm{SO}(1,2) \cap 
\mathrm{SO}(3)\bigr).
\]
It follows that $\mathcal{X}(\Gamma, \mathrm{Isom}(X))$ has 
exactly three connected components on which both $a$ and $b$ are mapped to nontrivial elements. For $i = 0,1,2$, we define
\[
\hom_i\bigl(\Gamma, \mathrm{Isom}(X)\bigr)
\]
to be the subset consisting of all 
$\rho \in \hom\bigl(\Gamma, \mathrm{Isom}(X)\bigr)$ 
such that $\rho(b)$ is nontrivial and
\[
\rho(a) \text{ is }
\begin{cases}
    \text{a rotation of order } 2, & \text{if } i = 0,\\[2pt]
    \text{an inversion,} & \text{if } i = 1,\\[2pt]
    \text{a reflection with respect to }\mathbb{H}^2, & \text{if } i = 2.
\end{cases}
\]
We denote the corresponding real GIT quotients (in the sense of \cite{RichardsonSlodowy}), taken with respect to the conjugation action of $\mathrm{Isom}(X)$, by
\[
\mathcal{X}_i\bigl(\Gamma, \mathrm{Isom}(X)\bigr)
=
\hom_i\bigl(\Gamma, \mathrm{Isom}(X)\bigr)
/\!/\mathrm{Isom}(X).
\]
The resulting character varieties carry a natural 
structure of semi-algebraic sets. From a topological 
point of view, 
$\mathcal{X}_i\bigl(\Gamma, \mathrm{Isom}(X)\bigr)$ 
is identified with the space of closed conjugacy orbits in 
$\hom_i\bigl(\Gamma, \mathrm{Isom}(X)\bigr)$. 
Further details concerning the structure of the (real) 
character varieties are provided in Section~\ref{Section:varietyofCharacters}.

\begin{Theorem} 
\label{Theorem:topology}
We have homeomorphisms
\begin{align*}
\mathcal{X}_0\bigl(\Gamma, \mathrm{Isom}(X)\bigr)
&\cong [0,1] \times [0,+\infty),\\
\mathcal{X}_1\bigl(\Gamma, \mathrm{Isom}(X)\bigr)
&\cong \mathbb{R}^3, \\
\mathcal{X}_2\bigl(\Gamma, \mathrm{Isom}(X)\bigr)
&\cong \mathbb{R}^3 .
\end{align*}
% Furthermore, for $i=1,2$  the corresponding component 
% of characters is also the topological quotient:
% \[
% \mathcal{X}_i\bigl(\Gamma, \mathrm{Isom}(X)\bigr)
% \cong \hom_i\bigl(\Gamma, \mathrm{Isom}(X)\bigr) \big/ \mathrm{Isom}(X),
% \qquad \textrm{ for }i=1,2.
% \]
\end{Theorem}

% Being the topological quotient is equivalent to the property that all 
% orbits by conjugation in the component 
% $\hom_i\bigl(\Gamma, \mathrm{Isom}(X)\bigr)$
% are closed, for $i=1,2$. This is not the case when $i=0$, there are non-closed orbits and the topological quotient is not Hausdorff.

The topology of
$ \mathcal{X}_1\bigl(\Gamma, \mathrm{Isom}(X)\bigr) $ 
(the component where $a$ is mapped to an inversion) was 
already described by
Richard Schwartz in \cite{Schwartz2025}, motivated by previous 
work on (relative) Anosov representations constructed
from Pappus' theorem \cite{Schwartz,Schwartz24}. 
% This was my first motivation (that started before I was aware of
% \cite{Schwartz2025} and yielded to a first version of this paper that computed only $ \mathcal{X}_1\bigl(\Gamma, \mathrm{Isom}(X)\bigr) $, and will remain unpublished forever).

\bigskip

A representation  in
$
\hom\bigl(\Gamma, \mathrm{Isom}(X)\bigr)
$
is called 
\emph{elliptic} if it has a global fixed point in $X$ and
\emph{Fuchsian} if it preserves a totally geodesic hyperbolic plane $\mathbb H^2\subset X$, up to homothety. 
Among Fuchsian representations we distinguish:
\begin{enumerate}[(I)]
    \item Fuchsian of type I if the invariant plane is congruent to 
$\mathrm{SL}_2(\mathbb{R})/\mathrm{SO}(2)\subset X$, induced by the reducible embedding 
$\mathrm{SL}_2(\mathbb{R})\subset \mathrm{SL}_3(\mathbb{R})$, and
    \item Fuchsian of type II if the invariant plane is congruent to 
$\mathrm{SO}(1,2)/(\mathrm{SO}(3)\cap \mathrm{SO}(1,2))\subset X$.
\end{enumerate}
We also refer to a totally geodesic plane that is isometric to $\mathbb{H}^2$ (up to homothety) as a \emph{Fuchsian plane}.
For ease of exposition, we systematically consider hyperbolic planes only up to homothety.

We denote by
$$
\mathcal{E}, \mathcal{F}_{\mathrm{I}}, \mathcal{F}_{\mathrm{II}}\subset \mathcal{X}\bigl(\Gamma, \mathrm{Isom}(X)\bigr)
$$
the loci of characters of elliptic and Fuchsian representations of type I and II, respectively.
We shall prove that the union 
$\mathcal{E}\cup  \mathcal{F}_{\mathrm{I}}\cup
\mathcal{F}_{\mathrm{II}}$ is a properly embedded graph in 
the nontrivial components of
$\mathcal{X}\bigl(\Gamma, \mathrm{Isom}(X)\bigr)$. 
Each component of $\mathcal E$ is either a point or a 
segment, though 
components of $\mathcal{F}_{\mathrm{I}}$ and 
$\mathcal{F}_{\mathrm{II}}$ are proper half-lines.
By analyzing the structure of the variety of characters 
we prove the following:

\begin{Proposition}
\label{Proposition:NonEllipticNonFuchsianSmooth}
   Away from $\mathcal{E}\cup  \mathcal{F}_{\mathrm{I}}\cup
\mathcal{F}_{\mathrm{II}}$,  
   $\mathcal{X}\bigl(\Gamma, \mathrm{Isom}(X)\bigr)$
   has a smooth analytic structure provided by $\mathbb{R}$-GIT.
\end{Proposition}

To determine the topology of the variety of characters
we use the fact that 
$\mathcal{X}\bigl(\Gamma, \mathrm{Isom}(X)\bigr)$ is
homeomorphic to the space of \emph{closed} orbits by conjugation in 
the variety of representations 
$\hom\bigl(\Gamma, \mathrm{Isom}(X)\bigr)$.
Since both $a$ and $b$ have finite order for any 
$\rho\in \hom\bigl(\Gamma, \mathrm{Isom}(X)\bigr)$, 
the fixed point sets $\mathrm{Fix}(\rho(a))$ 
and $\mathrm{Fix}(\rho(b))$
are nonempty. The first step in the proof is then:

\begin{Proposition}
\label{Proposition:mindis}
The orbit by conjugation of  
$\rho\in \hom(\Gamma,\mathrm{Isom}(X))$ is closed
if and only if there exists $x\in \mathrm{Fix}(\rho(a))$
that minimizes the distance between $\mathrm{Fix}(\rho(a))$ and $\mathrm{Fix}(\rho(b))$.
\end{Proposition}

% We establish that every non-closed orbit is contained in
% the component
% $\hom_0(\Gamma,\mathrm{Isom}(X))$
% by analyzing the ideal boundaries of 
% $\mathrm{Fix}(\rho(a))$ and $\mathrm{Fix}(\rho(b))$. In the case where 
% $\mathrm{Fix}(\rho(a))$ and $\mathrm{Fix}(\rho(b))$ are 
% asymptotic, we show that the orbit of $\rho$ accumulates on an 
% orbit for which the fixed point sets of $a$ and $b$ are parallel.

The description of the real character 
variety becomes a geometric problem in the symmetric 
space \(X\): the study, up to isometry, of configurations 
of the fixed-point sets of the two generators. More 
precisely, 
the proof of Theorem~\ref{Theorem:topology} is 
 based on a  study of the 
moduli space of pairs 
$\mathrm{Fix}(\rho(a))$ and $\mathrm{Fix}(\rho(b))$
that intersect or have a common perpendicular. 
The precise arguments vary according to the 
connected component under consideration, namely 
according to the dimension of $\mathrm{Fix}(\rho(a))$.

\medskip

The same geometric description of representations in 
terms of their fixed-point sets also allows us to 
locate Anosov representations inside these components.
Let $\rho \in \hom(\Gamma, \mathrm{Isom}(X))$
 be a representation with closed conjugation orbit, and such that $\rho(a)$ and $\rho(b)$ are nontrivial.
 Let $$
\gamma = \mathrm{Fix}(\rho(b)) \subset X,
$$
denote the fixed line of the rotation
$\rho(b)$.
 Define
$$
x \in \mathrm{Fix}(\rho(a)) \in X
$$
to be
the point that realizes the distance from 
$\mathrm{Fix}(\rho(a)) $ to $\gamma $, provided by Proposition~\ref{Proposition:mindis}. Namely:
$$
d(x ,\gamma)
=d(\mathrm{Fix}(\rho(a)) ,\gamma).
$$

Consider the \emph{parallel set} $P(\gamma)$, defined as 
the union of all geodesics in $X$ that are parallel to $\gamma$; 
see Definition~\ref{Def:parallel} below. 
The parallel set
$P(\gamma)$ is a totally geodesic submanifold, and it is isometric to the Riemannian product
$
P(\gamma) \cong \mathbb R \times \mathbb H^2
$, see Lemma~\ref{Lem:ParallelSet}.  

The rotation $\rho(b)$ admits two distinct families of invariant Fuchsian planes:
\begin{enumerate}[(a)]
    \item The factors $\{t_0\}\times \mathbb{H}^2$ contained in 
    the parallel set $P(\gamma)\cong 
    \mathbb{R}\times \mathbb{H}^2$, for $t_0\in \mathbb{R}$. %These are said to be of type~I.
    \item The planes orthogonal to $P(\gamma)$ at some point of $\gamma$. %These are said to be of type~II. 
\end{enumerate} 
Planes in the family (a) are
congruent to $\mathrm{SL}_2(\mathbb{R})/\mathrm{SO}(2)$
 whereas
those in the family (b) are 
congruent to 
$\mathrm{SO}(1,2)/(\mathrm{SO}(3)\cap \mathrm{SO}(1,2))$.
Planes in each family are related by 
\emph{transvections} along $\gamma$.

Let $\mathbb{H}^2$ be a $\rho(b)$–invariant Fuchsian plane and let 
$$
    \pi \colon X \to \mathbb{H}^2
$$
denote the nearest–point projection onto $\mathbb{H}^2$.
By construction $d(x, \pi(x))=d(x, \mathbb{H}^2)$.
With this notation, we show:

\begin{Theorem} 
\label{Thm:Anosov}
Let $\mathbb{H}^2 \subset X$ be a $\rho(b)$–invariant totally geodesic plane, and let $\pi \colon X \to \mathbb{H}^2$ and $x \in \mathrm{Fix}(\rho(a))$ be as  above. Then, for every  $C \geq 0$, there exists  $D = D(C) > 0$ such that, whenever
\[
d(x, \pi(x)) = d(x, \mathbb{H}^2) \leq C \qquad \text{and} \qquad d(\pi(x), \gamma) > D,
\]
the representation $\rho$ is Anosov.
\end{Theorem}

\begin{Remark}
    If $\rho$ is Fuchsian, then $d(x,\pi(x)) = 0$.
 In this case, 
    Theorem~\ref{Thm:Anosov} specializes to the following classical statement: a nontrivial 
    representation of $\Gamma$ 
    into $\mathrm{Isom}(\mathbb{H}^2)$
    is convex cocompact
 provided that the fixed point sets of $a$ and $b$ are sufficiently far apart. Since the space of Anosov representations forms an 
open subset, Theorem~\ref{Thm:Anosov} establishes a geometric property for this open neighborhood of the Fuchsian locus.
\end{Remark}

All components $\mathcal{X}_i(\Gamma,\mathrm{Isom}(X))$
for $i=0,1,2$
have nonempty  Fuchsian locus (which is non-compact), so 
every component has Anosov representations.

In \cite{Schwartz,Schwartz24}, Schwarz constructs a family of representations 
\(\mathcal S \subset \mathcal X_1(\Gamma,\mathrm{Isom}(X))\), arising from Pappus’ theorem, which consists of relatively Anosov representations. The submanifold \(\mathcal S\) is homeomorphic to a disk and decomposes \(\mathcal X_1(\Gamma,\mathrm{Isom}(X))\) into two connected components, one of which is entirely comprised of Anosov representations \cite{Schwartz2025}. Since \(\mathcal S\) also separates the ends of the Fuchsian locus, Theorem~\ref{Thm:Anosov} guarantees the existence of Anosov representations in each of the two connected components of the complement of \(\mathcal S\).

To prove Theorem~\ref{Thm:Anosov} we use the asymptotic methods 
and ideas of \cite{KLP25}, 
% a local-to-global principle, 
that rely 
on a local Morse lemma for symmetric spaces.

\paragraph{Organization of the paper:}
Sections~\ref{Sec:GeometryofX}
and~\ref{Section:asympotic}
provide some basic material on the symmetric space 
$X$, most of it standard, though we focus on some 
properties that we require. 
Section~\ref{Section:varietyofCharacters}
is devoted to the variety of characters, starting from 
the definition and basic properties 
to more specific results in our case, in particular we 
prove 
Propositions~\ref{Proposition:NonEllipticNonFuchsianSmooth}
and~\ref{Proposition:mindis}.
The topology of the components as well as  
the elliptic and Fuchsian locus 
is described in the next sections: 
$\mathcal{X}_1(\Gamma,\mathrm{Isom}(X))$ in 
Section~\ref{Section:X1}, 
$\mathcal{X}_0(\Gamma,\mathrm{Isom}(X))$ in 
Section~\ref{Section:X0}, and
$\mathcal{X}_2(\Gamma,\mathrm{Isom}(X))$ in 
Section~\ref{Section:X2}.
Finally Section~\ref{Sec:Anosov} is devoted to Theorem~\ref{Thm:Anosov} on Anosov representations.

\section{The geometry of the symmetric space
% $X=\mathrm{SL}_3(\mathbb{R})/\mathrm{SO}(3)$
$X$}
\label{Sec:GeometryofX}

In this section, we review the fundamental properties of 
the symmetric space $X$. The material is largely standard, 
with the possible exception of 
Subsections~\ref{Subsection:Fuchsian} 
and~\ref{Subsection:Rotations}, which concern Fuchsian 
planes and rotations. For these subsections we are not 
aware of an appropriate reference in the literature, 
although the arguments remain elementary.

\subsection{The symmetric space
$X=\mathrm{SL}_3(\mathbb{R})/\mathrm{SO}(3)$}

We recall some standard facts and fix notation; see
\cite[Appendix 5]{BGS}, 
 \cite[Chapter II.10]{BridsonHaefliger}, or \cite[Chapter~2]{Eberlein}. The symmetric space $X$ is a rank-two 
 Cartan-Hadamard manifold. 
 At the base-point 
$$
T_{\mathrm{[Id]}} X\cong\mathfrak p=
\{a\in \mathfrak{sl}_3\mathbb R\mid a= a^{\mathrm T}\},
$$
and the exponential map $
\exp\colon \mathfrak{p}\to X
$
is a diffeomorphism.
Totally geodesic submanifolds through 
    $[\mathrm{Id}]$
    are precisely subspaces of the form 
    $\mathrm{exp}(\mathfrak{h}\cap \mathfrak{p})$, 
    where $ \mathfrak{h}\subset\mathfrak{sl}_3\mathbb{R}  $ is a subalgebra \cite{Helgason}.
We set $x_0=[\mathrm{Id}]$. The inversion at $x_0$
is
$$
i_{x_0}([A])=[A^*], \qquad  A^*=(A^{\mathrm{T}})^{-1}=(A^{-1})^{\mathrm{T}}.
$$
{Maximal flats} have dimension two and  
$\mathrm{Isom}(X)$ acts transitively on
them.

\begin{Definition}
\label{Def:regular}
 A  line in $X$ is called \emph{regular} if it is
 contained in a single maximal flat. Otherwise it is called \emph{singular}.
\end{Definition}

A singular  line is contained in infinitely many maximal flats.
To visualize maximal flats and lines, up to the action of the
isometry group a maximal flat is represented by
$$
\exp \begin{pmatrix}
 \lambda_1 & 0 & 0 \\
 0 &  \lambda_2 & 0 \\
 0 & 0 & \lambda_3
\end{pmatrix}
\qquad \textrm{with } \lambda_1, \lambda_2, \lambda_3\in\mathbb R \textrm{ and } \lambda_1+\lambda_2+\lambda_3=0.
$$
Here a line through the origin is singular iff
$\lambda_i=\lambda_j$ for some $i\neq j$. 
Figure~\ref{Figure:MaxFlat} represents a maximal flat
and the three singular lines through a point.

Notice that when a
line is singular, the centralizer in $\mathrm{SO}(3)$
of the corresponding matrices 
(with $\lambda_i=\lambda_j$ for some $i\neq j$) 
is nontrivial. The action of this centralizer  
provides a one parameter family of maximal flats containing 
the   line. 

\begin{figure}[ht]
\begin{center}
  \begin{tikzpicture}[line join = round, line cap = round, scale=.7]
  \draw[very thick] (-2,0)--(2,0);
     \draw[very thick, rotate=60] (-2,0)--(2,0);
     \draw[very thick, rotate=-60] (-2,0)--(2,0);
     \node at (3,0) {${\lambda_2=\lambda_3}$};
      \node at (1.2,2) {${\lambda_1=\lambda_2}$};
      \node at (2,-2) {${\lambda_1=\lambda_3}$};
%   \node[black] at (1.5,1) {${\Delta}$};
%   \node at (6,2) {${\lambda_1+\lambda_2+\lambda_3=0}$};
  \end{tikzpicture}
  \end{center}
  \caption{The  plane {${\lambda_1+\lambda_2+\lambda_3=0}$} with the three singular lines through a point.}
  \label{Figure:MaxFlat}
  \end{figure}

\subsection{Parallel sets}
\label{SubSection:Parallel}

 \begin{Definition}
 \label{Def:parallel}
  For a line $\gamma\subset X$,  the \emph{parallel set} of $\gamma$
  is the union of all lines in $X$ parallel to $\gamma$. 
  It is denoted by   $P(\gamma)$.
 \end{Definition}

 When two lines are parallel, then  they lie in a common maximal flat, hence $P(\gamma)$ is also the union of all maximal flats containing $\gamma$:
 $$
 P(\gamma)= \bigcup_{\mu\textrm{ parallel to }\gamma} \mu
 =\bigcup_{{\gamma\subset F}} F
 $$
where the last union runs on maximal flats $F$ containing $\gamma$.

When $\gamma\subset X$ is a regular line, then $P(\gamma)$ is precisely
the maximal flat containing $\gamma$. The situation becomes nontrivial when
$\gamma$ is singular:

\begin{Lemma}
\label{Lem:ParallelSet}
 If $\gamma\subset X$ is a singular line, then $P(\gamma)$ is 
 a totally geodesic submanifold isometric to the product $\mathbb{R}\times \mathbb H^2$.
\end{Lemma}

See \cite{Eberlein} or \cite{LeebBonnerMS} for a proof.
Up to congruence and up to reparameterization,  
we may assume that  the singular geodesic 
$\gamma$
is 
$
\exp( t \, p_0)
$, $t\in\mathbb{R}$,
for 
\begin{equation}
    \label{eqn:p0}    
p_0=\begin{pmatrix}
    2&0&0\\
    0&-1&0\\
    0&0&-1
\end{pmatrix}.
\end{equation}
Consider the subalgebra
$\mathfrak{h}\subset \mathfrak{sl}_3\mathbb R$ of matrices
of the form
\begin{equation}
 \label{eqn:h}
\begin{pmatrix}
 * & 0 & 0 \\
 0 & * & * \\
 0 & * & *
\end{pmatrix},
 \end{equation}
so $\mathfrak{h}\cong\mathbb R\times \mathfrak{sl}_2\mathbb R$.
Then it is easy to check that
$$
P(\gamma)= \exp(\mathfrak{h}\cap \mathfrak{p})=\exp(\langle p_0,p_1,p_2\rangle),
$$
where $p_0$ is as in \eqref{eqn:p0} and 
\begin{equation}
    \label{eqn:p}    
p_1=\begin{pmatrix}
    0&0&0\\
    0&\frac12&0\\
    0&0&-\frac12
\end{pmatrix},
\quad
p_2=\begin{pmatrix}
    0&0&0\\
    0&0&\frac12\\
    0&\frac12&0
\end{pmatrix}.
\end{equation}
Here $p_0$ is tangent to the factor $\mathbb R$, and $p_1$ and $p_2$, to 
$\mathbb H^2$.

Consider the nearest point projection $\pi\colon X\to P(\gamma)$.

\begin{Lemma}
\label{Lem:Fiber}
For a singular line $\gamma\subset X$,
 each fiber of the nearest point projection $\pi\colon X\to P(\gamma)$
 is a totally geodesic 
 hyperbolic plane (up to homothety)
 congruent to 
 $\mathrm{SO}(1,2)/(\mathrm{SO}(3)\cap \mathrm{SO}(1,2))$.
\end{Lemma}

\begin{proof}
 By considering the action of isometries of $X$
 that preserve $P(\gamma)$, it suffices to discuss the fiber
 at the class of the identity.
 By orthogonality, this fiber is
 $\exp(\mathfrak{h}^\perp\cap \mathfrak{p} )$ where
 $\mathfrak{h}^\perp\cap \mathfrak{p}$ is the linear span of
 \begin{equation}
     \label{eqn:f}
     f_1=\begin{pmatrix}
  0 & 1& 0 \\
  1 & 0 & 0 \\
  0  & 0 & 0
 \end{pmatrix}
\quad\textrm{ and }\quad
 f_2=\begin{pmatrix}
  0 & 0& 1 \\
  0 & 0 & 0 \\
  1  & 0 & 0
 \end{pmatrix}.
\end{equation}
The matrices $f_1$ and $f_2$ together with
$$
[f_1,f_2]=
 \begin{pmatrix}
  0 & 0& 0 \\
  0 & 0 & 1 \\
  0  & -1 & 0
 \end{pmatrix}
$$
span the Lie subalgebra $\mathfrak{so}(1,2)$.
 So
 $\exp(\mathfrak{h}^\perp\cap \mathfrak{p} )$ is the totally geodesic hyperbolic plane $\mathrm{SO}(1,2)/\mathrm {SO}(2)\subset X$.
\end{proof}

With the notation of~\eqref{eqn:p0}, \eqref{eqn:p}, and~\eqref{eqn:f},
$$
\mathfrak{p}=\langle p_0, p_1, p_2, f_1, f_2 \rangle.
$$
We chose this notation  because $\langle p_0, p_1, p_2\rangle$ is tangent to the parallel set,
and $\langle f_1, f_2\rangle$, to the fiber of the nearest point projection.

\subsection{Fuchsian planes}
\label{Subsection:Fuchsian}

\begin{Definition}
A \emph{Fuchsian plane} is a totally geodesic submanifold 
\(\mathbb{H}^2 \subset X\) which is, up to homothety, isometric to the hyperbolic plane.
\end{Definition}

We say that two submanifolds on $X$ are \emph{congruent} if they are related by an isometry.

\begin{Proposition}
\label{Proposition:2Fuchsian}
There are two congruence classes of Fuchsian planes $\mathbb H^2\subset X$: 
\begin{enumerate}[(I)]
    \item 
$
\mathbb{H}^2_{\mathrm{I}} \; =\; \mathrm{SL}_2(\mathbb{R})/\mathrm{SO}(2) \;\subset\; X,
$
coming from the   reducible
embedding
$
\mathrm{SL}_2(\mathbb{R}) \subset \mathrm{SL}_3(\mathbb{R}).
$
\item 
$
\mathbb{H}^2_{\mathrm{II}} \;=\; \mathrm{SO}(1,2)\big/\bigl(\mathrm{SO}(3)\cap \mathrm{SO}(1,2)\bigr)
\;\subset\; X.
$
\end{enumerate}
\end{Proposition}

\begin{proof}
    Using \cite{Helgason}, for instance, the proof
    follows from the classification of representations of
    $\mathrm{SL}_2(\mathbb R)$ into $\mathrm{SL}_3(\mathbb R)$.
\end{proof}

We call these planes of type I and type II, 
respectively, up to congruence.
% We call the planes $\mathbb{H}^2_{\mathrm{I}}$
% and $\mathbb{H}^2_{\mathrm{II}}$ of type I and II respectively,
% and up to congruence.
Notice that type~I is not congruent to type~II 
as they have different sectional curvatures, 
$\operatorname{sec}(\mathbb{H}^2_{\mathrm{I}})
= 4 \operatorname{sec}(\mathbb{H}^2_{\mathrm{II}})$.
They are also distinguished by their ideal boundaries, see
Lemma~\ref{Lemma:boundaryFuchsian}.

\begin{Lemma}
\label{Lemma:UniqueParallelSet}
\begin{enumerate}[(a)]
    \item A Fuchsian plane of type I is contained in a unique parallel set. 
    \item A Fuchsian plane $\mathbb H^2$ of type II   is perpendicular to a unique parallel set  at each point $x\in\mathbb H^2$.
\end{enumerate}    
\end{Lemma}

\begin{proof}
    By homogeneity, except for uniqueness  the lemma is clear by construction. To prove uniqueness in (a) notice that each geodesic line in a Fuchsian plane of type I is regular, so it is contained in a unique maximal flat. The parallel set is the union of those maximal flats, hence unique. 
    For uniqueness of (b), assume by 
    Proposition~\ref{Proposition:2Fuchsian} that 
    the Fuchsian plane is $\exp(\langle f_1,f_2\rangle)$. 
    Then $(T_{x_0}\mathbb H^2_{\mathrm{II}})^\perp=\langle f_1,f_2\rangle^\perp=\mathfrak{p}\cap\mathfrak h$, for $\mathfrak h$ as in \eqref{eqn:h}. Hence 
    $\exp (\mathfrak{p}\cap\mathfrak h) $ is the unique parallel
    set perpendicular to $H^2_{\mathrm{II}}$ at $x_0$, as 
    a totally geodesic submanifold is determined by its tangent
    space at one point.
\end{proof}

\begin{Lemma}
\label{Lemma:uniquerotations}
    Given $\mathbb{H}^2\subset X$ a Fuchsian plane and $x\in \mathbb{H}^2$, there are only two rotations of order 
    three
    $r$ and $r^2$ that preserve both $\mathbb{H}^2$
    and $x$.
\end{Lemma}

\begin{proof}
Assume that $\mathbb{H}^2$ is of type I. Then $\mathbb{H}^2$ appears as the factor of a parallel set 
\( P \cong \mathbb{R} \times \mathbb{H}^2 \)
that is invariant under the action of $r$, by the uniqueness statement in Lemma~\ref{Lemma:UniqueParallelSet}. Since 
$r$ preserves 
the product decomposition 
\( P \cong \mathbb{R} \times \mathbb{H}^2 \), 
the fixed point set of the rotation $r$ 
coincides with one of the parallel geodesics determining the
$\mathbb{R}$–factor, namely with the parallel geodesic 
containing $x$. 
This geodesic axis uniquely determines the rotations $r$ and $r^2$.

If $\mathbb{H}^2$ is of type II, then it is orthogonal to a unique parallel set passing through $x$, and the desired conclusion follows directly from the preceding argument.
\end{proof}

\subsection{Isometric involutions of $X$}
\label{subsection:isometricinvolutions}
The basic isometric involutions of $X$ are inversions
with respect to a point $x\in X$, but there are more 
involutions that we describe next.

Start by recalling that the set of 
orientation-preserving isometries of $X$ can be identified to 
$\mathrm{SL}_3(\mathbb{R})$. If we denote elements in 
$X=\mathrm{SL}_3(\mathbb{R})/\mathrm{SO}(3)$ 
by equivalence classes of matrices 
$[M]\in X$, where $M\in \mathrm{SL}_3(\mathbb{R})$, 
then the 
isometry corresponding to $A\in \mathrm{SL}_3(\mathbb{R})$
maps $[M]\mapsto [A\,M]$.

To describe orientation-reversing isometries, we consider the 
\emph{contragredient} of a matrix $A$, denoted by 
$$
    A^*=(A^{-1})^{\mathrm{T}}=(A^{\mathrm{T}})^{-1},
$$
namely the composition of transposition and inverse.
From the isomorphism $\mathrm{Isom}^+(X)\cong\mathrm{SL}_3(\mathbb{R})$
we get:

\begin{Lemma}
    \label{Lemma:OR}
    Every orientation-reversing isometry of $X$ can be written as
    $$
    [M]\to [A\, M^*], \qquad \textrm{ for every class }
    [M]\in X,
    $$
    for some $A\in \mathrm{SL}_3(\mathbb{R})$.
\end{Lemma}

For instance, by taking $A=\mathrm{Id}$ in 
Lemma~\ref{Lemma:OR}, we get the inversion with respect to the class of the identity $[\mathrm{Id}]$.
More generally:

\begin{Remark}
The inversion with respect to a class $[B]\in X$
is the map $$
[M]\mapsto [B B^{\mathrm{T}}\, M^*],\qquad  \textrm{ for every class }[M]\in X.
$$
\end{Remark}

\begin{Lemma}
\label{Lem:4components}
The space of the isometric involutions
$$\{\phi\in\mathrm{Isom}(X)\mid \phi^2=\mathrm{Id} \}$$ 
has 4 connected  components, corresponding
to  conjugacy classes, that are the following:
\begin{enumerate}[(a)]
    \item $\phi$ is the identity.
    \item $\phi$ is orientation-preserving represented by a matrix 
    conjugate to $\mathrm{diag}(1,-1,-1)$. In particular 
    $\mathrm{Fix}(\phi)$ is a parallel set.
    \item $\phi$ is an inversion, in particular it
    is orientation-reversing and $\mathrm{Fix}(\phi)$
    is a point.
    \item $\phi$ is the reflection with respect to a 
    Fuchsian plane of type II, in particular $\phi$ is orientation-reversing and $\mathrm{Fix}(\phi)$
    is a Fuchsian plane of type II.
\end{enumerate}
\end{Lemma}

\begin{proof} When $\phi$ is orientation-preserving, it is represented by a matrix $A \in \mathrm{SL}_3(\mathbb{R})$
that satisfies $A^2=\mathrm{Id}$. Thus $A$ diagonalizes with eigenvalues $\pm 1$ and 
it falls into either case (a) or (b).

  When $\phi$ is  orientation-reversing,
  it maps each class
  $[M] \in X$ to
\[
\phi([M]) = [A M^*],
\]
for some matrix $A \in \mathrm{SL}_3(\mathbb{R})$. The condition 
$\phi^2 = \mathrm{Id}$ is equivalent to 
\(
A A^* = \mathrm{Id},
\)
which in turn can be written as \[
A = A^{\mathrm{T}}
.\] 
 It is straightforward to check that the space of symmetric matrices in
$\mathrm{SL}_3(\mathbb{R})$ consists of exactly two connected components, 
distinguished by the signature. Furthermore,
these components are conjugation orbits, because for $\phi$, $A$ and $M$ as above, and $B\in \mathrm{SL}_3(\mathbb{R})$,
$$
(B\circ\phi\circ B^{-1})([M])=[ B A B^{\mathrm{T}}\, M]
.$$
The involution $\phi$ is an inversion when $A$ is positive definite, 
and $\phi$ fixes  a Fuchsian plane of type~II
when $A$ has signature $(1,2)$.
\end{proof}

\subsection{Rotations}
\label{Subsection:Rotations}

We aim to describe the space of all possible rotations in 
$\mathrm{SO}(3)$ of angle $\frac{2\pi}{3}$.
For convenience, we fix a rotation
$$
R_{\frac{2\pi}{3}}=
\begin{pmatrix}
    1 & 0 & 0 \\
    0 & -\frac 12 & -\frac{\sqrt{3}}{2} \\
    0 & \frac{\sqrt{3}}{2} & -\frac{1}{2}
\end{pmatrix}.
$$
A general  rotation in $\mathrm{SO}(3)$ of angle $\frac{2\pi}{3}$  can be described by an oriented axis in $\mathbb R^3$, namely the space of such rotations is $S^2$.
Equivalently, we consider the conjugacy orbit of 
$R_{\frac{2\pi}{3}}$; this space is
$\mathrm{SO}(3)/\mathrm{SO}(2)\cong S^2$. We are 
interested in an explicit description of the homeomorphism 
with~$S^2$:

\begin{Lemma}
\label{Lemma:3rotations}
The map
\begin{equation}
    \label{eqn:Rot}
    \begin{array}{rcl}
            \mathrm{Rot}\colon S^2=\{v\in\mathbb R^3\mid |v|=1\} & \to   & \mathrm{SO}(3) \\
         R 
         \begin{pmatrix}
             1 \\ 0\\ 0
         \end{pmatrix} & \mapsto & R \,
         % \begin{pmatrix}
         %     1 & 0 & 0 \\ 0 & -\frac 12 & -\frac{\sqrt{3}}{2} \\
         %     0 & \frac{\sqrt{3}}{2} & \frac{1}{2} 
         % \end{pmatrix}
         R_{\frac{2\pi}{3}}\,
         R^{-1} 
    \end{array}
    \qquad \textrm{for all }R\in \mathrm{SO}(3)
\end{equation}
defines an equivariant homeomorphism between $S^2$ and the set of 
$\frac{2\pi}{3}$-rotations in $\mathrm{SO}(3)$.
\end{Lemma}

Being  equivariant means that
$
\mathrm{Rot}( A v)= \mathrm{Ad}_A(\mathrm{Rot}(v))
$
% \qquad \textrm{ 
for all 
% }
$
v\in S^2,\ A\in\mathrm{SO}(3)
$.

\begin{proof}
The map is well defined because $(1\ 0 \ 0)^{\mathrm{T}}$ and $R_{\frac{2\pi}{3}}$
have the same stabilizer. It is straightforward to
check that it yields a homeomorphism between $S^2$ and the space of 
$\frac{2\pi}{3}$ rotations.  It is equivariant by construction.
\end{proof}

Notice that 
$\mathrm{Rot}(-v)=\mathrm{Rot}(v)^{-1}=\mathrm{Rot}(v)^{-2}$.

\begin{Lemma}
    \label{Lemma:axis} For \(v = (x\, y\, z)^{\mathrm{T}}\), the fixed-point set of \(\mathrm{Rot}(v)\) in \(\mathfrak{p}\) is the one-dimensional subspace spanned by
\[
\mathrm{Axis}(\mathrm{Rot}(v)) 
= 3 v\, v^{\mathrm{T}}-\mathrm{Id}
= 3 
\begin{pmatrix}
x\\[2pt] y\\[2pt] z
\end{pmatrix} (x\, y\, z)
- \mathrm{Id}
=
\begin{pmatrix}
3x^2 - 1 & 3xy      & 3xz \\[3pt]
3yx      & 3y^2 - 1 & 3yz \\[3pt]
3zx      & 3zy      & 3z^2 - 1
\end{pmatrix}.
\]
\end{Lemma}

\begin{proof}
    The assertion is clear for $v=(1\ 0\ 0)^{\mathrm{T}}$, as $\mathrm{Rot}((1\ 0\ 0)^{\mathrm{T}})=
    R_\frac{2\pi}{3}$, and apply that it is equivariant for the action of $\mathrm{SO}(3)$.
\end{proof}

Observe that the polynomial entries of
the matrix appearing in Lemma~\ref{Lemma:axis}
are essentially those defining the Veronese embedding 
of $\mathbb{P}^2$ into $\mathbb{P}^5$.
In fact 
Lemmas~\ref{Lemma:3rotations} 
and~\ref{Lemma:axis}
together yield a homeomorphism between $\mathbb{P}^2$ and the locus of 
singular lines in $X$ passing through the point $x_0$. Examining the 
corresponding ideal boundary, 
Proposition~\ref{Prop:Flag} below
yields the same homeomorphism.

\section{Asymptotic  properties of $X$}
\label{Section:asympotic}

In this section
we review some asymptotic properties of the symmetric space 
$X=\mathrm{SL}_3(\mathbb{R})/\mathrm{SO}(3)$. These results are a special case of the theory of non-compact symmetric spaces.

\subsection{Visual and Tits boundaries}

We recall the notation and the features of the spherical 
building at infinity that will be used below; see
\cite{Ballmann, BGS, BridsonHaefliger}
for background.
The visual boundary $\partial_\infty X$ 
 has dimension four,
 $$
 \partial_\infty X  \cong S^4.
 $$
The same set, equipped with the Tits metric 
 $$
    d_\mathrm{Tits}(\xi,\zeta)=\sup_{x\in X} \angle_x(\xi,\zeta)
    $$
is the one-dimensional spherical building 
$\partial_\mathrm{Tits}X$ of type $A_2$.

The apartments of 
$\partial_\mathrm{Tits}X$ 
are the ideal boundaries of maximal flats and are circles of length \(2\pi\). The walls divide every apartment into spherical Weyl chambers of length \(\pi/3\). If \(\sigma\) is a spherical Weyl chamber and \(x\in X\), we denote by
$$ V(x,\sigma) $$
the Weyl sector based at \(x\) and asymptotic to \(\sigma\).

\begin{figure}[ht]
\begin{center}
  \begin{tikzpicture}[line join = round, line cap = round, scale=1]
   \fill[gray!30!white, opacity=.4]  (0,0)--(2.5,0)  arc[start angle=0, end angle=60, radius=2.5]--(0,0);
  \draw[black, thick] (-2,0)--(2.5,0);
     \draw[black, thick, rotate=60] (-2,0)--(2.5,0);
     \draw[black, thick, rotate=-60] (-2,0)--(2,0);
     \node at (2.7,-.3) {${\lambda_2=\lambda_3}$};
      \node at (0.4,2.3) {${\lambda_1=\lambda_2}$};
      \node at (1.8,-2) {${\lambda_1=\lambda_3}$};
    \draw[black, very thick] (2.6, 0)  arc[start angle=0, end angle=60, radius=2.6];  
    \draw[black, fill=black] (0,0) circle[radius=.06];
    \node[black] at (3.3, 1.2) {$\sigma\subset\partial_\infty X$};
    \node at (1.4,0.7) {$V(x,\sigma)$};
    \node at (-.29,.15) {$x$};
   \node at (-2.5,2) {$F$}; 
% %   \node[black] at (1.5,1) {${\Delta}$};
% %   \node at (6,2) {${\lambda_1+\lambda_2+\lambda_3=0}$};
  \end{tikzpicture}
  \end{center}
  \caption{ The spherical Weyl chamber $\sigma\subset\partial_\infty F$ and the 
   Weyl sector $V(x,\sigma)\subset F$,  in the same picture of a 
  maximal flat $F$ as in Figure~\ref{Figure:MaxFlat}. 
  In the plane $\lambda_1+\lambda_2+\lambda_3=0$, $V(x,\sigma)$  is defined 
  by $\lambda_1\geq\lambda_2\geq \lambda_3$.}
  \label{Figure:WeylChambers}
  \end{figure}

\subsection{Type map  and ideal orbit spaces}

% \subsection{The action of $\mathrm{SL}_3(\mathbb R)$ 
% at the ideal boundaries}

The actions of both $\mathrm{SL}_3(\mathbb R)$ and 
$\mathrm{Isom}(X)$ on $\partial_\mathrm{Tits} X$ are not transitive, and 
this   yields relevant geometric tools.

Let
$$ \sigma_{\mathrm{mod}} =SL_3(\mathbb R)\backslash\partial_{\mathrm{Tits}}X \cong[0,\pi/3] $$
be the model chamber. The quotient projection defines the type map
$$ \operatorname{type}:\partial_\infty X\to\sigma_{\mathrm{mod}}. $$
A geodesic is regular precisely when the types of its endpoints lie in \(\operatorname{int}(\sigma_{\mathrm{mod}})\). Inversion induces the canonical involution
$$ \iota:\sigma_{\mathrm{mod}}\to\sigma_{\mathrm{mod}}, $$
which exchanges the endpoints of the model chamber and fixes its center.

\begin{Proposition}
    \label{Prop:Flag} We have $\mathrm{SL}_3(\mathbb{R})$-equivariant homeomorphisms between:
    \begin{enumerate}[(a)]
        \item one orbit of vertices and the projective plane
        $$
        \mathrm{SL}_3(\mathbb{R}) \cdot v_1
        % \cong\mathrm{SL}_3(\mathbb{R})/P_1
        \cong \mathbb P^2 , 
        $$
        \item the other orbit of vertices and its dual
        $$
        \mathrm{SL}_3(\mathbb{R}) \cdot v_2
        % \cong \mathrm{SL}_3(\mathbb{R})/P_2
        \cong (\mathbb P^2)^*\textrm{, and } 
        $$
        \item the set of all (the orbit of any) chambers  and the flag manifold
        $$\mathrm{SL}_3(\mathbb{R}) \cdot \sigma 
        % \cong \mathrm{SL}_3(\mathbb{R})/B 
    \cong
    \mathrm{Flag}(\mathbb P^2)=\{(p,l)\in \mathbb P^2 
    \times (\mathbb P^2)^*\mid
    p\in l\}.
    $$ 
    \end{enumerate}
   
\end{Proposition}

 With the identifications of 
 Proposition~\ref{Prop:Flag},  a spherical apartment in $\partial_{\mathrm{Tits}}X$ 
    corresponds to three points in $\mathbb P ^2$ in generic position (not aligned)
    or, equivalently, to three projective lines in generic position.
     A parallel set corresponds to a point $p\in \mathbb P^2$ and a line $l\subset \mathbb P^2$ in generic position ($p\not\in l$).

% \subsection{More asymptotic notions}

% We keep on reviewing asymptotic notions, in particular opposition and the type map.

\bigskip

Two ideal points are opposite if they are the endpoints of a geodesic; equivalently,
$$ d_{\mathrm{Tits}}(\xi_-,\xi_+)=\pi. $$
In this case
$$ \operatorname{type}(\xi_-) =\iota(\operatorname{type}(\xi_+)). $$
Under the identifications of Proposition~\ref{Prop:Flag}, a point \(p\in\mathbb P^2\) and a line \(l\in(\mathbb P^2)^*\) are opposite precisely when \(p\notin l\), and two chambers are opposite precisely when the corresponding flags are in generic position.

% \paragraph{Opposition is a genericity condition}  
We conclude this subsection by examining the notion of opposition, understood as  a generic behavior.  
The symmetric space \(X\) is not a visibility space; in particular, there exist pairs of ideal boundary points that cannot be connected by any geodesic.

\begin{Proposition}
\label{Proposition:correspond}
\begin{enumerate}[(a)]
    \item 
    Given two opposite spherical chambers $\sigma_+,\sigma_-\subset \partial_\infty X$
    there exists a unique maximal flat $F(\sigma_+,\sigma_-)$ asymptotic 
    to both of them.
    \item
    Given $p\in \mathbb P^2$ and $l\in(\mathbb P^2)^*$ in generic position ($p\not\in l$),
    there exists a unique parallel set $P(\gamma)$ such that $\gamma$ (and all lines in 
    $P(\gamma)$ parallel to $\gamma$) is asymptotic to $p$ and $l$.
    \end{enumerate}
\end{Proposition}

%
% The following remark helps to understand Proposition~\ref{Proposition:correspond}:

\subsection{Ideal boundary of Fuchsian planes}

Next we discuss some basic properties of the ideal boundary of Fuchsian planes that we shall use.

\begin{Lemma}
    \label{Lem:FuchsianRegular}
    % Let $\zeta=\mathrm{center}(\sigma_{mod})$ be the fixed point of $\iota$.
For any Fuchsian $\mathbb H^2\subset X$,
$$
\mathrm{type}(\partial_\infty\mathbb H^2)=\{\mathrm{center}(\sigma_{mod})\}.
$$
\end{Lemma}

\begin{proof}
  A totally geodesic embedding $\mathbb H^2\subset X$ 
  is induced by a
  representation  $\varrho\colon \mathrm{SL}_2(\mathbb R)\to 
  \mathrm{SL}_3(\mathbb R)$, see \cite{Helgason}.
  Hence $\varrho(\mathrm{SL}_2(\mathbb R))\subset \mathrm{SL}_3(\mathbb R)$
  acts transitively on  $\partial_\infty\mathbb H^2$ and,  as the type map is
  $\mathrm{SL}_3(\mathbb R)$-invariant, 
  the ideal boundary $\partial_\infty\mathbb H^2$ has constant 
   type. Furthermore, the embedded $\mathbb H^2\subset X$
  is invariant by inversions centered at $\mathbb H^2\subset X$, thus
  the image of $\mathrm{type}(\partial_\infty\mathbb H^2)$ is $\iota$-invariant, where $\iota\colon \sigma_{mod}
  \to \sigma_{mod}$ is the canonical involution.
\end{proof}

As it   has constant type, $\partial_\infty\mathbb H^2$ is determined by its projection to the flag manifold.

\begin{Definition}
\begin{enumerate}[(a)]
    \item A closed curve  $S\subset \mathrm{Flag}(\mathbb{P}^2)$
    is called a \emph{pencil of flags} if there exist 
    $p\in \mathbb{P}^2$ and $L\in (\mathbb{P}^2)^*$ in generic position ($p\not\in L$) such that 
    $$
    \{(q,l)\in\mathrm{Flag}(\mathbb{P}^2)
    \mid q\in L \textrm{ and } p\in l\}.
    $$
    \item Given a conic $C\subset \mathbb{P}^2$, its 
    \emph{tangent flag} is
    $$
    \{(q,l)\in\mathrm{Flag}(\mathbb{P}^2)\mid q\in C
    \textrm{ and } l\textrm{ is tangent to } C \textrm{ at }q\}.
    $$
\end{enumerate}    
\end{Definition}

\begin{figure}[ht]
\begin{center}
  \begin{tikzpicture}[line join = round, line cap = round, scale=.8]
    \draw[very thick, gray, opacity=.5] (-3.5,-2)--(3.5,-2);
    \draw[thick] (0,.5)--(0,-2.5);
    \draw[thick]  (0.5,0.5)--(-2.5,-2.5);
    \draw[thick]  (-0.5,0.5)--(2.5,-2.5);
    \draw[thick]  (-0.2,0.5)--(1,-2.5);
    \draw[thick]  (0.2,0.5)--(-1,-2.5);
    \draw[thick]  (-0.7,0.5)--(3.5,-2.5);
    \draw[thick]  (0.7,0.5)--(-3.5,-2.5);
    \draw[thick]  (-0.35,0.5)--(1.75,-2.5);
    \draw[thick]  (0.35,0.5)--(-1.75,-2.5);
       \draw[very thick] (0,-2) circle[radius= .05]; 
       \draw[very thick] (2,-2) circle[radius= .05]; 
       \draw[very thick] (-2,-2) circle[radius= .05]; 
       \draw[very thick] (.8,-2) circle[radius= .05]; 
       \draw[very thick] (-.8,-2) circle[radius= .05]; 
       \draw[very thick] (2.8,-2) circle[radius= .05]; 
       \draw[very thick] (-2.8,-2) circle[radius= .05];
       \draw[very thick] (1.4,-2) circle[radius= .05]; 
       \draw[very thick] (-1.4,-2) circle[radius= .05];  
\begin{scope}[shift={(9,-1.25)}, xscale=1.6, yscale=1.2]
       \draw[ultra thick, gray, opacity=.5] (0,0) circle[radius=1];
       \draw[thick] (1.5,1)--(-1.5,1);
       \draw[thick, rotate=30] (1.5,1)--(-1.5,1);
       \draw[thick, rotate=45] (1.5,1)--(-1.5,1);
       \draw[thick, rotate=60] (1.5,1)--(-1.5,1);
       \draw[thick, rotate=-30] (1.5,1)--(-1.5,1);
       \draw[thick, rotate=-45] (1.5,1)--(-1.5,1);
       \draw[thick, rotate=-60] (1.5,1)--(-1.5,1);
       \draw[very thick] (0,1) circle[radius= .05];
       \draw[very thick, rotate=30] (0,1) circle[radius= .05];
       \draw[very thick, rotate=45] (0,1) circle[radius= .05];
       \draw[very thick, rotate=60] (0,1) circle[radius= .05];
       \draw[very thick, rotate=-30] (0,1) circle[radius= .05];
       \draw[very thick, rotate=-45] (0,1) circle[radius= .05];
       \draw[very thick, rotate=-60] (0,1) circle[radius= .05];
%         \draw[thick, rotate=60] (0,0)--(3,0);
%         \draw[thick, dashed] (0,0)--(0,3);
%         \draw[thick, dashed, rotate=-60] (0,0)--(0,3);
%         \draw[ultra thick] (3,0) arc[start angle=0, end angle=90, radius=3];
%     \node at (0,3.5) {${\lambda_1=0}$};
%     \node at (4,0) {${\lambda_2=\lambda_3}$};
%     \node at (3,2.5) {${Q}$};
\end{scope}
  \end{tikzpicture}
  \end{center}
  \caption{A pencil of flags and the tangent flag to a conic.}
  \label{Figure:PencialAndConic}
  \end{figure}

Figure \ref{Figure:PencialAndConic} illustrates the previous definition. 
The following lemma is straightforward from the
construction:

\begin{Lemma}
\label{Lemma:boundaryFuchsian}
Let $\mathbb{H}^2\subset X$ be a Fuchsian plane:
\begin{enumerate}[(a)]
    \item If $\mathbb{H}^2$ is of type I, then 
    $\partial_\infty \mathbb{H}^2$ is a pencil of flags. 
    \item If $\mathbb{H}^2$ is of type II, then 
    $\partial_\infty \mathbb{H}^2$ is the 
    tangent flag to a conic. 
\end{enumerate}
In addition: 
\begin{enumerate}[(a)]
  \setcounter{enumi}{2}
    \item All Fuchsian planes asymptotic to a given pencil of flags are parallel, and contained in a fixed parallel set.
    \item There is a unique 
    Fuchsian plane asymptotic to a given tangent flag to a  conic.    
\end{enumerate}
\end{Lemma}

The following lemma will be used when 
we describe the Fuchsian locus below.

\begin{Lemma}
\label{Lemma:fewplanes}
    Let $0\neq u\in\mathfrak{p}$ with 
    $\mathrm{type}(u)=\mathrm{center}(\sigma_{\mathrm{mod}})$.
    There exist precisely 4 totally geodesic planes 
    $\pi\subset X$
    through $x_0$ with $u\in T_{x_0}\pi$: a flat, a Fuchsian plane of type I and two Fuchsian planes of type II.
\end{Lemma}

\begin{proof}
 Up to multiplication and up to conjugacy, we may assume 
    $$
    u=\begin{pmatrix}
        1 & & \\
         & -1 & \\
         & & 0
    \end{pmatrix}.
    $$
    Then, for any $a,b,c,x,y,z\in\mathbb{R}$ we have  
    \begin{align*}
    \mathrm{ad}_u 
    \begin{pmatrix}
        a & x & y\\
        x & b & z\\
        y & z & c
    \end{pmatrix} &=
    \begin{pmatrix}
        0 & 2x & y\\
        -2x & 0 & -z\\
        -y & z & 0
    \end{pmatrix}
    \\    
     \mathrm{ad}_u\circ \mathrm{ad}_u 
    \begin{pmatrix}
        a & x & y\\
        x & b & z\\
        y & z & c
    \end{pmatrix} &=
    \begin{pmatrix}
        0 & 4x & y\\
        4x & 0 & z\\
        y & z & 0
    \end{pmatrix}
   \\
   \left[
    \begin{pmatrix}
        0 & 0 & y\\
        0 & 0 & z\\
        y & z & 0
    \end{pmatrix}
    ,
    \begin{pmatrix}
        0 & 0 & y\\
        0 & 0 & -z\\
        -y & z & 0
    \end{pmatrix}
   \right]
   &
   =
   2 \begin{pmatrix}
        y^2 & 0 &0\\
        0 & -z^2 & 0\\
        0 & 0 & z^2-y^2
    \end{pmatrix}.
    \end{align*}
    Using \cite{Helgason}, for instance, it
    follows from this computation that the planes $\pi$ are 
    those tangent to the span  $\langle u, v \rangle \subset \mathfrak p$, for
    $$
v=
    \begin{pmatrix}
        -1 & 0 &0\\
        0 & -1 & 0\\
        0 & 0 & 2
    \end{pmatrix},\quad 
    v=
    \begin{pmatrix}
        0 & 1 &0\\
        1 & 0 & 0\\
        0 & 0 & 0
    \end{pmatrix},\quad
    v=
    \begin{pmatrix}
        0 & 0 & 1\\
        0 & 0 & 1\\
        1 & 1 & 0
    \end{pmatrix}, \textrm{ and }
        v=
    \begin{pmatrix}
        0 & 0 & 1\\
        0 & 0 & -1\\
        1 & -1 & 0
    \end{pmatrix},
    $$
This yields the four planes of the statement.
\end{proof}

  We invoke Lemma~\ref{Lemma:boundaryFuchsian} to provide a geometric interpretation of Lemma~\ref{Lemma:fewplanes}. The limit $\lim_{t\to +\infty}\exp(t\,u)$ defines a flag $(p,l)\in \mathbb{P}^2\times (\mathbb{P}^2)^*$ with $p\in l$. Moreover, the involution $i_{x_0}$ induces a polarity on the flag manifold. Consequently, the lemma implies the following geometric characterization in $\mathbb{P}^2$: there exist
\begin{enumerate}[(a)]
    \item a unique triangle containing $p$ and $l$ that is invariant under the polarity $i_{x_0}$;
    \item a unique $i_{x_0}$-invariant pencil of flags containing $p$ and $l$;
    \item exactly two $i_{x_0}$-invariant tangent flags to conics containing $p$ and $l$.
\end{enumerate}

    % Up to conjugation and up to a scalar, we may write 
    % $$
    % u=\begin{pmatrix}
    %     1 & 0 & 0 \\
    %     0 & -1 & 0 \\
    %     0 & 0 & 0
    % \end{pmatrix}
    % $$
    % Finding totally geodesic planes is equivalent to finding $v\in\mathfrak{p}$ orthogonal to $u$ so that 
    % the space $E=\langle u, v\rangle$ satisfies the condition
    % $[[E,E],E]\subseteq E$. An easy calculation shows that the possibilities for $v$ are, up to a scalar.
    % $$
    % \begin{pmatrix}
    %     1 & 0 & 0 \\
    %     0 & 1 & 0 \\
    %     0 & 0 & -2
    % \end{pmatrix}, \quad
    % \begin{pmatrix}
    %     0 & 1 & 0 \\
    %     1 & 0 & 0 \\
    %     0 & 0 & 0
    % \end{pmatrix}, \quad
    % \begin{pmatrix}
    %     0 & 0 & 1 \\
    %     0 & 0 & 1 \\
    %     1 & 1 & 0
    % \end{pmatrix}, \quad
    % \begin{pmatrix}
    %     0 & 0 & 1 \\
    %     0 & 0 & -1 \\
    %     1 & -1 & 0
    % \end{pmatrix}
    % $$
    % They correspond respectively to a maximal flat, a Fuchsian plane of type I, and two Fuchsian planes of type II.
% \end{proof}

\section{The variety of characters}
\label{Section:varietyofCharacters}

We first review the definition of the variety of characters and its main
properties in general, and then we focus on 
the modular group $\Gamma$ and $\mathrm{Isom}(X)$.

\subsection{Definition of the variety of characters}

In this subsection $\Upsilon$ denotes a finitely generated group and $G$ a real reductive algebraic group. 
The notation $\Upsilon$ is not standard, but in this paper 
$\Gamma$ denotes $\mathrm{PSL}_2(\mathbb{Z})$. 
Notice also that we use $\mathcal X$ for the variety of characters, as 
$X$ denotes $\mathrm{SL_3(\mathbb R)}/\mathrm{SO}(3)$.

Let 
$
\hom^{\mathrm{closed}}(\Upsilon, G)
$
denote the subset of representations 
$\rho\in\hom(\Upsilon,G)$ whose $G$-orbit by conjugation 
is closed. 

\begin{Definition} The variety of characters is the 
space of closed orbits:
$$
\mathcal{X}(\Upsilon,G)=\hom^{\mathrm{closed}}(\Upsilon, G)
\big/G.
$$
\end{Definition}

Let $K\subset G$ denote a maximal compact subgroup. 

\begin{Theorem}[Richardson-Slodowy]
\label{Theorem:RS}
    There exists a $K$-invariant algebraic subvariety 
    $\mathcal C \subset \hom(\Upsilon,G)$ such that, for 
    every closed $G$-orbit in $\hom(\Upsilon,G)$, its 
    intersection with $\mathcal C$ consists 
    of exactly one $K$-orbit.

Moreover, each $G$-orbit in $\hom(\Upsilon,G)$ has a 
unique closed $G$-orbit in its Zariski closure; 
equivalently, every $G$-orbit accumulates to a single 
closed $G$-orbit.
\end{Theorem}

This result is proved by Richardson and Slodowy  
\cite{RichardsonSlodowy}, see also  B\"ohm and Lafuente 
\cite{BohmLafuente}.

The subset $\mathcal C\subset \hom(\Upsilon,G)$ in  
Theorem~\ref{Theorem:RS}
is sometimes called the critical set, but we will not discuss its 
definition here. We only use the fact that
$$
\mathcal{X}(\Upsilon,G)\cong\mathcal{C}/K.
$$
This has some consequences:
\begin{enumerate}[(a)]
     \item The space $\mathcal{X}(\Upsilon,G)$ is Hausdorff, because $K$ is compact.
     In addition, by  the second assertion of Theorem~\ref{Theorem:RS}, 
     $\mathcal{X}(\Upsilon,G)$ is the  \emph{Hausdorffication} of the topological quotient
     $\hom(\Upsilon,G)\big/G$.
    \item By a theorem of Procesi and 
    Schwarz~\cite{ProcesiSchwarz} applied to 
    the action of $K$ on
    $\mathcal{C}$, 
    $\mathcal{X}(\Upsilon,G)$ is semi-algebraic.
\end{enumerate}

Following~\cite{JohnsonMillson}, we define:

\begin{Definition}
\label{Def:good}
Let 
%$\rho $ be a representation in 
$\rho\in\hom(\Upsilon,G)$:
\begin{enumerate}[(a)]
    \item The representation $\rho$ is said to be \emph{stable} if its conjugation orbit 
$G\cdot\rho$ is closed and its stabilizer
$\mathrm{Stab}(\rho)=\{\varphi\in G\mid \varphi\cdot\rho=\rho\} $
is finite.
   \item The representation $\rho$ is said to be
    \emph{good} if it is stable and  $\mathrm{Stab}(\rho)$ is trivial.
\end{enumerate}
\end{Definition}

% \begin{Definition}
%     A representation $\rho\in\hom(\Upsilon,G)$ is
%     \begin{enumerate}[(a)]
%             \item     \emph{stable} if  
%     its orbit is closed and its stabilizer is finite,
%         \item     \emph{good} if  
%     its orbit is closed and its stabilizer is trivial,
%     \end{enumerate}
    Denote by $\mathcal{X}^{\mathrm{stable}}(\Upsilon,G)$
    and $\mathcal{X}^{\mathrm{good}}(\Upsilon,G)$ the respective projections in the character variety.   Notice that the subsets   $\mathcal{X}^{\mathrm{good}}(\Upsilon,G)$
  and  $\mathcal{X}^{\mathrm{stable}}(\Upsilon,G)$ are open in 
   $\mathcal{X}(\Upsilon,G)$,
   because the conditions of having finite or trivial stabilizers define $G$-invariant open subsets in  both
$\hom(\Upsilon, G)$
and $\hom^{\mathrm{closed}}(\Upsilon, G)$.
% \end{Definition}

\begin{Proposition}
\label{Proposition:smoothness}
    The subset $\mathcal{X}^{\mathrm{good}}(\Upsilon,G)$ has a natural structure of analytic variety, and 
   $\mathcal{X}^{\mathrm{stable}}(\Upsilon,G)$, of analytic orbifold. 

   If $\rho$ is a smooth point of 
   $ \hom(\Upsilon, G)$, then 
   this structure at the class of $\rho$ is smooth.
\end{Proposition}

\begin{proof} We employ a slice 
$\mathcal{S}\subset \hom(\Upsilon, G)$ provided by \cite[Theorem~3.1]{Porti25}. For a given $\rho\in \hom(\Upsilon, G)$,
the set $\mathcal{S}$ is an analytic subvariety defined in a neighborhood of $\rho$, invariant under the action of the stabilizer
$\mathrm{Stab}(\rho)$, such that 
$\mathcal{S}\big/\mathrm{Stab}(\rho)$ constitutes a neighborhood
of the class of $\rho$ in $\mathcal \mathcal{X}(\Upsilon,G)$.
Moreover, $\mathcal{S}$ is smooth whenever $\hom(\Upsilon, G)$
is smooth at $\rho$.
The naturality of this structure follows from the uniqueness of
$\mathcal{S}$ up to bi-analytic $\mathrm{Stab}(\rho)$-equivariant equivalence, as stated in \cite[Theorem~3.1]{Porti25}.
\end{proof}

\subsection{Fuchsian and elliptic representations}

In this subsection we deal again with the modular group
$\Gamma=\mathrm{PSL}_2(\mathbb{Z})$ and 
the isometry group of the symmetric space
$X=\mathrm{SL}_3(\mathbb{R})/\mathrm{SO}(3)$.

Proposition~\ref{Proposition:NonEllipticNonFuchsianSmooth}
is proved at the end of the subsection.
For convenience, we recall the definitions from the introduction.

% A representation 
% $
% \rho\in \hom\bigl(\Gamma, \mathrm{Isom}(X)\bigr)
% $
% is called 
% \emph{elliptic} if it has a global fixed point in $X$ and
% \emph{Fuchsian} if it preserves a totally geodesic hyperbolic plane.
% 
% 
\begin{Definition}
A representation $\rho\in
\hom_0(\Gamma, \mathrm{Isom}(X))
$ is  called 
\begin{enumerate}[(a)]
\item \emph{elliptic} if it has a global fixed point in $X$,
    \item \emph{Fuchsian} if
its image preserves 
a  totally geodesic hyperbolic plane
$\mathbb{H}^2\subset X$.

It is \emph{Fuchsian of type I} when the invariant $\mathbb H^2$ is congruent to 
$\mathrm{SL}_2(\mathbb R)/\mathrm{SO}(2)$, and \emph{of
type II} when it is congruent to 
$\mathrm{SO}(1,2)/\bigl(\mathrm{SO}(1,2)\cap\mathrm{SO}
(3)\bigr)$.
\end{enumerate}
% In the Fuchsian case it is called \emph{Fuchsian of type I or II} 
% according to whether the invariant plane is of type I or II.
\end{Definition}

Notice that a representation can be simultaneously elliptic 
and Fuchsian (of any type or both).

\begin{Definition}
 A representation in $\mathrm{SL}_n(\mathbb{R})$ is
 \emph{reducible} if its image has a proper invariant subspace in
 $\mathbb R^n$. 
\end{Definition}

Let $\Gamma_2\subset \Gamma$ denote the unique subgroup of index
$2$. Namely 
$$
\Gamma_2=\ker\left(\Gamma\cong \mathbb
{Z}/2\mathbb{Z}* \mathbb{Z}/3\mathbb{Z}\twoheadrightarrow
\mathbb
{Z}/2\mathbb{Z}\right)
$$
and $\Gamma_2\cong \mathbb{Z}/3\mathbb{Z}* \mathbb{Z}/3\mathbb{Z}$

\begin{Lemma}
\label{lemma:FuchsianParallel}
For $\rho 
\in 
\hom(\Gamma, \mathrm{Isom}(X))
$: 
\begin{enumerate}[(a)]
    \item $\rho$ is Fuchsian of type I if and only if 
     $\rho$ preserves a parallel set $P$.
     \item If $\rho$ is Fuchsian of type I, then  $\rho$ restricted to $\Gamma_2$ is reducible. 
     \item If the conjugation orbit of $\rho$ is closed and
      $\rho$ restricted to $\Gamma_2$ is reducible,
      then $\rho$ is Fuchsian of type I.
\end{enumerate}
\end{Lemma}

The proof of (c) is postponed to Lemma~\ref{lemma:FuchsianParallelc}.

\begin{proof}[Proof of (a) and (b)]  
(a)
 Assume that $\rho$ preserves the Fuchsian plane of type I,  $\mathbb H^2\subset X$. 
By Lemma~\ref{Lemma:UniqueParallelSet} (a) the invariant $\mathbb H^2$
 is contained in a single parallel set, that must be preserved. 

Assume now that $\rho$ preserves a parallel set $P$.  
Using the product isometric decomposition 
\(
P\cong \mathbb{R}\times \mathbb{H}^2,
\)
and noting that $\Gamma\cong \mathbb
{Z}/2\mathbb{Z}* \mathbb{Z}/3\mathbb{Z}$,
it follows that the action of $\rho(\Gamma)$ on the $\mathbb{R}$–factor admits a global fixed point $t_0\in \mathbb{R}$. Consequently, 
$\rho$ preserves the Fuchsian plane $\{t_0\}\times \mathbb{H}^2$.

(b) Assume that $\rho$ is Fuchsian of type I, by (a) it preserves a parallel set $P$.
This parallel set $P$ corresponds to a pair $(p, l)$ with 
 $p \in \mathbb{P}^2$,  $l \in (\mathbb{P}^2)^*$ and $p \notin l$.
This pair is $\rho$-invariant and, since 
the restriction $\rho\vert_{\Gamma_2}$ preserves the type, 
$\rho(\Gamma_2)$ fixes both $p$ and $l$. Thus 
$\rho\vert_{\Gamma_2}$ is reducible. 
% (c)$\Rightarrow$(b).
%  Since it is reducible, the restriction 
% $ \rho\vert_{\Gamma_2} $ 
% preserves either a point $p \in \mathbb{P}^2$ or a projective line 
% $l \in (\mathbb{P}^2)^*$. If $ \rho\vert_{\Gamma_2} $ preserves a point 
% $p \in \mathbb{P}^2$, then it also preserves the line $\rho(a)(p) = i_x(p)$, 
% which is in generic position with respect to $p$. By applying the same 
% reasoning in the dual situation, we may assume that $ \rho\vert_{\Gamma_2} $ 
% preserves both a point $p \in \mathbb{P}^2$ and a line
% $l \in (\mathbb{P}^2)^*$, where $p \notin l$ and $i_x$ 
% interchanges $p$ and $l$.
% So the parallel set $P$ determined by the pair $(p, l)$ 
% % ($P$ is the union of all geodesics asymptotic to $p$ and $l$)
% is $\rho$-invariant.
\end{proof}

\begin{Lemma}
\label{Lem:CRirreducible}
 A representation in $\mathrm{SL}_3(\mathbb{R})$
 has a proper invariant subspace in $\mathbb R^3$
 if and only if it 
 has a proper invariant subspace in $\mathbb C^3$.
Namely, it is $\mathbb{R}$-irreducible iff it is $\mathbb{C}$-irreducible.
\end{Lemma}

\begin{proof}
One of the implications is obvious. For the other one, assume that
a representation in $\mathrm{SL}_3(\mathbb{R})$ has a proper invariant 
$\mathbb C$-subspace
 $V\subset \mathbb C^3$.
Let $\overline V$ denote its complex conjugate.
If $V=\overline V$, then $V\cap \mathbb{R}^3 $ is 
a proper invariant subspace
of $\mathbb R^3$.
Otherwise, $V\oplus \overline V$ (when $\dim V=1$)
or $V\cap \overline V$ (when $\dim V=2$) yield a
proper invariant subspace
of $\mathbb R^3$.
\end{proof}

\begin{Remark}
 The very same argument yields that Lemma~\ref{Lem:CRirreducible}
 holds true for
 representations in  $\mathrm{SL}_n(\mathbb{R})$ for $n$ odd.
 It fails for $n$ even.
\end{Remark}

\begin{Proposition}
\label{Proposition:notFEsmooth}
Let $\rho\in\hom(\Gamma,\mathrm{Isom}(X))$ be a representation with closed conjugacy orbit.
\begin{enumerate}[(a)]
    \item If $\rho$ is not Fuchsian of type I then $\rho$
    is stable.
    \item If $\rho$ is not Fuchsian (of any type) nor elliptic, then $\rho$
    is good.
\end{enumerate}
\end{Proposition}

\begin{proof} (a)
   Let $\rho \in \hom(\Gamma, \mathrm{Isom}(X))$ be a representation that is not Fuchsian of type I. A representation into $\mathrm{SL}_3(\mathbb{R})$ that is $\mathbb{C}$-irreducible is stable by \cite{JohnsonMillson}. Therefore, by Lemmas~\ref{lemma:FuchsianParallel} and~\ref{Lem:CRirreducible}, the restriction of $\rho$ to $\Gamma_2$ is stable. It then follows, for example by applying \cite[Theorem~1.1]{JohnsonMillson}, that $\rho$ itself is stable.

(b) Suppose now that $\rho$ is neither Fuchsian nor elliptic. We prove that the stabilizer $\mathrm{Stab}(\rho)$ is trivial. Let $\phi \in \mathrm{Stab}(\rho)$. If $\phi$ preserves orientation, then it is represented by a matrix $A \in \mathrm{SL}_3(\mathbb{R})$ that commutes with $\rho(\gamma)$ for every $\gamma \in \Gamma_2$. Since the restriction of $\rho$ to $\Gamma_2$ is $\mathbb{C}$-irreducible, Schur’s lemma implies that $A$ must be the identity matrix, and consequently $\phi$ is trivial. 

If $\phi$ reverses orientation, then by the previous 
argument $\phi^2$ is trivial. Hence, by 
Lemma~\ref{Lem:4components}, $\phi$ must fix either a 
point $x \in X$ or a totally geodesic subspace 
$\mathbb{H}^2 \subset X$ of type II. This, however, would 
contradict the assumption that $\rho$ is neither elliptic nor Fuchsian.
\end{proof}

\begin{proof}[Proof of 
Proposition~\ref{Proposition:NonEllipticNonFuchsianSmooth}]
Since
 $\mathrm{PSL}_2(\mathbb Z)$ is virtually free,
 by Goldman's theory \cite{Goldman} there are no obstructions 
 to smoothness of the variety of representations.
Then the claim follows from 
Propositions~\ref{Proposition:smoothness} 
and~\ref{Proposition:notFEsmooth}.
\end{proof}

% At the Fuchsian locus we can only say that there is a 
% semi-algebraic structure. 

\subsection{Counting components}

In this subsection we compute the number of components of $\mathcal{X}(\Gamma,\mathrm{Isom}(X))$.

Recall  the presentation in the introduction of the modular group
\begin{equation}
\label{eqn:presentation}
\Gamma =\mathrm{PSL}_2(\mathbb Z)\cong \langle a,b\mid a^2=b^3=1\rangle  .
\end{equation}
In Lemma~\ref{Lem:4components} we have shown that there are 4 connected components of the subset of involutions of 
$\mathrm{Isom}(X)$, which gives the possibilities for the image of $a$. For the image of $b$, we recall the notation:
\begin{equation}
 \label{eqn:rotation}
R_\frac{2\pi}{3} = \begin{pmatrix}
          1 & 0 & 0 \\
          0 & -\frac{1}{2} & -\frac{\sqrt{3}}{2} \\
          0 & \frac{\sqrt{3}}{2} & -\frac{1}{2}
         \end{pmatrix}.
\end{equation}
The 
fixed point set $\mathrm{Fix}(R_\frac{2\pi}{3})$ is a \emph{singular}
geodesic $\exp( t\, p_0)$,
where $p_0=\mathrm{diag}(2,-1,-1)$,  
see Equation~\eqref{eqn:p0}.

The following lemma is straightforward: 

\begin{Lemma}
    \label{Lemma:orderthree}
    A nontrivial isometry of $X$ of order three is represented, up to conjugacy, by a matrix
    conjugate to $R_{\frac{2\pi}{3}}$.
\end{Lemma}

Thus it follows that the subset of elements 
$B \in \mathrm{Isom}(X)$
that 
satisfy $B^3=\mathrm{Id}$ has two connected components, 
one of 
them trivial. As 
$\Gamma\cong \mathbb{Z}_2*\mathbb{Z}/3\mathbb{Z}$, by combining
Lemmas~\ref{Lem:4components} (the space of involutions has 4 components
that are also conjugacy classes) and Lemma~\ref{Lemma:orderthree}, we deduce: 

\begin{Corollary}
    The space of representations $\hom(\Gamma, \mathrm{Isom}(X) )$ has 8 components. Five components  map either $a$ or $b$
    to the identity and  those components project to single points
    in the variety of characters.
\end{Corollary}

We are interested in the 3 nontrivial components. 
As in the introduction, 
for $i = 0,1,2$, define
\[
\hom_i\bigl(\Gamma, \mathrm{Isom}(X)\bigr)
\]
to be the subset consisting of all $\rho \in \hom\bigl(\Gamma, \mathrm{Isom}(X)\bigr)$ such that $\rho(b)$ is nontrivial and
\[
\rho(a) \text{ is }
\begin{cases}
    \text{a rotation of order } 2 & \text{for } i = 0,\\[2pt]
    \text{an inversion} & \text{for } i = 1,\\[2pt]
    \text{a reflection with respect to }\mathbb{H}^2 & \text{for } i = 2.
\end{cases}
\]
The corresponding components in the variety of characters 
% real GIT quotients (in the sense of \cite{RichardsonSlodowy}) with respect to the conjugacy action of $\mathrm{Isom}(X)$ 
are denoted by
\[
\mathcal{X}_i\bigl(\Gamma, \mathrm{Isom}(X)\bigr)
=
\hom_i^{\textrm{closed}}\bigl(\Gamma, \mathrm{Isom}(X)\bigr)
\big/\mathrm{Isom}(X),
\]
for $i=0,1,2$.

\subsection{Closed orbits}

In this section we prove  Proposition~\ref{Proposition:mindis} 
in the introduction, which we split into two statements,
Proposition~\ref{Proposition:PerpendicularClosed} and 
Corollary~\ref{Corollary:ClosedPerpendicular}.
We start with a couple of geometric lemmas.

\begin{Lemma}
\label{Lemma:NoAsymptotic}
    For $\rho\in \hom_2(\Gamma,\mathrm{Isom}(X)) $,  $\mathrm{Fix}(\rho(a))$ and
    $\mathrm{Fix}(\rho(b))$ are not asymptotic.
    In particular either $\mathrm{Fix}(\rho(a))$ and
    $\mathrm{Fix}(\rho(b))$ intersect, or they have a common perpendicular.
\end{Lemma}

\begin{proof}
    By Lemma~\ref{Lem:FuchsianRegular} every point in 
    $\partial_\infty \mathrm{Fix}(\rho(a))$
    has type the center of the model spherical Weyl chamber $\sigma_{\mathrm{mod}}$.
    As $\mathrm{Fix}(\rho(b))$ is a singular geodesic, its ideal end-points
    $\partial_\infty\mathrm{Fix}(\rho(b))$
    are singular, hence they have type in the boundary
    $\partial \sigma_{\mathrm{mod}}$.
    So $\partial_\infty \mathrm{Fix}(\rho(a))
    \cap 
    \partial_\infty\mathrm{Fix}(\rho(b))=\emptyset$, which means 
    that they are not asymptotic. 
    The second assertion follows from convexity of the distance function.
\end{proof}

\begin{Definition}
A \emph{transvection} along a geodesic $\gamma$ is defined as the composition of two inversions  whose centers lie in 
 $\gamma$.
The resulting isometry has translational displacement equal to twice the distance between these centers.    
\end{Definition}

\begin{Lemma}
\label{Lemma:asymptotic}
    Let $\gamma_1$ and $\gamma_2$ be geodesic lines in $X$.
    \begin{enumerate}[(a)]
        \item If $\gamma_1$ and $\gamma_2$ are asymptotic at both ends, then they are parallel.
        \item If $\gamma_1$ and $\gamma_2$ are asymptotic at a 
        single end, and if $\tau_t$ denotes the transvection along $\gamma_1$ of displacement $t>0$ (in the direction of the 
        common ideal end), then $\tau_{-t}(\gamma_2)$ converges to a line parallel to $\gamma_1$ as $t\to +\infty$.
    \end{enumerate}
\end{Lemma}

Here convergence stands for pointed Hausdorff convergence: 
$\tau_{-t}(\gamma_2)\to \gamma$
as $t\to+\infty$
means that for every compact 
subset
$K\subset X$, $\tau_{-t}(\gamma_2)\cap K$ 
Hausdorff converges to $\gamma\cap K$.

\begin{proof}
   Assertion (a) follows directly from the convexity of the distance function: any convex and bounded function on $\mathbb{R}$ must be constant.
    
For Assertion (b), observe that
$$d(\gamma_1(s),\tau_{-t}(\gamma_2))=d(\gamma_1(s+t),\gamma_2).$$
On the other hand, the function $\varphi(s)=d(\gamma_1(s),\gamma_2)$ is convex and bounded at $+\infty$,
 which implies that both $\varphi$ and its derivative $\varphi'$ are non-increasing. Consequently, $\lim_{s\to+\infty}\varphi'(s)=0$, and $\varphi$ approaches a constant function as $s\to+\infty$. Convergence  then follows easily.
\end{proof}

% \begin{Corollary}
% \label{Corollary:closed_orbits}
%     For $i=1,2$ all orbits in 
%     $\hom_i\bigl(\Gamma, \mathrm{Isom}(X)\bigr)$ are closed.    
% \end{Corollary}

% \begin{proof}
% The case $i=1$, as $\mathrm{Fix}(\rho(a))$ is a single point, 
% this is a straightforward consequence of 
% Proposition~\ref{Proposition:PerpendicularClosed}.
% When $i=2$, this is now a corollary of
% Proposition~\ref{Proposition:PerpendicularClosed}
% together wilt Lemma~\ref{Lemma:NoAsymptotic}.
% \end{proof}

% In the component $\hom_0(\Gamma,\mathrm{Isom}(X))$,
% $\mathrm{Fix}(\rho(a))$ is a parallel set
% (and $\mathrm{Fix}(\rho(b))$ is a singular line 
% for all components considered here). 
% They do not have always a common perpendicular, in particular 
% as 
% $\mathrm{Fix}(\rho(a))$ and
% $\mathrm{Fix}(\rho(b))$
% may be asymptotic. So we need a lemma on asymptotic lines.

% \subsection{Asymptotic lines and non closed orbits}
\label{Subsection:Asymptoticlines}
% In this subsection we show that when one end
% of the line $\mathrm{Fix}(\rho(b))$ is asymptotic to 
% the parallel set 
% $\mathrm{Fix}(\rho(a))$ 
% then the orbit of $\rho$ is not closed. 
% We start with a lemma on asymptotic lines.

\begin{Proposition}
\label{Proposition:Accumulate}
For a representation $\rho \in \hom_0(\Gamma,\mathrm{Isom}(X))$, assume that precisely one end of $\mathrm{Fix}(\rho(b))$ is asymptotic to the parallel set $\mathrm{Fix}(\rho(a))$. Then the $\mathrm{Isom}(X)$–orbit of $\rho$ in $\hom_0(\Gamma,\mathrm{Isom}(X))$ accumulates at a representation $\rho_{\infty} \in \hom_0(\Gamma,\mathrm{Isom}(X))$ such that either
\(
\mathrm{Fix}(\rho_{\infty}(b)) \subset \mathrm{Fix}(\rho_{\infty}(a))
\)
or $\mathrm{Fix}(\rho_{\infty}(b))$ is parallel to a 
geodesic line contained in $\mathrm{Fix}(\rho_{\infty}
(a))$. 

In particular $\rho_{\infty}$ does not lie in 
the orbit of $\rho$.
\end{Proposition}

\begin{proof}
    The fact that precisely one end of $\mathrm{Fix}(\rho(b))$ is asymptotic to $\mathrm{Fix}(\rho(a))$ means that there is a geodesic $\gamma$ in the parallel set 
    $\mathrm{Fix}(\rho(a))$ so that $\mathrm{Fix}(\rho(b))$ and 
    $\gamma$ are asymptotic at precisely one end-point. 
    By Lemma~\ref{Lemma:asymptotic} there is a family of transvections along $\gamma$ such that 
    $\tau_{-t}(\mathrm{Fix}(\rho(b)))$ converges to a geodesic 
    parallel to $\gamma$. As $\gamma\subset \mathrm{Fix}(\rho(a))$, $\tau_t$ preserves the parallel set
    $\mathrm{Fix}(\rho(a))$. 
    Now consider the family of conjugate representations 
    $$\rho_t=\tau_{-t}\,\rho\,\tau_t.$$ 
    Then $\mathrm{Fix}(\rho_t(a))=\mathrm{Fix}(\rho(a))$ remains constant  and 
    $\mathrm{Fix}(\rho_t(b))=\tau_{-t}(\mathrm{Fix}(\rho(b)))$
    converges to a geodesic parallel to $\gamma$. This yields that  $\rho_t$ converges to $\rho_{\infty}\in \hom_0(\Gamma,\mathrm{Isom}(X))$, with $\mathrm{Fix}(\rho_\infty(a))=\mathrm{Fix}(\rho(a))$ and 
    $\mathrm{Fix}(\rho_\infty(b))$
     parallel to $\gamma\subset \mathrm{Fix}(\rho(a))$
     (possibly equal).  
\end{proof}

\begin{Corollary}
\label{Corollary:ClosedPerpendicular}
    If the orbit of $\rho$ is closed, then there exists a 
    point $x\in \mathrm{Fix}(\rho(a))$ that minimizes in 
    $\mathrm{Fix}(\rho(a))$
    the distance to $\mathrm{Fix}(\rho(b))$.
\end{Corollary}

\begin{proof}
  Such a point 
    $x\in \mathrm{Fix}(\rho(a))$ always exists if 
    $\rho\in\hom_1(\Gamma,\mathrm{Isom}(X))$, by
    construction, or if $\rho\in\hom_2(\Gamma,\mathrm{Isom}(X))$ by Lemma~\ref{Lemma:NoAsymptotic}. 
    Therefore, we may assume that 
    $\rho\in\hom_0(\Gamma,\mathrm{Isom}(X))$.
    In this case, if the orbit of $\rho$ is closed, then  Proposition~\ref{Proposition:Accumulate}
    implies that either both ends of  
      $\mathrm{Fix}(\rho(b))$ are 
      asymptotic to
      $\mathrm{Fix}(\rho(a))$ 
      (in which case the fixed point sets are parallel),
      or else no end of $\mathrm{Fix}(\rho(b))$ is
      asymptotic to
      $\mathrm{Fix}(\rho(a))$ (in which case the fixed point sets admit a common perpendicular). In both cases there exists
      such an $x\in \mathrm{Fix}(\rho(a))$ that minimizes 
      the distance to $\mathrm{Fix}(\rho(b))$.
\end{proof}

\begin{Proposition}
\label{Proposition:PerpendicularClosed}
If there exists a point \(x\in \mathrm{Fix}(\rho(a))\) that realizes the distance between \(\mathrm{Fix}(\rho(a))\) and \(\mathrm{Fix}(\rho(b))\), then the conjugacy orbit of \(\rho\in \hom(\Gamma,\mathrm{Isom}(X))\) is closed.
\end{Proposition}

\begin{proof}
Let \((\rho_n)\) be a sequence of representations, each conjugate to \(\rho\), and suppose that \(\rho_n \to \rho_\infty\). Assume that for each \(n\) there exists a point \(x_n\in \mathrm{Fix}(\rho_n(a))\) that realizes the distance from \(\mathrm{Fix}(\rho_n(a))\) to \(\mathrm{Fix}(\rho_n(b))\). Our objective is to show that \(\rho_\infty\) is conjugate to each \(\rho_n\), which will imply that the conjugacy orbit of \(\rho\) is closed.

First, note that the convergence 
\(\rho_n \to \rho_\infty\) entails 
that the fixed-point sets 
\(\mathrm{Fix}(\rho_n(a))\) and \(\mathrm{Fix}
(\rho_n(b))\) converge (in the sense of Hausdorff 
convergence on bounded subsets) to \(\mathrm{Fix}
(\rho_\infty(a))\) and \(\mathrm{Fix}
(\rho_\infty(b))\), respectively.

For each \(n\), the point 
\(x_n\in \mathrm{Fix}(\rho_n(a))\) 
minimizing the distance to \(\mathrm{Fix}
(\rho_n(b))\) may or may not be unique. Suppose first 
that the minimizer is not unique. Then, by 
Lemma~\ref{Lemma:NoAsymptotic}, 
$\rho_n\in\hom_0(\Gamma,\mathrm{Isom}(X))$ and  
\(\mathrm{Fix}(\rho_n(a))\) is a parallel set. 
Furthermore, 
 \(\mathrm{Fix}(\rho_n(b))\) is a line parallel to 
 \(\mathrm{Fix}(\rho_n(a))\), 
 not necessarily contained in it, by convexity of the distance function. Passing to the limit, we deduce that \(\mathrm{Fix}(\rho_\infty(b))\) is parallel to \(\mathrm{Fix}(\rho_\infty(a))\), and from 
 the description of fixed point sets, one readily concludes that \(\rho_\infty\) is conjugate to each \(\rho_n\).

Hence, for the remainder of the proof we assume that for each \(n\) the minimizer \(x_n\) is unique. After passing to a subsequence, we may suppose that either \((x_n)\) converges in \(X\) or escapes to infinity.

Assume first that \(x_n \to x_\infty\in X\). By conjugating each \(\rho_n\) by an isometry contained in a fixed compact subset of \(\mathrm{Isom}(X)\), we may arrange that \(x_n = x_\infty\) for all \(n\). In this situation, the conjugating isometries lie in the stabilizer \(\mathrm{Stab}(x_\infty)\), which is compact; therefore, after extracting a subsequence, these conjugating isometries converge. It follows that \(\rho_\infty\) is conjugate to each \(\rho_n\).

Now consider the case in which \((x_n)\) diverges. Since \(x_n\) minimizes the distance from \(\mathrm{Fix}(\rho_n(a))\) to \(\mathrm{Fix}(\rho_n(b))\), 
these fixed point sets have a common perpendicular with endpoints \(x_n\in \mathrm{Fix}(\rho_n(a))\) and \(y_n\in \mathrm{Fix}(\rho_n(b))\). In particular, \(y_n\) realizes the distance from \(\mathrm{Fix}(\rho_n(b))\) to \(\mathrm{Fix}(\rho_n(a))\). Since the distance function to  \(\mathrm{Fix}(\rho_n(a))\) is convex along \(\mathrm{Fix}(\rho_n(b))\) and attains its minimum at \(y_n\), this distance function is either (i) constant along \(\mathrm{Fix}(\rho_n(b))\) or (ii) unbounded in both directions along \(\mathrm{Fix}(\rho_n(b))\).
In case (i), \(\mathrm{Fix}(\rho_n(b))\) remains at a constant distance from a line contained in \(\mathrm{Fix}(\rho_n(a))\). This contradicts the assumption that \(x_n\) is the unique minimizer.

Next we prove that case (ii) is also impossible. Parameterize the line \(\mathrm{Fix}(\rho_n(b))\) by arc length as
$$ \gamma_n:\mathbb R\to \mathrm{Fix}(\rho_n(b)), 
\qquad 
\textrm{ with }
\gamma_n(0)=y_n. $$
Since all the representations \(\rho_n\) are conjugate, we may choose these parameterizations so that the function
$$ f(t)=d\bigl(\gamma_n(t),\mathrm{Fix}(\rho_n(a))\bigr) $$
is independent of \(n\). In case (ii),
$ \lim_{t\to\pm\infty}f(t)=+\infty$. 
Since \((x_n)\) diverges and \(d(x_n,y_n)\) is constant, it follows that \((y_n)\) diverges as well. Therefore, for every compact subset \(K\subset X\), the distance from
$ \mathrm{Fix}(\rho_n(b))\cap K$ 
to \(\mathrm{Fix}(\rho_n(a))\) tends to infinity,
because $\lim_{t\to\pm\infty} f(t)=+\infty$ and $(y_n)$ diverges.
This is incompatible with the convergence, on bounded subsets, of both \(\mathrm{Fix}(\rho_n(a))\) and \(\mathrm{Fix}(\rho_n(b))\). The proof is therefore complete.
\end{proof}

Proposition~\ref{Proposition:mindis} follows from 
Corollary~\ref{Corollary:ClosedPerpendicular}
and 
Proposition~\ref{Proposition:PerpendicularClosed}.
We also state a couple of corollaries:

\begin{Corollary}
\label{Corollary:hom12closed}
  The orbit of any representation in 
   $\hom_1(\Gamma,\mathrm{Isom}(X))$ or
   $\hom_2(\Gamma,\mathrm{Isom}(X))$
    is closed.
\end{Corollary}

\begin{proof}
For $\rho\in \hom_1(\Gamma,\mathrm{Isom}(X))$,
$\mathrm{Fix}(\rho(a))$ is a single point, and for
$\rho\in \hom_2(\Gamma,\mathrm{Isom}(X))$, use 
Lemma~\ref{Lemma:NoAsymptotic} to find points that realize the distance between both fixed point sets.
\end{proof}

\begin{Corollary}
     \label{Corollary:Fix}
    Let $\rho \in \hom_0(\Gamma,\mathrm{Isom}(X))$. Then the orbit of $\rho$ is closed if and only if  either none or both ideal end-points of $\mathrm{Fix}(\rho(b))$ are asymptotic to $\mathrm{Fix}(\rho(a))$.
\end{Corollary}

\begin{proof}
The condition that either none of both ideal-end points
    of $\mathrm{Fix}(\rho(b))$ are asymptotic to $\mathrm{Fix}(\rho(a))$ is equivalent to the fact that the distance between both fixed point sets is realized.
\end{proof}

Next 
we can prove Lemma~\ref{lemma:FuchsianParallel}(c), that we restate.

\begin{Lemma}
\label{lemma:FuchsianParallelc}
If  $\rho 
\in 
\hom(\Gamma, \mathrm{Isom}(X))
$ has a closed orbit and
      $\rho$ restricted to $\Gamma_2$ is reducible,
      then $\rho$ is Fuchsian of type I.
\end{Lemma}

\begin{proof}
   If $\rho(a)$ or $\rho(b)$ is trivial, then the statement follows immediately. We therefore assume that both
$\rho(a)$ and $\rho(b)$ are nontrivial.

The hypothesis that the restriction of $\rho$ to $\Gamma_2$
 is reducible implies that $\rho(\Gamma_2)$ admits a fixed point
 in $\mathbb{P}^2$ or in its dual $(\mathbb{P}^2)^*$ (or in 
both). This fixed point is unique, since $\rho(b)$ has a unique
 fixed point in both $\mathbb{P}^2$ and $(\mathbb{P}^2)^*$.
 We now distinguish cases according to whether 
$\rho(a)$ reverses or preserves the orientation.

First, suppose
 that $\rho(a)$ is orientation-reversing, namely 
    $\rho\in\hom_i(\Gamma,\mathrm{Isom}(X))$ for $i=1$ or $2$.
In this case, $\rho(a)$ interchanges
 $\mathbb{P}^2$ and $(\mathbb{P}^2)^*$, which implies that 
$\rho(\Gamma_2)$
 has a fixed point in both $\mathbb{P}^2$ and $(\mathbb{P}^2)^*$.
 Consequently, the parallel set 
$P(\mathrm{Fix}(\rho(b))) \cong \mathbb{R}\times \mathbb{H}^2$
(consisting of parallels to $\mathrm{Fix}(\rho(b))$)
is $\rho$-invariant.
    By Lemma~\ref{lemma:FuchsianParallel} (a) it
 follows that $\rho(a)$ has an invariant Fuchsian plane of type~I.

Now assume that $\rho(a)$ is orientation-preserving, 
    namely $\rho\in\hom_0(\Gamma,\mathrm{Isom}(X))$.
Because
$\rho(a)$ preserves the type, it also preserves the (unique) fixed point
 of $\rho(\Gamma_2)$ in $\mathbb{P}^2$ 
or in $(\mathbb{P}^2)^*$. Hence,
$\mathrm{Fix}(\rho(a))$ and
 $\mathrm{Fix}(\rho(b))$ are asymptotic. By 
Corollary~\ref{Corollary:ClosedPerpendicular}, 
as the orbit of $\rho$ is closed,
the fixed point set $\mathrm{Fix}(\rho(b))$
 is either parallel to, or contained in $\mathrm{Fix}(\rho(a))$. Therefore, the parallel set
\( P(\mathrm{Fix}(\rho(b))) \cong \mathbb{R}\times
 \mathbb{H}^2 \)
is invariant under $\rho$, and we apply  
Lemma~\ref{lemma:FuchsianParallel} (a) again.
\end{proof}

% We start analyzing the component 
% $\mathcal X_1\bigl(\Gamma, \mathrm{Isom}(X)\bigr)
% $
% in the following section, as it is simpler than 
% $\mathcal X_0\bigl(\Gamma, \mathrm{Isom}(X)\bigr)
% $.

%
\section{The component $\mathcal{X}_1(\Gamma, \mathrm{Isom}(X))$}
\label{Section:X1}

In this section we describe the topology of 
$\mathcal{X}_1(\Gamma, \mathrm{Isom}(X))$
and its Fuchsian and elliptic locus. Recall that, 
for $\rho\in \hom_1(\Gamma,\mathrm{Isom}(X))$,
$\rho(a)$ is an inversion and $\rho(b)$ 
a nontrivial rotation.

To determine the topology of $\mathcal{X}_1(\Gamma, \mathrm{Isom}(X))$  consider the composition
$$
X\overset \varphi\longrightarrow \hom_1(\Gamma, \mathrm{Isom}(X)) 
\overset {\mathrm{pr}}\longrightarrow \mathcal{X}_1(\Gamma, \mathrm{Isom}(X)),
$$
where, for $x\in X$, $\varphi(x)$ is the representation 
\begin{align*}
  \varphi(x)\colon \Gamma & \to \mathrm{Isom}(X)\\ 
    a & \mapsto i_x\\
    b & \mapsto R_{\frac{2\pi}{3}}
\end{align*}
where $i_x$ denotes the inversion with respect to $x$,  $R_{\frac{2\pi}{3}}$ is a fixed 
rotation as in~\eqref{eqn:rotation},
and $\mathrm{pr}$ denotes the natural projection to the 
character variety (all conjugacy orbits in this 
component
are closed,
by
Corollary~\ref{Corollary:hom12closed}).
Consider also the centralizer of~$R_{\frac{2\pi}{3}}$
$$
Z(R_{\frac{2\pi}{3}})=\{
g\in\mathrm{Isom}(X)\mid g R_{\frac{2\pi}{3}}=R_{\frac{2\pi}{3}} g
\}.
$$

\begin{Lemma}
\label{Lemma:homeoZ}
    The composition $\mathrm{pr}\circ\varphi$ factors to a homeomorphism
    $$
    Z(R_{\frac{2\pi}{3}})\big\backslash X\cong \mathcal{X}_1(\Gamma, \mathrm{Isom}(X)).
    $$
\end{Lemma}

\begin{proof}
By construction, $\mathrm{pr}\circ\varphi$ descends to a continuous bijection
$$
Z(R_{\frac{2\pi}{3}})\backslash X\to \mathcal{X}_1(\Gamma,\mathrm{Isom}(X)).
$$
Openness follows by elementary arguments. \end{proof}

\begin{Proposition}
\label{Prop:homeocharacter}
 We have a natural homeomorphism
 $$
 \mathcal{X}_1( \Gamma, \mathrm{Isom}(X) )
 % \big/ \mathrm{Isom}(X)
 \cong
 (S^1\times\mathbb{Z}_2)\big\backslash \big( \mathbb H^2\times \mathbb H^2 \big),
 $$
 where an element $e^{i\theta}\in S^1$
 acts as a rotation of angle $\theta$ on one of the factors $\mathbb H^2$,
 and angle $2\theta$ on the other factor $\mathbb H^2$,
 whilst $\mathbb Z_2$ acts as a rotation of angle $\pi$ 
 on both factors.
\end{Proposition}

\begin{proof}
By  Lemma~\ref{Lemma:homeoZ}
it suffices to determine the quotient
$Z(R_{\frac{2\pi}{3}}) \big\backslash X$, by using that 
$Z(R_{\frac{2\pi}{3}})\cong \mathbb R\times S^1\times \mathbb{Z}_2$.
To describe this quotient we consider the fibration induced by the
nearest point projection to the parallel set
$$
\mathbb H^2\to X\overset\pi\to P(\gamma)\cong
\mathbb R\times\mathbb H^2,
$$
where $\gamma=\mathrm{Fix}(R_{\frac{2\pi}{3}})$.
The fiber is a Fuchsian plane of type II,
Lemma~\ref{Lem:Fiber}. 
The factor $\mathbb R$ in the product
$Z(R_{\frac{2\pi}{3}})\cong \mathbb R\times S^1\times \mathbb{Z}_2$ acts by translation on
the first factor of the product
$P(\gamma)\cong \mathbb R\times\mathbb H^2$. So, by taking a 
$Z(R_{\frac{2\pi}{3}})$-invariant trivialization of the bundle, we have
$$
  \mathcal{X}_1( \Gamma, \mathrm{Isom}(X) )
 % \big/ \mathrm{Isom}(X)
 \cong
(S^1\times\mathbb Z_2)\big\backslash \big( \mathbb H^2\times \mathbb H^2 \big).
$$
To determine the action of $S^1$ on each factor $\mathbb H^2$,
it suffices to consider the adjoint action of the one-parameter group
$$
R_\theta=
\begin{pmatrix}
 1 & 0 & 0 \\
 0 &  \cos\theta & -\sin \theta \\
 0 & \sin\theta & \cos\theta
\end{pmatrix}
\qquad \theta\in  [0,2\pi ]
$$
on the tangent space to:
\begin{enumerate}[(a)]
    \item the factor $\mathbb H^2$ of the parallel set, namely $\langle p_1,p_2\rangle$, and
    \item the factor $\mathbb H^2$ of the fiber,
which equals  $\langle f_1,f_2\rangle$. 
\end{enumerate}
% on the
% tangent space of the factor $\mathbb H^2$ of the parallel set, which is the span $\langle p_1,p_2\rangle$,
% % $$
% % \begin{pmatrix}
% %  0 & 0 & 0 \\
% %  0 &  1 & 0 \\
% %  0 & 0 & -1
% % \end{pmatrix}
% % \quad\textrm{ and }\quad
% % \begin{pmatrix}
% %  0 & 0 & 0 \\
% %  0 &  0 & 1 \\
% %  0 & 1 & 0
% % \end{pmatrix},
% % $$
% and on the tangent space to the fiber $\mathbb H^2$,
% which equals the span $\langle f_1,f_2\rangle$. 
Here $p_1$ and $p_2$ are as in~\eqref{eqn:p} and
$f_1$ and $f_2$, as in~\eqref{eqn:f},
and the action of $R_\theta$
is the adjoint action.
% $$
% \begin{pmatrix}
%  0 & 1 & 0 \\
%  1 & 0 & 0 \\
%  0 & 0 & 0
% \end{pmatrix}
% \quad\textrm{ and }\quad
% \begin{pmatrix}
%  0 & 0 & 1 \\
%  0 & 0 & 0 \\
%  1 & 0 & 0
% \end{pmatrix}.
% $$
The matrix $R_\theta$ acts by a rotation of angle $2\theta$
on $\langle p_1,p_2\rangle$, 
and of angle $\theta$ on   
 $\langle f_1,f_2\rangle$.
 The generator of $\mathbb Z_2$ is an inversion in a point
 of $\gamma$ and acts by a rotation of angle $\pi$  on each factor.
This proves the proposition.
\end{proof}

\begin{Corollary}
 We have a homeomorphism
$$ 
% \hom_1( \Gamma, \mathrm{Isom}(X) )
%  \big/ \mathrm{Isom}(X)\cong 
     \mathcal{X}_1(\Gamma, \mathrm{Isom}(X) )
     \cong{\mathbb{R}^3}.
$$
\end{Corollary}

\begin{proof}
In Proposition~\ref{Prop:homeocharacter} we pass to the orbit space associated with the action of $S^1\times\mathbb{Z}_2$. Identifying each copy of $\mathbb{H}^2$ with the complex plane $\mathbb{C}$, it suffices to analyse the quotient
\[
(S^1\times\mathbb{Z}_2)\backslash \mathbb{C}^2 .
\]
Since the action preserves the standard norm on $\mathbb{C}^2$, the resulting orbit space is naturally the metric cone over the quotient of the unit sphere $S^3\subset \mathbb{C}^2$. We first take the quotient by the $S^1$-subgroup. Consider the map
\[
\begin{array}{rcl}
     S^3&\to & \mathbb{CP}^1  \\
     (z_1,z_2)& \mapsto & [z_1^2: z_2] ,
\end{array}
\]
which is constant on $S^1$-orbits. A direct verification shows that it induces a homeomorphism
\[
S^1\backslash S^3 \cong \mathbb{CP}^1 \cong S^2 .
\]
Moreover, the residual $\mathbb{Z}_2$-action on $\mathbb{CP}^1$ is given by
\[
[w_1,w_2]\longmapsto [w_1,-w_2],
\]
which corresponds, under the identification $\mathbb{CP}^1\cong S^2$, to a rotation of order two. Consequently, the further quotient by $\mathbb{Z}_2$ is again homeomorphic to $S^2$.
\end{proof}

\begin{Lemma}
\label{lemma:Fuchsian1}
For $\rho 
\in 
\hom_1(\mathrm{PSL}_2(\mathbb Z), \mathrm{Isom}(X))
$:
\begin{enumerate}[(a)]
    \item $\rho$ is Fuchsian of type I if and only if $x\in P(\gamma)$.
    \item  $\rho$ is Fuchsian of type II if and only if $\pi(x)\in \gamma$, where $\pi\colon X\to P(\gamma)$ is the nearest point projection.  
\end{enumerate}
\end{Lemma}

\begin{proof}

    To prove (a), assume first that  $x\in P(\gamma)$.
    As $P(\gamma)\cong\mathbb R\times \mathbb{H}^2 $, $x$ is contained in one of the factors $\{t\}\times\mathbb{H}^2$: this is a Fuchsian plane of type I preserved by both $\rho(a)$ and $\rho(b)$.

    For the other implication, if $\rho$ is Fuchsian of type I it preserves a parallel set $P$ by Lemma~\ref{lemma:FuchsianParallel}, and this parallel set is precisely $P(\gamma)$, as it is the unique 
    parallel set preserved by $\rho(b)$.
    As $a$ has finite order and $P(\gamma)$ is totally geodesic,   
    $\mathrm{Fix}(\rho(a))=\{x\}$ intersects $ P(\gamma)$.
    
    To prove (b) one implication is easy, because the fibre of the nearest point projection at any point of $\gamma$ is a Fuchsian plane of type II invariant by $\rho(b)$.

    For the other implication, let $\mathbb{H}^2_{\mathrm{II}}$ denote the $\rho$-invariant plane of type II. As $a$ has finite order
    $x=\mathrm{Fix}(\rho(a))\in 
    \mathbb{H}^2_{\mathrm{II}}$ 
    and $\mathrm{Fix}(\rho(b))$
    intersects perpendicularly  
    $\mathbb{H}^2_{\mathrm{II}}$. 
    By looking at the action $\rho(b)$ 
    at the tangent space of any point in $\gamma$, 
    and in particular at the $\rho(b)$-invariant subspaces,
    it follows that $\mathbb{H}^2_{\mathrm{II}}$ is 
    perpendicular to 
    $P(\gamma)$.
\end{proof}

Let $\mathcal{E}_1$ denote the intersection 
$\mathcal{E}\cap \mathcal{X}_1(\Gamma,\mathrm{Isom}(X))$, and 
similarly for
$\mathcal{F}_{\mathrm{I},1}$
and
$\mathcal{F}_{\mathrm{II},1}$.

\begin{Corollary}
    In the homeomorphism $\mathcal{X}_1(\Gamma,\mathrm{Isom}(X))\cong\mathbb R^3$, $\mathcal{F}_{\mathrm{I},1}$ corresponds 
    to $[0,\infty)\times \{O\}$ and $\mathcal{F}_{\mathrm{II},1}$,  
    to $(-\infty,0]\times \{O\}$, where $O$ denotes the origin in $\mathbb R^2$. In addition
    $\mathcal{E}_1= \mathcal{F}_{\mathrm{I},1}\cap \mathcal{F}_{\mathrm{II},1
    }$ corresponds to the origin in $\mathbb R^3$.
\end{Corollary}

\begin{proof}
We use the homeomorphism of Proposition~\ref{Prop:homeocharacter}
$$
 \mathcal{X}_1( \Gamma, \mathrm{Isom}(X) )
 \cong
 (S^1\times\mathbb{Z}_2)\big\backslash \big( \mathbb H^2\times \mathbb H^2 \big).
 $$
By Lemma~\ref{lemma:Fuchsian1}
Fuchsian characters of type II correspond to  
$\mathbb H^2\times \{p_2\}$, and those of type I to the second 
$\{p_1\}\times \mathbb H^2$-factor, where $(p_1,p_2)\in
\mathbb{H}^2\times \mathbb{H}^2$ is the fixed point by the action of 
$S^1\times\mathbb{Z}_2$.
The assertion on the Fuchsian loci follows easily.

A representation $\rho $ has a global  fixed point in $X$ only when 
$\mathrm{Fix}(\rho(a))\in \mathrm{Fix}(\rho(b))$. Thus via the 
previous homeomorphism the elliptic locus
corresponds precisely to the fixed point 
$(p_1,p_2)\in
\mathbb{H}^2\times \mathbb{H}^2$.
\end{proof}

\begin{Remark}
    Every Fuchsian representation in 
    $\hom_1(\mathrm{PSL}_2(\mathbb Z), 
    \mathrm{Isom}(X))$
    preserves the orientation of the invariant Fuchsian 
    plane, because an inversion of $\mathbb H^2$ is 
    orientation-preserving. We shall see that this is 
    not
    the case for every component of the variety 
    of representations.
\end{Remark}

\section{The component $\mathcal{X}_0(\Gamma,\mathrm{Isom}(X))$}
\label{Section:X0}

This section is organized into two subsections. In the first, we describe the topology of the component $\mathcal{X}_0(\Gamma,\mathrm{Isom}(X))$, while in the second we analyze its elliptic and Fuchsian loci. Recall that for $\rho \in \hom_0(\Gamma,\mathrm{Isom}(X))$, the element $\rho(a)$ is an orientation-preserving isometry that fixes a parallel set.

\subsection{The topology of $\mathcal{X}_0(\Gamma,\mathrm{Isom}(X))$}

\begin{Proposition}
    \label{Proposition:TopologyX_0}
    $\mathcal{X}_0(\Gamma,\mathrm{Isom}(X))\cong [0,1]\times [0,+\infty) .
$
\end{Proposition}

\begin{proof}
By Corollary~\ref{Corollary:Fix}, the subset of closed orbits in 
$\hom_0(\Gamma,\mathrm{Isom}(X))$ is
\begin{multline*}
\hom^{\mathrm{closed}}_0(\Gamma,\mathrm{Isom}(X))=
\\
\{\rho\in \hom_0(\Gamma,\mathrm{Isom}(X))\mid
\partial_\infty\mathrm{Fix}(\rho(a))\cap  
\partial_\infty \mathrm{Fix}(\rho(b))\neq 
\textrm{one-point set}\}    
\end{multline*}
The condition that
 $
 \partial_\infty\mathrm{Fix}(\rho(a))\cap
 \partial_\infty \mathrm{Fix}(\rho(b))
 $
is not a one-point set means that either neither end of
\(\mathrm{Fix}(\rho(b))\) is asymptotic to \(\mathrm{Fix}(\rho(a))\), or both are. In the latter case, the two lines are parallel by Lemma~\ref{Lemma:asymptotic}.

We aim to compute the quotient:
$$\mathcal{X}_0(\Gamma,\mathrm{Isom}(X))\cong
\hom^{\mathrm{closed}}_0(\Gamma,\mathrm{Isom}(X))\big/
\mathrm{Isom}(X).
$$
We consider the level sets of the distance function
$$
L_t=\{\rho\in\hom^{\mathrm{closed}}_0(\Gamma,\mathrm{Isom}(X))
\mid d (\mathrm{Fix}(\rho(a)),  
\mathrm{Fix}(\rho(b)))=t \}\qquad 
\textrm{ for }t\geq 0.
$$
Each $L_t$ is  $\mathrm{Isom}(X)$-invariant and 
$$\hom^{\mathrm{closed}}_0(\Gamma,\mathrm{Isom}(X))=\bigcup_{t\geq 0} L_t,$$ so we are interested in the quotient
$L_t\big/\mathrm{Isom}(X)$. 
Since we are assuming that 
$\mathrm{Fix}(\rho(a))$  and $  
\mathrm{Fix}(\rho(b))$ are not asymptotic at precisely one
ideal point, by Lemma~\ref{Lemma:asymptotic}:
\begin{itemize}
 \item For $\rho\in L_0 $,  $\mathrm{Fix}(\rho(a))\cap \mathrm{Fix}(\rho(b)) \neq \emptyset$.
    \item When $t>0$, for $\rho\in L_t $,  
    $\mathrm{Fix}(\rho(a))$ and $\mathrm{Fix}(\rho(b))$ have a common perpendicular.
\end{itemize}
Consider $L_t'$ the subset of representations $\rho$ in $L_t$
such that 
$\mathrm{Fix}(\rho(a))$ is a fixed parallel set
$$
\mathrm{Fix}(\rho(a))=P_0=\exp(\mathfrak{h}\cap\mathfrak{p}),
$$
for $\mathfrak{h}$   as in \eqref{eqn:h},
and that the class of the identity $x_0\in P_0$ minimizes 
in $P_0$
the distance
to $
\mathrm{Fix}(\rho(b))
$, so that $T_{x_0}X\cong \mathfrak{p}$.
Namely
$$
L'_t=\big\{\rho\in L_t\mid \mathrm{Fix}(\rho(a))=P_0\textrm{ and } d\big(x_0, \mathrm{Fix}(\rho(b)\big)=t \big\}.
$$

\begin{Lemma}
    View $\mathrm{O}(2)$ as a subgroup of $\mathrm{SO}(3)$
and let $\langle \mathrm{O}(2), i_{x_0}\rangle$ be the group 
generated by $\mathrm{O}(2)$ and the involution $i_{x_0}$.
Then 
$$
L_t\big/\mathrm{Isom}(X)=L_t'\big/\langle \mathrm{O}(2), i_{x_0}\rangle
$$
\end{Lemma}

\begin{proof}
    Every representation in $L_t$ can be conjugated to a representation in $L_t'$, and the lemma follows because 
    $\langle \mathrm{O}(2), i_{x_0}\rangle$ is the stabilizer of the pair $(P_0,x_0)$.
\end{proof}

% Since the group $\langle \mathrm{O}(2), i_{x_0}\rangle$ is compact,
% it follows from this lemma that all $\mathrm{Isom}(X)$-orbits in 
% $
% \hom^{\mathrm{closed}}_0(\Gamma,\mathrm{Isom}(X))
% $
% are closed. Combined with Proposition~\ref{Proposition:Accumulate}, we get:

% \begin{Corollary}
% The subspace of representations in 
% $
% \hom_0(\Gamma,\mathrm{Isom}(X))
% $ with closed orbits  is 
%   $
% \hom^{\mathrm{closed}}_0(\Gamma,\mathrm{Isom}(X))
% $.
% \end{Corollary}

The condition $\mathrm{Fix}(\rho(a))=P_0$ implies that,
for $\rho\in L_t'$, $\rho(a)=\mathrm{diag}(1,-1,-1)$. So we are interested in the expression of $\rho(b)$. When $\rho\in L_0'$, then $\rho(b)\in \mathrm {SO}(3)$, and when $\rho\in L_t'$, 
$\rho(b)$ is conjugate to a rotation by a transvection 
in a direction perpendicular to both 
$\mathrm{Fix}(\rho(a))$ and 
$\mathrm{Fix}(\rho(b))$. We describe it in the following 
lemma, whose proof is elementary (see 
Lemmas~\ref{Lemma:3rotations} and \ref{Lemma:axis} for 
notation).

\begin{Lemma}
\label{Lemma:parametersL'}
\begin{enumerate}[(a)]
    \item For $t> 0$, a representation $\rho \in L'_t$ satisfies
     $$
    \rho(b)=\exp( t \, m) \operatorname{Rot}(v)\exp(-t\, m)
    $$
    for uniquely determined $m\in(\mathfrak{h}\cap \mathfrak{p}) ^\perp$
    and $v\in S^2$,  such that  
      $|m|=1$ and 
    $$\mathrm{Axis}(\mathrm{Rot(v))\perp m}.$$
     \item A representation $\rho \in L'_0$ satisfies
     $$
    \rho(b)= \operatorname{Rot}(v)
    $$
    for a unique  $v\in S^2$.
\end{enumerate}   
\end{Lemma}

We need to take into account the action of the stabilizer 
$\langle \mathrm{O}(2), i_{x_0}\rangle$. Given $\rho(b)$ as in the lemma and  $A\in \mathrm{O}(2)\subset \mathrm{SO}(3)$:
\begin{align}
% \textrm{For }A\in \mathrm{O}(2),\qquad 
A\rho(b) A^{-1} &=\exp(t\operatorname{Ad}_A (m))
    \operatorname{Rot}(A v)
    \exp(-t\operatorname{Ad}_A (m)), \\
    i_{x_0}\rho(b) i_{x_0}^{-1} & = \exp( t \, (-m)) \operatorname{Rot}(v)\exp(-t\, (-m)).
\end{align} 
Lemma~\ref{Lemma:parametersL'} gives a way to describe the topology of each $L'_t$:

\begin{Corollary}
% \label{Corollary:Lt}
\begin{enumerate}[(a)]
    \item For $t> 0$, we have an equivariant homeomorphism
    $$L'_t\cong  
    \{(m,v)\in 
    (\mathfrak{h}^\perp\cap\mathfrak{p})\times S^2\mid
|m|=1, \ \mathrm{Axis}(\mathrm{Rot(v))\perp m}
\}
$$
where the group $\langle\mathrm{O}(2),i_{x_0}\rangle$ acts as follows:
\begin{align*}
 \textrm{For }A\in \mathrm{O}(2),\qquad 
A\cdot ( m,v)&=(\operatorname{Ad}_A (m),A v),
     \\
    i_{x_0}\cdot (m,v) & = (-m, v).
\end{align*} 

% $\big/ \langle \mathrm{SO}(1,2)\cap\mathrm{SO}(3), i_{x_0}\rangle
%     $$
\item We have an equivariant homeomorphism
$L'_0\cong S^2$, where 
the action of  $\mathrm{O}(2)$ on $S^2$ is standard and the action of $i_{x_0}$ is trivial.
\end{enumerate}    
\end{Corollary}

Notice that $\mathfrak{h}^\perp\cap\mathfrak{p}=\langle f_1,f_2\rangle$, where
$f_1$ and $f_2$ are as in \eqref{eqn:f}. Since $\mathrm{O}(2)$
acts transitively on the unit circle of $\mathfrak{h}^\perp\cap\mathfrak{p}$, 
up to conjugation we may choose $m= f_1$ and for 
$v=(x\; y\; z)^{\mathrm{T}}\in S^2$, the condition
$\mathrm{Axis}(\mathrm{Rot(v))\perp f_1}$
becomes $xy=0$.
Thus for $t\geq 0$
\begin{multline}
\label{eqn:quptoentLt}
L_t'\big/ \langle \mathrm{O}(2), i_{x_0}\rangle
\cong\{
v\in S^2\mid
\mathrm{Axis}(\mathrm{Rot}(v)) \perp f_1 
\}\big/K_2
\\
\cong\{
(x,y,z)\in\mathbb R^3\mid x^2+y^2+z^2=1,\ xy=0
\}\big/K_2    
\end{multline}
where 
\begin{equation}
    \label{eqn:K2}
    K_2=\left\{
\begin{pmatrix}
    1 & & \\
    & 1 & \\
    & & 1
\end{pmatrix},\ 
\begin{pmatrix}
    1 & & \\
    & -1 & \\
    & & -1
\end{pmatrix},\ 
\begin{pmatrix}
    -1 & & \\
    & 1 & \\
    & & -1
\end{pmatrix},\ 
\begin{pmatrix}
    -1 & & \\
    & -1 & \\
    & & 1
\end{pmatrix}
\right\}
\end{equation}
% $K_2\cong \mathbb{Z}_2 \times \mathbb{Z}_2$
% is as in ***. 
The set $\{
(x,y,z)\in\mathbb R^3\mid x^2+y^2+z^2=1,\ xy=0
\}$ is precisely the union of two maximal circles
$x=0$ and $y=0$, both invariant by $K_2$.
Each quotient 
$$
(S^2\cap\{x=0\})/K_2\qquad\textrm{ and }\qquad 
(S^2\cap\{y=0\})/K_2
$$
is  isometric to a quarter-circle.
More precisely consider the  meridian arcs of length $\pi/2$:
$$
Q_1=\{(x,y,z)\in S^2\mid x=0, y\geq 0, z\geq 0\}
\quad\textrm{ and }\quad
Q_2=\{(x,y,z)\in S^2\mid x\geq 0, y= 0, z\geq 0\}.
$$
The intersection is $Q_1\cap Q_2=\{(0,0,1)\}$, so
$Q_1\cup Q_2$ is homeomorphic to an interval. From elementary considerations we have:

\begin{Lemma}
\label{Lemma:Q1Q2}
The inclusion $Q_1\cup Q_2\subset S^2$ 
composed with \eqref{eqn:quptoentLt}
induces a homeomorphism 
$$Q_1\cup Q_2\cong L'_t\big/\langle \mathcal{O}(2), i_{x_0}\rangle.$$  
In particular $L'_t\big/\langle \mathcal{O}(2), i_{x_0}\rangle$
is homeomorphic to an interval. 
\end{Lemma}

% Thus $L'_t\big/\mathrm{O}(2)$ is the union of two intervals
% $$L'_t\big/\mathrm{O}(2)\cong \left( (S^2\cap\{x=0\})/K_2\right) \, \cup \,
% \left( (S^2\cap\{y=0\})/K_2\right) \cong [0,1]
% $$

When $t=0$ the quotient is also an interval: 
$$L_0'/\mathrm{O}(2)\cong S^2/\mathrm{O}(2)\cong [0,1].$$ 
We aim to know how the set 
$$
\bigcup_{t>0}{L'}_t\big/\langle\mathrm{O}(2), i_{x_0}\rangle\cong [0,1]\times (0,+\infty) 
$$
is attached to $L_0'/\mathrm{O}(2)$. Consider the 
 first coordinate map $x\colon S^2\subset 
 \mathbb R^3\to [-1,1]$ and 
its absolute value $|x|\colon S^2\to [0,1]$; it induces 
 a homeomorphism 
 $$|x|\colon S^2/\mathrm{O}(2)\xrightarrow{\cong} [0,1].
 $$
Because of Lemma~\ref{Lemma:Q1Q2}, 
the attaching of 
$L'_t\big/\langle \mathcal{O}(2), i_{x_0}\rangle$ to
$L_0'/\mathrm{O}(2)\cong [0,1]$ when $t\searrow 0$ is described by the map
 $$|x|\colon Q_1\cup Q_2\to  [0,1].
 $$
The map $\vert x\vert$ collapses $Q_1$ to a point and restricts to a homeomorphism $Q_2\cong [0,1]$.
 It follows that the total space is homeomorphic to 
 $[0,1]\times [0,+\infty)$. 
This finishes the proof of Proposition~\ref{Proposition:TopologyX_0}.
\end{proof}

\subsection{Fuchsian and elliptic loci of 
$\mathcal{X}_0(\Gamma,\mathrm{Isom}(X))$}

Set 
\(
\mathcal{E}_0 = \mathcal{E} \cap \mathcal{X}_0(\Gamma,\mathrm{Isom}(X))
\), that is 
 $\mathcal{E}_0$ denotes the elliptic locus in 
$\mathcal{X}_0(\Gamma,\mathrm{Isom}(X))$. Define $\mathcal{F}_{\mathrm{I},0}$ and $\mathcal{F}_{\mathrm{II},0}$ analogously. By construction, each of these sets is a closed subset of 
$\mathcal{X}_0(\Gamma,\mathrm{Isom}(X))$.

\begin{Proposition}
\label{Proposition:Fuchsian0}
    The homeomorphism $\mathcal{X}_0(\Gamma,\mathrm{Isom}(X))\cong [0,1]\times [0,+\infty)$ maps:
    \begin{enumerate}[(i)]
        \item $\mathcal{E}_0$ to $[0,1]\times \{0\}$,
        \item $\mathcal{F}_{\mathrm{II},0}$ to $\{0,1\}\times [0,+\infty)$, and
        \item $\mathcal{F}_{\mathrm{I},0}$ to  a proper half-line 
        that intersects $\partial\big( [0,1]\times [0,+\infty)\big)$ at $(0,0)$ and projects homeomorphically to the second coordinate $[0,+\infty)$.
    \end{enumerate}    
\end{Proposition}

We know that points away from the Fuchsian and elliptic 
locus have a natural analytic structure, Proposition~\ref{Proposition:NonEllipticNonFuchsianSmooth}. Besides,  points of the interior of 
$\mathcal{F}_{\mathrm{II},0}$ and the interior of 
$\mathcal{E}_0$ have a local structure of 
analytic orbifold (and the isotropy group equals the stabilizer of the representation, that can be checked to be  $ 
\mathbb{Z}_2$).

\begin{figure}[ht]
\begin{center}
  \begin{tikzpicture}[line join = round, line cap = round, scale=1]
\draw [thick] (-1,3)--(-1,0)--(1,0)--(1,3); 
\draw [thick] (-1,0) to [out=45, in=270] (0,3); 
 \fill[color=gray!20!white, opacity=.3] (-1,3)--(-1,0)--(1,0)--(1,3);
\node at (0,-.3) {$\mathcal{E}_0$};
 \node at (1.5,1.5) {$\mathcal{F}_{\mathrm{II},0}$};
 \node at (-1.5,1.5) {$\mathcal{F}_{\mathrm{II},0}$};
\node at (-.4,2.5) {$\mathcal{F}_{\mathrm{I},0}$};
  \end{tikzpicture}
  \end{center}
  \caption{The elliptic and Fuchsian locus in
  $\mathcal{X}_0(\Gamma,\mathrm{Isom}(X))$.}
  \label{Figure:EF0}
  \end{figure}

\begin{proof}
    We use the constructions in the proof of Proposition~\ref{Proposition:TopologyX_0}. In particular 
    the elliptic locus is precisely $L_0'/\mathrm{O}(2)$, thus (i) is clear.

    We describe the Fuchsian locus by analyzing
    its intersection with the quotient of the level sets $L_t$, where  
    $t=d(\mathrm{Fix}(\rho(a)),\mathrm{Fix}(\rho(b)))>0$.
    In the notation of the proof of 
    Proposition~\ref{Proposition:TopologyX_0}, we reduce to considering representations $\rho$ in 
    $L'_t$ for $t>0$. 
    An invariant Fuchsian plane $\mathbb{H}^2$ must contain 
    the geodesic through $x_0$ tangent to the vector $m$ 
    (see Lemma~\ref{Lemma:parametersL'} for the notation), 
    because  $\mathbb{H}^2$ intersects  both $\mathrm{Fix}
    (\rho(a))$ and 
    $\mathrm{Fix}(\rho(b))$ perpendicularly.  Then we apply 
    Lemma~\ref{Lemma:fewplanes}: there are precisely three Fuchsian planes that contain this line, one of type I and two of type II. In the proof of Proposition~\ref{Proposition:TopologyX_0}, we reduce to $$m=f_1=\begin{pmatrix}
        0& 1& 0 \\ 1 & 0 & 0 \\ 0 & 0 & 0
    \end{pmatrix}.
    $$ 
    So we may give explicitly the Fuchsian planes containing the line tangent to $m=f_1$: they are the planes tangent to  $\langle f_1, w_i\rangle$ 
    for $i=1,2,3$, 
    where
\begin{align*}
    w_1&%=\tfrac{1}{2}p_0-p_1
    =\begin{pmatrix}
        1& 0& 0 \\ 0 & -1 & 0 \\ 0 & 0 & 0
    \end{pmatrix}\textrm{ when }\mathbb{H}^2\textrm{ is of type I, and} \\
    w_2 &%f_2
    = \begin{pmatrix}
        0& 0& 1 \\ 0 & 0 & 0 \\ 1 & 0 & 0
    \end{pmatrix}
    \quad\textrm{ or }\quad
    w_3= %2 p_2 = 
    \begin{pmatrix}
        0& 0& 0 \\ 0 & 0 & 1 \\ 0 & 1 & 0
    \end{pmatrix}\textrm{ when }\mathbb{H}^2\textrm{ is of type II.}
\end{align*}
As $\rho(b)$ preserves $\mathbb{H}^2$, 
by Lemma~\ref{Lemma:uniquerotations}
we know that there is only one possible $v\in S^2$ (up to sign) for 
which $\rho(b)=\exp(-t\, m)\mathrm{Rot}(v)\exp(t\, m)$. We can give 
these rotations explicitly: $\mathrm{Rot}(v_i)$ preserves the Fuchsian plane tangent to $\langle f_1, w_i\rangle$, where   
$$
v_1=\pm \begin{pmatrix}
    0 \\ 0\\1
\end{pmatrix},
\qquad
v_2=\pm \begin{pmatrix}
    1 \\ 0\\0
\end{pmatrix},
\textrm{ and }
\qquad
v_3=\pm \begin{pmatrix}
    0 \\ 1\\0
\end{pmatrix}.
$$
Let $Q_1$ and $Q_2$ be as in Lemma~\ref{Lemma:Q1Q2}; 
$v_2$ and $v_3$ are the 
end-points of the 
union $Q_1\cup Q_2$ and $\{v_1\}=Q_1\cap Q_2$. As $Q_1$ collapses to a 
point when $t=0$, the assertion follows.    
\end{proof}

By examining the adjoint action of $\rho(a)=\mathrm{diag}(1,-1,-1)$ on the elements $m$, $w_1$, $w_2$, and $w_3$ introduced in the proof of Proposition~\ref{Proposition:Fuchsian0}, we obtain the following observation:

\begin{Remark}
In the Fuchsian setting considered in Proposition~\ref{Proposition:Fuchsian0}, the corresponding action preserves the orientation of the invariant Fuchsian plane $\mathbb{H}^2$ when it is of type~II, and does not preserve the orientation when it is of type~I.
\end{Remark}

\section{The component $\mathcal{X}_2(\Gamma,\mathrm{Isom}(X))$}
\label{Section:X2}

This section is divided into three subsections. In the first two, we compute the variety of characters; in the final one, we describe the elliptic and Fuchsian locus.

The component $\hom_2(\Gamma,\mathrm{Isom}(X))$ is precisely the set of representations sending $a$ to a reflection across a totally geodesic hyperbolic plane of type II. In particular, $\mathrm{Fix}(\rho(a))$ is a Fuchsian plane of type II, while $\mathrm{Fix}(\rho(b))$ is a singular geodesic.

By Corollary~\ref{Corollary:hom12closed}
 all orbits in 
 $\hom_2(\Gamma,\mathrm{Isom}(X))$ are closed.

\subsection{Topology of the level sets}

% The fact that all orbits are closed relies in the following observation, that is also used to construct the
% level sets.

% \begin{Lemma}
% \label{Lemma:NoAsymptotic}
%     For $\rho\in \hom_2(\Gamma,\mathrm{Isom}(X)) $,  $\mathrm{Fix}(\rho(a))$ and
%     $\mathrm{Fix}(\rho(b))$ are not asymptotic.
%     In particular either $\mathrm{Fix}(\rho(a))$ and
%     $\mathrm{Fix}(\rho(b))$ intersect, either they have a common perpendicular.
% \end{Lemma}

% \begin{proof}
%     By Lemma~\ref{Lem:FuchsianRegular} every point in 
%     $\partial_\infty \mathrm{Fix}(\rho(a))$
%     has type the center of the model spherical Weyl chamber $\sigma_{\mathrm{mod}}$.
%     As $\mathrm{Fix}(\rho(b))$ is a singular geodesics, its ideal end-points
%     $\partial_\infty\mathrm{Fix}(\rho(b))$
%     are singular, hence they have type in the 
%     $\partial \sigma_{\mathrm{mod}}$.
%     So $\partial_\infty \mathrm{Fix}(\rho(a))
%     \cap 
%     \partial_\infty\mathrm{Fix}(\rho(b))=\emptyset$, which means 
%     that they are not asymptotic. 
%     The second assertion follows from convexity of the distance function.
% \end{proof}

As in the proof of Proposition~\ref{Proposition:TopologyX_0},
consider the level subsets of the distance function between
$\mathrm{Fix}(\rho(a))$ and $\mathrm{Fix}(\rho(b))$. Namely, 
for $t\geq 0 $ define the $\mathrm{Isom}(X)$-invariant subsets
$$
L_t=\{\rho\in\hom_2(\Gamma,\mathrm{Isom}(X))\mid 
d\bigl(\mathrm{Fix} (\rho(a)),\mathrm{Fix} (\rho(b))\bigr)=t\}.
$$
As for $\mathcal{X}_0(\Gamma,\mathrm{Isom}(X))$, 
we     
shall study every quotient $L_t\big/\mathrm{Isom}(X)$.
By Corollary~\ref{Corollary:hom12closed}
all conjugation orbits  in 
$\hom_2(\Gamma,\mathrm{Isom}(X))$
are closed, so
$$\mathcal{X}_2(\Gamma,\mathrm{Isom}(X))=\bigcup_{t\geq 0}
L_t\big/\mathrm{Isom}(X)
.
$$
In this subsection we aim to determine the topology of
$$
L_t\big/\mathrm{Isom}(X).
$$

By Lemma~\ref{Lemma:NoAsymptotic}, $\mathrm{Fix} (\rho(a))$ and 
$\mathrm{Fix} (\rho(b))$ are not asymptotic; 
thus if they are at distance zero, then they intersect and we have:
\begin{equation}
    \label{eqn=NoAsymptotics}
    L_0=\{\rho\in\hom_2(\Gamma,\mathrm{Isom}(X))\mid 
\mathrm{Fix} (\rho(a))\cap \mathrm{Fix} (\rho(b))\neq\emptyset\}.
\end{equation}
To determine the quotients 
$L_t\big/\mathrm{Isom}(X)$, for $t\geq 0$
we consider $L_t'$ the subset of representations $\rho$ in $L_t$ that satisfy the following two properties:
\begin{enumerate}[(a)]
    \item  $\mathrm{Fix}(\rho(a))$ is a \emph{fixed}
    Fuchsian plane of type II 
    (see \eqref{eqn:f} for notation):
    $$\mathrm{Fix}(\rho(a))=\mathbb H^2_{\mathrm{II},0}=\exp(\langle f_1,f_2\rangle ).
    $$
    \item The class of the identity $x_0=[\mathrm{Id}]$ minimizes the distance to $\mathrm{Fix}(\rho(b))$ on $\mathrm{Fix}(\rho(a))=\mathbb H^2_{\mathrm{II},0}$
    (and equals the intersection when the distance is zero).
\end{enumerate}
So the  subset of representations that satisfy (a) and (b) is
$$
L_t'=\{ \rho\in L_t\mid \mathrm{Fix}(\rho(a))=  \mathbb H^2_{\mathrm{II},0}\textrm{ and } d(x_0, \mathrm{Fix}(\rho(b)))= t\}.
$$

    Let $\langle \mathrm{SO}(1,2)\cap \mathrm{SO}(3),i_{x_0}\rangle$ 
    denote the group 
generated by $\mathrm{SO}(1,2)\cap \mathrm{SO}(3)\cong \mathrm{O}(2)$ and the involution $i_{x_0}$.

\begin{Lemma}
With the previous notation, for $t\geq 0$:
$$
L_t/\mathrm{Isom}(X)=L_t'/\langle \mathrm{SO}(1,2)\cap \mathrm{SO}(3),i_{x_0}\rangle .
$$
\end{Lemma}

\begin{proof}
    Every representation in $L_t$ can be conjugated to a representation in $L_t'$, because there is a unique congruence class of Fuchsian planes of type II, and the lemma follows because 
    $\langle \mathrm{SO}(1,2)\cap \mathrm{SO}(3), i_{x_0}\rangle$ is the stabilizer of the pair $(\mathbb H^2_{\mathrm{II},0},x_0)$.
\end{proof}

The condition $\mathrm{Fix}(\rho(a))=
{\mathbb H^2_{\mathrm{II},0}}$ implies that,
for $\rho\in L_t'$, 
$$\rho(a)([M])=[\mathrm{diag}(1,-1,-1) M^*]\qquad 
\text{ for every class }[M]\in X.$$ So we are interested in the expression of $\rho(b)$. When $\rho\in L_0'$, 
$\rho(b)\in \mathrm {SO}(3)$, and when $\rho\in L_t'$, 
$\rho(b)$ is conjugate to a rotation in 
$\mathrm{SO}(3)$
by a transvection in a direction perpendicular to both $\mathrm{Fix}(\rho(a))$ and 
$\mathrm{Fix}(\rho(b))$. We describe it in the following lemma, 
as we did in Lemma~\ref{Lemma:parametersL'}.

\begin{Lemma}
\label{Lemma:parametersL'2}
\begin{enumerate}[(a)]
    \item For $t>0$, a representation $\rho \in L'_t$ satisfies
\[
\rho(b)=\exp(t\,m)\,\operatorname{Rot}(v)\,\exp(-t\,m),
\]
for uniquely determined elements $m\in\mathfrak{h}\cap\mathfrak{p}=(\mathbb H^2_{\mathrm{II},0})^\perp$ and $v\in S^2$ such that $|m|=1$ and the rotation axis of $\operatorname{Rot}(v)$ is orthogonal to $m$, i.e.
\[
\mathrm{Axis}(\operatorname{Rot}(v)) \perp m.
\]
     \item A representation $\rho \in L'_0$ satisfies
     $$
    \rho(b)= \operatorname{Rot}(v)
    $$
    for a unique  $v\in S^2$.
\end{enumerate}   
\end{Lemma}

Again we take into account the action of 
$\langle \mathrm{SO}(1,2)\cap\mathrm{SO}(3), 
i_{x_0}\rangle$, 
the stabilizer of the pair 
$((\mathbb H^2_{\mathrm{II},0}),
x_0)$. Given $\rho(b)$ as in
Lemma~\ref{Lemma:parametersL'2} and  $A\in \mathrm{SO}(1,2)\cap\mathrm{SO}(3)$:
\begin{align}
% \textrm{For }A\in \mathrm{O}(2),\qquad 
A\rho(b) A^{-1} &=\exp(t\operatorname{Ad}_A (m))
    \operatorname{Rot}(A v)
    \exp(-t\operatorname{Ad}_A (m)), \\
    i_{x_0}\rho(b) i_{x_0}^{-1} & = \exp( t \, (-m)) \operatorname{Rot}(v)\exp(-t\, (-m)).
\end{align} 
% Lemma~\ref{Lemma:parametersL'2} gives a way to describe the topology each $L'_t$:

\begin{Corollary}
\label{Corollary:Lt}
\begin{enumerate}[(a)]
    \item For $t> 0$, we have an equivariant  homeomorphism
    $$L'_t\cong  
    \{(m,v)\in \left(\mathbb{H}^2_{\mathrm{II},0} \right)^\perp\times S^2\mid
|m|=1, \ \mathrm{Axis}(\mathrm{Rot(v))\perp m}
\},
$$
where the group $\langle\mathrm{SO}(1,2)\cap\mathrm{SO}(3),i_{x_0}\rangle$ acts as follows:
\begin{align*}
 \textrm{For }A\in \mathrm{SO}(1,2)\cap\mathrm{SO}(3),\qquad 
A\cdot ( m,v)&=(\operatorname{Ad}_A (m),A v),
     \\
    i_{x_0}\cdot (m,v) & = (-m, v).
\end{align*} 

% $\big/ \langle \mathrm{SO}(1,2)\cap\mathrm{SO}(3), i_{x_0}\rangle
%     $$
\item We have an equivariant homeomorphism
$L'_0\cong S^2$, where 
the action of  $\mathrm{O}(2)$ on $S^2$ is standard and the action of $i_{x_0}$ is trivial.
\end{enumerate}    
\end{Corollary}

Up to this point, our argument has closely mirrored the proof of Proposition~\ref{Proposition:TopologyX_0}, but here the paths diverge and differences begin. We start 
describing the space
\begin{equation}
    \label{eqn:quotient2}
   \{m\in  \mathfrak{p}\mid
m\perp \mathbb{H}^2_{\mathrm{II},0},\ 
|m|=1 
\}\big/ \langle \mathrm{SO}(1,2)\cap\mathrm{SO}(3), i_{x_0}\rangle
\end{equation}
Consider the (tangent space of the) maximal flat
\begin{equation}
\label{eqn}
F=\{\mathrm{diag}(\lambda_1,\lambda_2,\lambda_3)\mid
\lambda_1+\lambda_2+\lambda_3=0\}\subset\mathfrak{p}.
\end{equation}
We identify the spherical apartment $\partial_\infty F$ with the unit circle in $F$ given by
\begin{equation}
    \label{eqn:F}
    F\cap\{\lambda_1^2+\lambda_2^2+\lambda_3^2=1\}.
\end{equation}
Furthermore, consider the quarter-circle
\[
Q=F\cap\{\lambda_1^2+\lambda_2^2+\lambda_3^2=1\}
\cap\{\lambda_2\geq\lambda_3\}\cap\{\lambda_1\geq 0\},
\]
which is an arc of length $\pi/2$, whose endpoints are given by the intersections with the hyperplanes
$\{\lambda_1=0\}$ and $\{\lambda_2=\lambda_3\}$, respectively.
Let $\alpha\in[0,\pi/2]$ denote the arc-length parameter on $Q$, and let 
$q_\alpha\in Q$ be the point at arc-length $\alpha$, so that 
$$
q_0=Q\cap\{\lambda_1=0\}
\qquad\text{and}\qquad
q_{\pi/2}=Q\cap\{\lambda_2=\lambda_3\}.
$$
See Figure~\ref{Figure:Q}.

\begin{figure}[ht]
\begin{center}
  \begin{tikzpicture}[line join = round, line cap = round, scale=.7]
    \draw[very thick] (-2,0)--(2,0);
    \draw[very thick, rotate=60] (-2,0)--(2,0);
    \draw[very thick, rotate=-60] (-2,0)--(2,0);
        \node at (3,0) {${\lambda_2=\lambda_3}$};
        \node at (1.2,2) {${\lambda_1=\lambda_2}$};
        \node at (1.3,-2.2) {${\lambda_1=\lambda_3}$};
    \draw[thick, dashed] (0,-2.5)--(0,2.5);
    \draw[thick, dashed, rotate=60] (0,-2.5)--(0,2.5);
    \draw[thick, dashed, rotate=-60] (0,-2.5)--(0,2.5);
        \node at (0,3) {${\lambda_1=0}$};
        \node at (3.1,1.2) {${\lambda_2=0}$};
        \begin{scope}[shift={(8,-1.5)}]
        \draw[thick] (0,0)--(3,0);
        \draw[thick, rotate=60] (0,0)--(3,0);
        \draw[thick, dashed] (0,0)--(0,3);
        \draw[thick, dashed, rotate=-60] (0,0)--(0,3);
        \draw[ultra thick] (3,0) arc[start angle=0, end angle=90, radius=3];
        \draw[black, fill=black] (3,0) circle[radius=.09];
        \draw[black, fill=black] (0,3) circle[radius=.09];
    \node at (-1,1.75) {${\lambda_1=0}$};
    \node at (1.5,-.5) {${\lambda_2=\lambda_3}$};
    \node at (0,3.5) {$q_0$};
    \node at (3.7,0) {$q_{\pi/2}$};
    \node at (3,2.5) {${Q}$};
\end{scope}
  \end{tikzpicture}
  \end{center}
  \caption{The  plane 
  {${\lambda_1+\lambda_2+\lambda_3=0}$} with the three 
  singular lines through the origin identified 
  with the 
  maximal flat $F$, on the left, and the 
  quarter-circle $Q$ bounded by $\lambda_1=0$ and
  $\lambda_2=\lambda_3$ on the right.}
  \label{Figure:Q}
  \end{figure}

\begin{Lemma}
\label{Lemma:Qhomeo}
    The inclusion $Q\subset
        \{m\in  \mathfrak{p}\mid
m\perp \mathbb{H}^2_{\mathrm{II},0},\ 
|m|=1 
\}
$
induces a homeomorphism
    $$ Q\cong 
     \{m\in  \mathfrak{p}\mid
m\perp \mathbb{H}^2_{\mathrm{II},0},\ 
|m|=1 
\}\big/ \langle \mathrm{SO}(1,2)\cap\mathrm{SO}(3), i_{x_0}\rangle.
$$
\end{Lemma}

\begin{proof}
    The subspace $(\mathbb{H}^2_{\mathrm{II},0})^\perp$ coincides with the tangent space to a parallel set of the form $\mathfrak{h}\cap\mathfrak{p}$ as described in \eqref{eqn:h}. Exploiting the action of the subgroup $SO(2)\subset SO(1,2)$, we may, for the purpose of computing the orbit space in \eqref{eqn:quotient2}, assume without loss of generality that $m$ is diagonal. More precisely, the point $m$ lies in an apartment isometric to $S^1$, which is the ideal boundary of a maximal flat $F$ as in \eqref{eqn:F}. 

This spherical apartment $\partial_\infty F$ is canonically identified with the unit circle
\[
F\cap\{\lambda_1^2+\lambda_2^2+\lambda_3^2=1\},
\]
and its stabilizer in the group
\[
\langle \mathrm{SO}(1,2)\cap\mathrm{SO}(3), i_{x_0}\rangle
\]
is the order–four subgroup
\[
\{\mathrm{Id},\varsigma, i_{x_0},  \varsigma i_{x_0}\},
\]
where
\begin{equation}
    \label{eqn:varsigma}
    \varsigma=\begin{pmatrix}
        1 & 0 & 0\\[0.2em]
        0 & 0 & -1\\[0.2em]
        0 & 1 & 0
    \end{pmatrix}.
\end{equation}
% So $\varsigma $ is an element of the Weyl group, and  
% $\varsigma i_{x_0}$ is 
% another reflection 
% along a perpendicular axis and inducing the canonical 
% involution of the spherical Weyl chamber 
% $\sigma_{\mathrm{mod}}$.  
Thus the inclusion induces a homeomorphism
\begin{multline*}
F\cap\{\lambda_1^2+\lambda_2^2+\lambda_3^2=1\}
\big/\{\mathrm{Id},\varsigma, i_{x_0},  \varsigma i_{x_0}\}
\\
\cong
\{m\in  \mathfrak{p}\mid
m\perp \mathbb{H}^2_{\mathrm{II},0},\ 
|m|=1 
\}\big/
\langle \mathrm{SO}(1,2)\cap\mathrm{SO}(3), i_{x_0}\rangle  .  
\end{multline*}
The group
$
\{\mathrm{Id},\varsigma, i_{x_0},  \varsigma i_{x_0}\}
$
acts on $F\cap\{\lambda_1^2+\lambda_2^2+\lambda_3^2=1\}$
as a group of reflections with fundamental domain $Q$.
More precisely, $\varsigma$ is an element of a Weyl group 
(that fixes $\lambda_2=\lambda_3$) and $i_{x_0}$ is the antipodal map, so that $\varsigma i_{x_0}=i_{x_0}\varsigma$ is a reflection that fixes $\lambda_1=0$.
Thus the inclusion 
$Q\subset F\cap\{\lambda_1^2+\lambda_2^2+\lambda_3^2=1\}$
induces  a homeomorphism
$$
Q\cong F\cap\{\lambda_1^2+\lambda_2^2+\lambda_3^2=1\}
\big/\{\mathrm{Id},\varsigma, i_{x_0},  \varsigma i_{x_0}\}.
$$
\end{proof}

\begin{Lemma}
    For $t>0$, $L'_t\big/\langle \mathrm{SO}(1,2)\cap \mathrm{SO}(3),i_{x_0}\rangle    \cong S^2$.
\end{Lemma}

\begin{proof}
By Lemma~\ref{Lemma:Qhomeo}, we have a homeomorphism
\[
L'_t\big/\langle \mathrm{SO}(1,2)\cap 
\mathrm{SO}(3),i_{x_0}\rangle 
\cong \{ (q_\alpha,v)\in Q\times S^2 \mid 
q_\alpha\perp \mathrm{Axis}(\mathrm{Rot}(v)) \}
\big/\sim,
\]
where the equivalence relation $\sim$ is defined by
\[
(q_\alpha,v)\sim (q_\alpha',v') \iff 
\bigl(q_\alpha=q_\alpha' \text{ and } 
v'=\gamma v \text{ for some } \gamma\in \mathrm{Stab}(q_\alpha)\bigr).
\]
Here, $\mathrm{Stab}(q_\alpha)$ denotes the stabilizer subgroup of $q_\alpha\in Q$ in 
\(
\langle \mathrm{SO}(1,2)\cap\mathrm{SO}(3), i_{x_0}\rangle.
\)
     
    To describe  $\mathrm{Stab}(q_\alpha)$ consider
     the group $K_2$ defined in~\eqref{eqn:K2}.
% $$
% K_2=\left\{
% \begin{pmatrix}
%     1 & & \\
%     & 1 & \\
%     & & 1
% \end{pmatrix},\ 
% \begin{pmatrix}
%     1 & & \\
%     & -1 & \\
%     & & -1
% \end{pmatrix},\ 
% \begin{pmatrix}
%     -1 & & \\
%     & 1 & \\
%     & & -1
% \end{pmatrix},\ 
% \begin{pmatrix}
%     -1 & & \\
%     & -1 & \\
%     & & 1
% \end{pmatrix}
% \right\}
% $$
% The stabilizer $\mathrm{Stab}(q)$ in 
% $\langle \mathrm{SO}(1,2)\cap\mathrm{SO}(3),i_{i_0}\rangle$
% of a point in $q\in Q$ is:
\begin{enumerate}[(a)]
    \item When $q_\alpha\in \mathrm{int}(Q)$,
    e.g.~$0<\alpha<\pi/2$,
    $\mathrm{Stab}(q_\alpha)=K_2$.
    \item  When $\alpha=\pi/2$,
    % $\{q\}=Q\cap \{\lambda_2=\lambda_3\}$,
    $\mathrm{Stab}(q_{\pi/2})=\mathrm{SO}(1,2)\cap\mathrm{SO}(3)
    \cong O(2)$, in particular $K_2\subset \mathrm{Stab}(q_{\pi/2})$.
    \item When $\alpha=0$, 
    the stabilizer $\mathrm{Stab}(q_{0})$ is the subgroup generated by $K_2$ and the composition
\[
\varsigma i_{x_0} = i_{x_0} \varsigma,
\]
where $\varsigma$ denotes the Weyl reflection introduced in~\eqref{eqn:varsigma}, and $i_{x_0}$ denotes inversion about $x_0$. In particular, $\mathrm{Stab}(q_{0})$ is isomorphic to the dihedral group of order $8$.
\end{enumerate}
Notice that $S^2/K_2$ is an orbifold with underlying space $S^2$
and three branching points of order~$2$.

We first claim that
$
\{ (q_\alpha,v)\in Q\times S^2\mid 
    q_\alpha\perp\mathrm{Axis}(\mathrm{Rot}(v)) \}
$
projects in $Q\times (S^2/K_2)$ to 
a subset fiber-wise homeomorphic to  
a cylinder 
$$
Q\times S^1.
$$
To prove that claim, notice that for any $q_\alpha=\mathrm{diag}(\lambda_1,\lambda_2,\lambda_3)\in Q$
the set 
\begin{equation}
    \label{eqn:conic}
\{v\in S^2\mid  \mathrm{Axis}(\mathrm{Rot}(v))\perp q_\alpha\}=
\{(x,y,z)\in\mathbb R^3\mid  
x^2+y^2+z^2=1,\
\lambda_1 x^2+
\lambda_2 y^2+\lambda_3z^2=0\}    
\end{equation}
is a conic in $S^2$. It is degenerate precisely when the type 
of $q$ is the center of $\sigma_{\mathrm{mod}}$, namely when 
$\alpha=0$ or $\alpha=\pi/3$, as
$q_{0}=Q\cap\{\lambda_1=0\} $
and
$q_{\pi/3}=Q\cap\{\lambda_2=0\} $,
see Figure~\ref{Figure:Q}. When it is not degenerate, it is
homeomorphic to the disjoint union $S^1\sqcup S^1$. When it is degenerate, it is the union of two maximal circles in $S^2$ (that meet at two antipodal points).
The quotient $S^2/K_2$ is a
2-orbifold with underlying space again $S^2$, and the conic
\eqref{eqn:conic} (singular or not) projects to a circle $S^1$.
The projection of the conic meets the branching locus of 
$S^2/K_2$ precisely when it is degenerate.

After the claim, to obtain 
$L'_t\big/\langle \mathrm{SO}(1,2)\cap \mathrm{SO}(3),i_{x_0}\rangle$ we need to quotient further
the cylinder
$Q\times S^1$ at its boundary circles. Namely:
% For that purpose chose a parameterization of 
% $Q$ induced by the isometry $Q\cong [0,\pi/2]$, where 
% $q_0=Q\cap\{\lambda_2=\lambda_3\}$ and 
% $q_{{\pi}/{2}}=Q\cap \{\lambda_1=0\}$.
% So
$$
L'_t\big/\langle \mathrm{SO}(1,2)\cap \mathrm{SO}(3),i_{x_0}\rangle\cong
Q\times S^1/\!\sim ,$$ 
where $\sim$ is the relation induced by
the action of $\mathrm{SO}(2)$ on $\{q_{\pi/2}\}\times S^1$
and the action of the reflection $\varsigma$ on $\{q_0\}\times S^1$. Topologically, the end 
$\{q_{\pi/2}\}\times S^1$ collapses to a point, and the quotient of
$\{q_0\}\times S^1$ is homeomorphic to $[0,1]$, 
so  $Q\times S^1/\!\sim\;\cong S^2$.
\end{proof}

\subsection{The topology of 
$\mathcal{X}_2(\Gamma,\mathrm{Isom}(X))$}

\begin{Proposition}
    \label{Proposition:TopologyX_2}
    $\mathcal{X}_2(\Gamma,\mathrm{Isom}(X))\cong \mathbb R^3
$.
\end{Proposition}

\begin{proof}
We continue with the notation of the previous subsection.
For $t> 0$ we know that 
$L'_t\big/\langle \mathrm{SO}(1,2)\cap \mathrm{SO}(3),i_{x_0}\rangle\cong S^2$ and we aim to understand the attaching map to 
$L_0'\big/\mathrm{SO}(3)\cap\mathrm{SO}(1,2)$ when 
$t\searrow0$.
When $t=0$, we know by Corollary~\ref{Corollary:Lt} that 
$$
L_0'\big/\mathrm{SO}(3)\cap\mathrm{SO}(1,2)
\cong S^2\big/\mathrm{SO}(3)\cap\mathrm{SO}(1,2)\cong S^2\big/\mathrm{O}(2)\cong [0,1].
$$
Viewing $S^2$ as the unit sphere in $\mathbb R^3$, 
the homeomorphism $S^2\big/\mathrm{O}(2)\cong [0,1]$
is induced by the 
absolute value of the first coordinate map
$$
\begin{array}{rcl}
   |x|\colon S^2=\{x^2+y^2+z^2=1\} & \to & [0,1]  \\
     (x,y,z)& \mapsto & |x| 
\end{array}
$$
We look at the composition 
$$L_t'\hookrightarrow \mathfrak{p}\times  S^2
\xrightarrow{\mathrm{proj}}
S^2
\xrightarrow{|x|} [0,1]
$$
and denote  by 
$$
\phi
\colon L_t'\big/\mathrm{SO}(3)\cap\mathrm{SO}(1,2)\to [0,1]
$$ the induced map. 

\begin{Proposition} 
    \label{Proposition:Fibre}
    The fiber of $\phi\colon L_t'\big/\mathrm{SO}(3)\cap\mathrm{SO}(1,2)
    \to [0,1]$ is
    $$
    \phi^{-1}(x)\cong 
    \begin{cases}
        [0,1] & \textrm{ if }x=0 \\
        S^1 & \textrm{ if } 0<x<1 \\
        \textrm{a point}  & \textrm{ if }x=1
    \end{cases}
    $$
\end{Proposition}

\begin{proof}
    Write explicitly the parameterization of the 
    quarter-circle $Q$ by arc-length $\alpha\in [0,\pi/2]$:
$$
q_\alpha=
\begin{pmatrix}
    a_\alpha & & \\
    & b_\alpha & \\
    & & c_\alpha
\end{pmatrix}
=
\sin(\alpha) \frac{1}{\sqrt{6}} 
\begin{pmatrix}
    2 & & \\
    & -1 & \\
    & & -1
\end{pmatrix}
+
\cos(\alpha) \frac{1}{\sqrt{2}} 
\begin{pmatrix}
    0 & & \\
    & -1 & \\
    & & 1
\end{pmatrix}
$$
For $(x,y,z)\in S^2$, the condition 
$\mathrm{Axis}(\mathrm{Rot}((x\, y\, z)^{\mathrm{T}})\perp q_\alpha$
becomes the system of equations
\begin{equation}
    \label{eqn:system}
    \begin{cases}
        a_\alpha x^2+b_\alpha y^2+c_\alpha z^2 & = 0 \\
        x^2+y^2+z^2& =1
    \end{cases}
\end{equation}
For fixed $x\in [0,1]$ (the coordinate for $S^2/\mathrm{O}(2)$)
this system of equations on $y$ and $z$ (with parameters $x$ and $\alpha$) reads
\begin{equation}
    \label{eqn:systemy2z2}
    \begin{cases}
        y^2= & \frac{1}{2}+  \left(1-x^2-\frac{3 x^2-1}
        {\sqrt{3}}\tan(\alpha)\right)  \\
        z^2= & \frac{1}{2}+  \left(1-x^2+\frac{3 x^2-1}
        {\sqrt{3}}\tan(\alpha)\right) 
    \end{cases}
\end{equation}
We describe when this system has solutions, by distinguishing two cases:
\begin{enumerate}[(a)]
    \item For $0\leq x\leq \frac{1}{\sqrt{3}}$,
    the condition $z^2\geq 0$ becomes 
    $
    \tan(\alpha)\leq \frac{1-x^2}{1-3 x^2}\sqrt{3}
    $.
    \item For $\frac{1}{\sqrt{3}}\leq x \leq 1$,
    the condition $y^2\geq 0$ becomes 
    $
    \tan(\alpha)\leq \frac{1-x^2}{3 x^2-1}\sqrt{3}
    $.
\end{enumerate}
So define $A\colon [0,1]\to [0,\frac\pi2]$ 
by the formula:
$$
A(x)=
\begin{cases}
    \arctan\left( \frac{1-x^2}{1-3 x^2}\sqrt{3}  \right) &
    \textrm{ for }  0\leq x\leq \frac{1}{\sqrt{3}}
    \\
    \arctan\left( \frac{1-x^2}{3 x^2-1}\sqrt{3}  \right) &
    \textrm{ for }  \frac{1}{\sqrt{3}}\leq x\leq 1
\end{cases}
$$
For $0\leq \alpha \leq A(x)$ there are solutions $(y,z)$ to the system of equations \eqref{eqn:system}.

\begin{Lemma}
    \label{Lemma:projectRA}
    The projection to the coordinates $(x,\alpha)\in [0,1]\times [0,\tfrac{\pi}{2}]$ gives a map from  the space 
    $L_t'\big/\mathrm{SO}(3)\cap\mathrm{SO}(1,2)$ onto the  region 
    \[
        R_A=\{(x,\alpha)\in [0,1]\times [0,\tfrac{\pi}{2}] \mid
        \alpha\leq A(x)\}.
    \]
    See Figure~\ref{Figure:Ra}.
    This projection is a $2\!:\!1$ covering map over the interior of $R_A$ and is $1\!:\!1$ on the boundary $\partial R_A$.
\end{Lemma}

\begin{proof}[Proof of Lemma~\ref{Lemma:projectRA}]
The map $A$ has been built so that for each $(x,\alpha)\in R_A$
there exist solutions $(y,z)$ to \eqref{eqn:systemy2z2}, which
is  equivalent 
to~\eqref{eqn:system}. For $(x,\alpha)\in \mathring{R}_A$ there are
4 solutions 
for the pair $(y,z)$, namely $(\pm y, \pm z)$. The stabilizer of 
$q_\alpha$ is $K_2$, but since $x>0$ is fixed, 
we only quotient by the change of sign $(y,z)\mapsto (-y,-z)$. 
Therefore
this gives two solutions.
For $(x,\alpha)\in \partial {R}_A$ the discussion is 
similar, but several cases need to be considered (for 
instance  when $x=\frac{1}{\sqrt{3}}$ 
and $\alpha=\frac{\pi}{2}$) to show that here the 
projection is 1:1.
\end{proof}

\begin{figure}[ht]
\begin{center}
  \begin{tikzpicture}[line join = round, line cap = round, scale=2]
\fill[color=gray!20!white, opacity=.3] (0,pi/3) to [out=0, in=255] (0.577,pi/2) to [out=-65, in=115] (1,0)--(0,0);
  \draw[thick, -stealth] (0,0)--(1.2,0);
  \draw[thick, -stealth] (0,0)--(0,pi/2+.2);
  % %  %
    \node at (0.577,-.3) {$\frac{1}{\sqrt{3}}$};   
    \node at (0,-.2) {$0$};
    \node at (1,-.2) {$1$};
  %   %
     \node at (-.2,0) {$0$};
  %   \node at (-.7,1) {$\pi/6$};
     \node at (-.3,pi/3) {$\pi/3$};
     \node at (-.3,pi/2) {$\pi/2$};
    \node at (1.3,0) {$x$};
    \node at (0,pi/2+.3) {$\alpha$};
  \draw[thick] (0,pi/3) to [out=0, in=255] (0.577,pi/2) to [out=-65, in=115] (1,0);
  \draw[thin, dashed, opacity=.3] (0.577,0)--(0.577,pi/2)--(0,pi/2);
  \node at (.4,.6) {$R_A$};
% \end{scope}
  \end{tikzpicture}
  \end{center}
  \caption{The region $R_A=\{(x,\alpha)\in [0,1]\times [0,\tfrac{\pi}{2}] \mid
        \alpha\leq A(x)\}$.}
  \label{Figure:Ra}
  \end{figure}

\noindent \emph{End of the proof of Proposition~\ref{Proposition:Fibre}.}
Lemma~\ref{Lemma:projectRA} describes the projection from
$L_t'\big/\mathrm{SO}(3)\cap\mathrm{SO}(1,2)$  to 
$R_A$ (equivalent to the 
projection from $S^2\subset \mathbb R^3$ to the unit disc in $\mathbb{R}^2$) and $\vert x\vert$ is precisely a coordinate, so 
Proposition~\ref{Proposition:Fibre} follows easily.
\end{proof}

\noindent \emph{End of the proof of Proposition~\ref{Proposition:TopologyX_2}.}
We have described the quotient of each level subset
$$
L_t\big/\mathrm{Isom}(X)\cong L_t'/\langle\mathrm{SO}(1,2)\cap \mathrm{SO}(3),
i_{x_0} \rangle
$$
for $t\geq 0$. 
The map $\phi$ of Proposition~\ref{Proposition:Fibre} describes 
how $L_t'/\langle\mathrm{SO}(1,2)\cap \mathrm{SO}(3),
i_{x_0} \rangle\cong S^2$ collapses to a segment
$L_0$ when $t\searrow0$, so that the result is 
homeomorphic to $\mathbb{R}^3$.
\end{proof}

\subsection{The elliptic and Fuchsian Locus of 
$\mathcal{X}_2(\Gamma,\mathrm{Isom}(X))$}

Let $\mathcal{E}_2$ denote the elliptic locus in 
$\mathcal{X}_2(\Gamma,\mathrm{Isom}(X))$, that is,
\(
\mathcal{E}_2 = \mathcal{E} \cap \mathcal{X}_2(\Gamma,\mathrm{Isom}(X)),
\)
and define $\mathcal{F}_{\mathrm{I},2}$ and $\mathcal{F}_{\mathrm{II},2}$ analogously. By construction, all these sets are closed subsets of 
$\mathcal{X}_2(\Gamma,\mathrm{Isom}(X))$.
\begin{Proposition}
We have the following homeomorphisms
$$
\mathcal{E}_2\cong[0,1],\qquad \mathcal{F}_{\mathrm{I},2}\cong [0,+\infty), 
\quad\textrm{ and }\quad
\mathcal{F}_{\mathrm{II},2}\cong [0,+\infty) \sqcup [0,+\infty).
$$
In addition 
% $\mathcal{E}\cap \mathcal{F}_{\mathrm{II},2}=\partial \mathcal{E}_2=\partial \mathcal{F}_{\mathrm{II},2}$
% and $\mathcal{E}\cap \mathcal{F}_{\mathrm{I},2}
% =\mathcal{F}_{\mathrm{I},2}\cap\mathcal{F}_{\mathrm{II},2}$ is one of the end-points of either 
those components intersect at the end-points of the intervals or half-intervals as in Figure~\ref{Figure:F2}.
\end{Proposition} 

\begin{figure}[ht]
\begin{center}
  \begin{tikzpicture}[line join = round, line cap = round, scale=1]
    \draw[very thick]  (0,0)--(2,0)--(5,0);
    \draw[very thick]  (-2,-1)--(0,0)--(-2,1);
    \draw[black, fill=black] (0,0) circle[radius=.06];
    \draw[black, fill=black] (2,0) circle[radius=.06];
    \node at (1,.3) {$\mathcal{E}_2$};
  \node at (-2.5,1) {$\mathcal{F}_{\mathrm{I},2}$};
  \node at (-2.5,-1) {$\mathcal{F}_{\mathrm{II},2}$};
  \node at (4,.3) {$\mathcal{F}_{\mathrm{II},2}$};
  \end{tikzpicture}
  \end{center}
  \caption{The Fuchsian locus of $\mathcal{X}_2(\Gamma,\mathrm{Isom}(X))$.}
  \label{Figure:F2}
  \end{figure}

\begin{Remark}
    Along the proof we check that representations in 
    $\mathcal{F}_{\mathrm{I},2}$
    (and in the adjacent component of $\mathcal{F}_{\mathrm{II},2}$)
    do not preserve the orientation of the 
    invariant Fuchsian plane, while representations in 
    the other component of 
    $\mathcal{F}_{\mathrm{II},2}$ preserve it.
\end{Remark}

\begin{proof}
    The elliptic locus $\mathcal{E}_2 $ corresponds to 
    $L_0'\big/\mathrm{SO}(3)\cap\mathrm{SO}(1,2)
$ in the proof of 
Proposition~\ref{Proposition:TopologyX_2}, hence it is 
homeomorphic to a closed interval.

We now describe the intersection of the Fuchsian locus 
with the space
\(
L_t'\big/\bigl(\mathrm{SO}(3)\cap\mathrm{SO}
(1,2)\bigr)\cong S^2
\)
for each parameter value \(t>0\). Since this 
intersection is in fact independent of \(t\), every 
connected component of the Fuchsian locus 
is a closed subset and is 
homeomorphic to a half-interval. The attaching point to 
\(\mathcal E_2\) is determined by the value of \(\lvert 
x\rvert\).

An invariant $\mathbb{H}^2$ is perpendicular to both 
$\mathrm{Fix}(\rho(a))$ and  $\mathrm{Fix}(\rho(b))$, so 
if $m\in Q$ is as in Lemma~\ref{Lemma:parametersL'2},
then $m$ is tangent to the invariant $\mathbb{H}^2$.
By Lemma~\ref{Lem:FuchsianRegular} $\mathrm{type}(m)=\mathrm{center}(\sigma_{\mathrm{mod}})$, thus 
using the parameterization in the proof of 
Lemma~\ref{Proposition:Fibre} it corresponds to $\alpha=0$ and $\alpha=\pi/3$.

When $\alpha=0$, 
$$m=\frac{1}{\sqrt{2}}
    \begin{pmatrix}
        0 & 0 & 0 \\
        0 & 1 & 0 \\
        0 & 0 & -1
    \end{pmatrix}
% \mathrm{diag}(0,1,-1)
.
$$  By Lemma~\ref{Lemma:fewplanes} $m$ is tangent to a unique Fuchsian plane of type I and two of type II. We look at each case:
\begin{enumerate}[(a)]
    \item The Fuchsian plane of type I is tangent to 
the plane spanned by $m$ and 
$$
\begin{pmatrix}
    0 & 0 & 0 \\
    0 & 0 & 1 \\
    0 & 1 & 0
\end{pmatrix}
$$
which is invariant under the rotations 
$R_{\pm \frac{2\pi}{3}}=\mathrm{Rot}( \pm (1 \, 0\, 0)^{\mathrm{T}})
$. Those rotations are unique, by
Lemma~\ref{Lemma:uniquerotations}: 
thus we get the coordinate   $|x|=1$. Namely, in 
Figure~\ref{Figure:Ra}, it corresponds to the point
with coordinates $(|x|, \alpha)=(1,0)$. Notice that this 
Fuchsian plane is perpendicular to 
$\mathrm{Fix}(\rho(a))$, so the orientation of the invariant Fuchsian plane is preserved.

\item The Fuchsian planes of type II are  those tangent to $m$ and one of the following
$$
\begin{pmatrix}
    0 & 1 & 1 \\
    1 & 0 & 0 \\
    1 & 0 & 0
\end{pmatrix}
\qquad \textrm{ or }\qquad
\begin{pmatrix}
    0 & 1 & -1 \\
    1 & 0 & 0 \\
    -1 & 0 & 0
\end{pmatrix}
$$
The unique rotations that preserve the plane 
(provided by Lemma~\ref{Lemma:uniquerotations})
are respectively 
$\mathrm{Rot}(\pm\frac{1}{\sqrt{2}} (0\ 1\ 1)^{\mathrm{T}}) $  and
$\mathrm{Rot}(\pm\frac{1}{\sqrt{2}} (0\  1\  -1)^{\mathrm{T}}) $.
In both cases $x=0$ (they are related by the 
stabilizer of $m$). Thus this yields the point with 
coordinates $(|x|,\alpha)=(0,0)$. In this case the 
invariant Fuchsian plane and $\mathrm{Fix}(\rho(a))$
intersect at a geodesic, so the orientation of the 
Fuchsian plane is reversed by $\rho(a)$.

\end{enumerate}

When $\alpha=\frac{\pi}{3}$ then 
$$m=\frac{1}{\sqrt{2}}
    \begin{pmatrix}
        1 & 0 & 0 \\
        0 & 0 & 0 \\
        0 & 0 & -1
    \end{pmatrix}
% \mathrm{diag}(1, 0, -1)
.
$$
Again $m$ is tangent to a Fuchsian plane of type I and two of type II, and we distinguish two cases:

\begin{enumerate}[(a)]
    \item The type I plane is the plane tangent to $m$ and
$$
\begin{pmatrix}
    0 & 0 & 1 \\
    0 & 0 & 0 \\
    1 & 0 & 0
\end{pmatrix}
$$
The unique rotations that preserve this plane are 
$\mathrm{Rot}( \pm (0\, 1\, 0)^{\mathrm{T}})$, hence $x=0$, 
thus the point
has coordinates 
$(|x|,\alpha)=(0,\pi/3)$.
The invariant Fuchsian plane and $\mathrm{Fix}(\rho(a))$
meet at a geodesic, so the orientation of the Fuchsian plane is reversed by $\rho(a)$. 
\item 
The type II planes are tangent to the span of $m=\mathrm{diag}(1, 0, -1)$ and one of the following
$$
\begin{pmatrix}
    0 & 1 & 0 \\
    1 & 0 & 1 \\
    0 & 1 & 0
\end{pmatrix}
\qquad \textrm{ or }\qquad
\begin{pmatrix}
    0 & 1 & 0 \\
    1 & 0 & -1 \\
    0 & -1 & 0
\end{pmatrix}
$$
In this case the Fuchsian plane does not intersect $\mathrm{Fix}(\rho(a))$ perpendicularly, so the representation is not Fuchsian. 
\end{enumerate}

Summarizing, we have obtained the points with coordinates 
$(|x|,\alpha)$ equal to $(1,0)$, $(0,0)$
and $(0,\frac{\pi}{3})$. As $\vert x|$ describes the attaching map to the segment $L_0\big/\mathrm{SO}(3)\cong  [0,1]$, the proposition follows.
\end{proof}

\section{Finding Anosov representations}
\label{Sec:Anosov}

  In the present section we prove Theorem~\ref{Thm:Anosov}. Starting from the presentation
\[
\Gamma= \mathrm{PSL}_2(\mathbb{Z}) = \langle a, b \mid a^2 = b^3 = 1 \rangle,
\]
for  
$\rho \in \hom(\Gamma, \mathrm{Isom}
(X))$ a representation whose orbit by conjugation is closed, 
we define \(x\) to be the point of $\mathrm{Fix}(\rho(a))$
that minimizes the distance
to the singular geodesic $\gamma=\mathrm{Fix}(\rho(b))$,
namely $$
d(x,\gamma)=d(\mathrm{Fix}(\rho(a)),\gamma).
$$

Let \(\mathbb{H}^2\subset X\) 
denote a $\rho(b)$-invariant 
totally geodesic hyperbolic plane (perpendicular to $\gamma$), and consider the nearest–point projection
\[
\pi \colon X \to \mathbb H^2.
\]
We 
use the distances
\[
s = d(x, \mathbb{H}^2) \quad \text{and} \quad t = d(\gamma, \pi(x)).
\]
Theorem~\ref{Thm:Anosov} asserts that, for a given \(C \geq 0\), there exists a threshold \(T = T(C)\) such that for all 
$s\leq C$ and 
\(t \geq 
T\), the representation \(\rho\) is Anosov.

Recall that $X$ has non-positive sectional curvature, a fact we shall exploit throughout this  section.

% \begin{Theorem}[Theorem~\ref{Thm:Anosov}]
% Given $s>0$, for $t>0$ sufficiently large, depending on $s$,  representations with coordinates 
% $(s,\theta,t)$ are Anosov.
% \end{Theorem}

\subsection{Triangles with small angles}

In this subsection
we establish some  geometric lemmas, 
preliminary to the proof of Theorem~\ref{Thm:Anosov}.

Let $\rho\in \hom(\Gamma,\mathrm{Isom}(X))$,
  and let $x$, $\gamma$, $\mathbb{H}^2$, $s$ and $t$ be as above. Besides $x$, we consider the following points 
 that minimize the distance to $\gamma$ in the respective
 fixed point sets:
\begin{align*}
    &x\in \mathrm{Fix} (\rho(a)), \\
    &y=  \rho(b)(x)\in \mathrm{Fix} (\rho(bab^2)), \\
    &z=  \rho(b^2) (x)\in \mathrm{Fix} (\rho(b^2ab)).
\end{align*}

\begin{Lemma}
    \label{lem:smallangles}
    For any $\varepsilon>0$, if $t$ is sufficiently large (depending on $C$ and $\varepsilon$, where $s\leq C$),
    then the angles within the triangle $xyz$ are less than $\varepsilon$.
\end{Lemma}

\begin{proof}
To simplify, assume that the sectional curvature 
of $\mathbb H^2$ is $-1$, as there are only 
two congruence classes of such $\mathbb{H}^2$.

Let $ \bar x=\pi(x)$,  
$ \bar y=\pi(y)$, and $ \bar z=\pi(z)$ denote the 
respective
projections of $x$, $y $ and $z$  to $\mathbb{H}^2$, 
they define an equilateral triangle in 
$\mathbb H^2$ 
centered at 
$$
c=\gamma\cap \mathbb H^2.
$$
By  construction
% \begin{align*}
% &
\begin{equation*}
d(c,\bar{x})=d(c,\bar{y})=d(c,\bar{z})=t.    
\end{equation*}
% \qquad\qquad \textrm{ and }\\\
% &\angle_c(\bar{x},\bar{y})=\angle_c(\bar{x},\bar{z})=\angle_c(\bar{y},\bar{z})=\tfrac{3\pi}{2}. &
% \end{align*}
By elementary hyperbolic geometry, as $t\to+\infty$, $\angle_{\bar x}(\bar y,\bar z)\to 0$.

\begin{figure}[ht]
\begin{center}
  \begin{tikzpicture}[line join = round, line cap = round, scale=.7]

   % \fill[gray!30!white, opacity=.4]  (0,0)--(2.5,0)  arc[start angle=0, end angle=60, radius=2.5]--(0,0);
  \draw[dgray, thick] (0,0)--(0,2);
  \fill[black]  (0,2) circle[radius=.07];
     \draw[dgray,  thick, rotate=120] (0,0)--(0,2);
     \draw[dgray,  thick, rotate=-120] (0,0)--(0,2);
    \draw[thick]  (0,2) arc[start angle = 195, end angle = 225, radius= 6.6];
  \fill[black]  (0,-.78) circle[radius=.07];
    \begin{scope}[rotate=120]
     \draw[thick]  (0,2) arc[start angle = 195, end angle = 225, radius= 6.6];
       \fill[black]  (0,2) circle[radius=.07];
       \node at (0,2.3) {$\bar z$};
       \fill[black]  (0,-.78) circle[radius=.07];
       \node at (-.25,-1.4) {$ m_{\bar x\bar y}$};
     \end{scope}
    \begin{scope}[rotate=-120]
     \draw[thick]  (0,2) arc[start angle = 195, end angle = 225, radius= 6.6];
       \fill[black]  (0,2) circle[radius=.07];
       \node at (0,2.3) {$\bar x$};
        % \fill[black]  (0,-.78) circle[radius=.07];
     \end{scope}
     \fill[black]  (0,0) circle[radius=.07];
      \node at (0.2,.15) {$c$};
      \node[dgray] at (0.7,-.13) {$t$};
      \node at (0,2.4) {$\bar y$};
       \node at (0.3,-1.2) {$m_{\bar x\bar z}$};
%       \node at (0.5,2.2) {${\lambda_1=\lambda_2}$};
%       \node at (1.8,-2) {${\lambda_1=\lambda_3}$};
%     \draw[red, very thick] (2.6, 0)  arc[start angle=0, end angle=60, radius=2.6];      
%     \node[red] at (2.7, 1.2) {$\sigma$};
%     \node at (1.4,0.7) {$V(x,\sigma)$};
% % %   \node[black] at (1.5,1) {${\Delta}$};
% % %   \node at (6,2) {${\lambda_1+\lambda_2+\lambda_3=0}$};
  \end{tikzpicture}
  \end{center}
  \caption{The triangle $\bar{x}\bar{y}\bar{z}$ in the plane $\{ h_0\}\times\mathbb H^2$.}
  \label{Figure:HyperbolicTriangle}
  \end{figure}

To prove the lemma, we estimate $\angle_x(y,z)$ from three  inequalities:

\begin{itemize}
    \item  Firstly consider the midpoints 
    $m_{\bar x \bar y}=\mathrm{mid}(\bar x, \bar y)$ and 
    $m_{\bar x \bar z}=\mathrm{mid}(\bar x, \bar z)$. 
By thinness and 
symmetry of the hyperbolic triangle $\bar x\bar y\bar z$,
Figure~\ref{Figure:HyperbolicTriangle},
$$
d(m_{\bar x \bar y},m_{\bar x \bar z} )
<d(m_{\bar x \bar y},c)+d(c,m_{\bar x \bar z} )<\log 3.
$$
On the other hand:
$$
d(x, m_{\bar x \bar y})\geq d(\bar x, c)-d(x,\bar x)-d(c,m_{\bar x \bar y})
\geq t-s-\tfrac12\log 3.
$$
% d(x, m_{\bar x \bar z})&\geq d(c, \bar x)-d(x,\bar x)-d(m_{c,\bar x \bar z})
% \geq t-s-\tfrac12\log 3 
% \end{align*}
As $s\leq C$, both inequalities yield:
\begin{equation}
    \label{eqn:step1}
    \angle_x(m_{\bar x \bar y}, m_{\bar x \bar z}) <\varepsilon,
    \qquad\text{ for }t\textrm{ large enough.}
\end{equation}

\item
Secondly, for $t$ sufficiently large we have 
$\angle_{m_{\bar x \bar y}}(x,\bar x) < \varepsilon/2$, because 
$d(x, \bar x)=s$ remains bounded whilst $d( m_{\bar x \bar y},x)$ and
$d( m_{\bar x \bar y}, \bar x)$ are arbitrarily large with $t$. 
By symmetry, we also have
$\angle_{m_{\bar x \bar y}}(y,\bar y) < 
\varepsilon/2$. 
From the previous two inequalities, we get that 
$\angle_{m_{\bar x \bar y}}(x,y)> \pi-\varepsilon$, because 
 $m_{\bar x \bar y}$ lies in the segment $\bar x\bar y$.
 % Figure~\ref{Figure:Lines}.
By considering the
triangle $m_{\bar x \bar y}xy$, we get
\begin{equation}
    \label{eqn:step2}
    \angle_x(m_{\bar x \bar y}, y) <\varepsilon
    \qquad\text{ for }t\textrm{ large enough,}
\end{equation}
because in non-positive curvature the addition of angles of a triangle 
is  at most $\pi$.

\item
Finally, by symmetry, the same argument as in \eqref{eqn:step2} yields
\begin{equation}
    \label{eqn:step3}
    \angle_x(m_{\bar x \bar z}, z) <\varepsilon
    \qquad\text{ for }t\textrm{ large enough.}
\end{equation}
\end{itemize}
The three inequalities \eqref{eqn:step1}, \eqref{eqn:step2}, 
and \eqref{eqn:step3} imply the lemma (replace $\varepsilon$
by $\varepsilon/3$).
\end{proof}

\begin{Lemma}
    \label{Lemma:almostperp}
    Assume $\mathrm{Fix}(\rho(a))$ is not a single point
    (namely $\mathrm{Fix}(\rho(a))$ is  either
    a Fuchsian plane of type II or a parallel set).
        For $\varepsilon >0$, if $t$
    is sufficiently large (depending on $C$ and $\varepsilon$, 
    where $s\leq C$),
    then the angle at $x$ between $xy$ and $\mathrm{Fix}(\rho(a))$ is at
    least $\frac{\pi}{2}-\varepsilon$.
\end{Lemma}

\begin{proof}
Let \(p \in \gamma\) be the point such that the segment 
${xp}$ is orthogonal to both   
\(\mathrm{Fix}(\rho(a))\) and \(\gamma\). 
Our goal is to establish that 
\begin{equation}
    \label{eqn:claimed}
\angle_x(p,y) < \varepsilon  ,  
\end{equation}
because, since $xp$ is perpendicular to 
\(\mathrm{Fix}(\rho(a))\), \eqref{eqn:claimed} 
yields the lemma.
% Once we have \eqref{eqn:claimed} we conclude because 
% for any 
% $p\in \mathrm{Fix}(\rho(a))$, $p\neq x$, we have
% $$
% \angle_x(p,y)\geq
% \angle_x(p,w)-\angle_x(w,y)=\frac{\pi}{2}-\angle_x(w,y)
% >\frac{\pi}{2}-\varepsilon.
% $$

Analogously 
to the argument in Lemma~\ref{lem:smallangles}, 
we shall derive \eqref{eqn:claimed} as a consequence 
of three auxiliary inequalities.
        We recover the notation of Lemma~\ref{lem:smallangles}:
        let  $ \bar x=\pi(x)$,  
$ \bar y=\pi(y)$, and $ \bar z=\pi(z)$ denote the projections
to $\mathbb{H}^2$, they 
 define an equilateral triangle 
centered at 
$
c=\gamma\cap \mathbb H^2
$, so that \(
d(c,\bar{x})=d(c,\bar{y})=d(c,\bar{z})=t    \).

For the first inequality, observe that
\[
d(p,c) \leq d(x,\bar x) = s,
\]
since both segments \(px\) and \(c\bar x\) are perpendicular to 
\(\gamma\), and the nearest point projection to 
$\gamma$ is distance non-increasing. Moreover, as
\[
d(x,c)\geq d(\bar x,c)-d(x,\bar x)=t-s,
\]
we obtain
\begin{equation}
\label{eqn:newstep1}
    \angle_x(p,c)<\varepsilon
    \qquad\text{for } t \text{ sufficiently large.}
\end{equation}

Next, recall that \(m_{\bar x\bar y}\) denotes the midpoint 
of \(\bar x\) and \(\bar y\). Since
\(
d(c, m_{\bar x\bar y})<\frac{1}{2}\log 3,
\) and $d(x,c)\geq t-s$,
it follows that
\begin{equation}
\label{eqn:newstep2}
    \angle_x(c,m_{\bar x\bar y})<\varepsilon
    \qquad\text{for } t \text{ sufficiently large.}
\end{equation}

The claimed inequality~\eqref{eqn:claimed} 
then follows from  
inequalities~\eqref{eqn:newstep1} 
and~\eqref{eqn:newstep2}, together with inequality~\eqref{eqn:step2}
in the proof of Lemma~\ref{lem:smallangles} (namely \(\angle_x(m_{\bar x\bar y},y)<\varepsilon\) for large \(t\)).
\end{proof}

By Lemma~\ref{Lemma:almostperp} when 
$\mathrm{Fix}(\rho(a))$ is larger than a point (and trivially when it is just a point):

\begin{Corollary}
    \label{Corollary:xyay}
     For $\varepsilon >0$, if $t$
    is sufficiently large (depending on $C$ and $\varepsilon$, where $s\leq C$), then 
    $$
    \angle_x(y,\rho(a)y)\geq \pi-2\varepsilon.
    $$
\end{Corollary}

For the next lemma we need a couple of definitions:
% (see Remark~\ref{Remark:TypeRegularity}):

\begin{Definition}
\label{Definition:regular}
Given $x\neq y\in X$:
\begin{enumerate}[(a)]
    \item The \emph{type of the segment}
$xy$ is $\mathrm{type}(xy)=\mathrm{type}(r(+\infty))\in \sigma_{\mathrm{mod}}$,
where $r$ is the ray such that $r(0)=x$ and $y\in\mathrm{Image}(r)$.
    \item The segment $xy$ is \emph{regular} 
if $\mathrm{type}(xy)\in\mathrm{int}(\sigma_{\mathrm{mod}})$.
\item Given $\Theta\subset \mathrm{int}(\sigma_{\mathrm{mod}})$ a relatively compact  subinterval 
in $\mathrm{int}(\sigma_{\mathrm{mod}})$
that is $\iota$-invariant, the segment $xy$ is  $\Theta$-\emph{regular} if
$\mathrm{type}(xy)\in\Theta$.
\end{enumerate}
\end{Definition}

Notice that $\mathrm{type}(yx)=\iota( \mathrm{type}(xy) ) $.

\begin{Lemma} 
\label{Lem:regular}
For sufficiently large $t>0$ (depending on $C$, where $s\leq C$), 
the type of $xy$ 
is arbitrarily close
to  the center   of the spherical Weyl chamber
$\zeta=\mathrm{center}(\sigma_{\mathrm{mod}})$.

In particular, given a subinterval 
$\Theta\subset\mathrm{int} (\sigma_{mod})$ with 
$\iota(\Theta)=\Theta$ and 
$\mathrm{int}(\Theta)\neq\emptyset$, 
for sufficiently large $t$ (depending on $C$ and $\Theta$)
$xy$ is $\Theta$-regular. 
\end{Lemma}

\begin{proof}
% Here we apply continuity. Let $\zeta\in \sigma_{mod}$ be the fixed point of $\iota$,
% so that $\zeta\in\mathrm{int}(\Theta)$.
% Recall the 
% isometry  $\mathbb R\times\mathbb H^2\cong P(\gamma)$, we claim that 
% for any $h_0\in\mathbb R$, 
% \begin{equation}
% \label{eqn:constanttype}
% \mathrm{type} (\partial_\infty (\{h_0\}\times \mathbb H^2) )=\{\zeta\}.
% \end{equation}
% To prove \eqref{eqn:constanttype}, notice that using the description in 
% Subsection~\ref{SubSection:Parallel}
% every asymptotic point in 
% $\partial_\infty (\{h_0\}\times \mathbb H^2)$
% can be seen as the limit when $t\to+\infty$ of the geodesic
% $$
% t\mapsto A_t=
% \begin{pmatrix}
% e^{\tau_0} & & \\
%  & e^{-\tau_0/2} & \\
%  & & e^{-\tau_0/2}
% \end{pmatrix}
% \begin{pmatrix}
% 1 & & \\
%  & \cos\theta  & -\sin\theta\\
%  & \sin\theta & \cos\theta
% \end{pmatrix}
% \begin{pmatrix}
% 1 & & \\
%  & e^{t} & \\
%  & & e^{-t}
% \end{pmatrix}.
% $$
% The singular values of $A_t$, namely the square root of the eigenvalues of $\left(A_t\right)^{\mathrm T} A_t$, are
% $\{ e^{t-\tau_0/2}, e^{\tau_0}, e^{-t-\tau_0/2} \}$, which implies \eqref{eqn:constanttype}.

By Lemma~\ref{Lem:FuchsianRegular},  
% Once we have checked \eqref{eqn:constanttype}, 
the argument follows from continuity
of the type map. Use  the fact that, for bounded $s$,  both triangles
$xyz$ and $\bar{x}\bar{y}\bar{z}$ converge to the same ideal triangle 
(for convergence in the pointed Hausdorff topology, e.g.~Hausdorff 
convergence of the intersection with any metric ball).
The second assertion is clear from the first one  as $\zeta=\mathrm{center}(\sigma_{\mathrm{mod}})\in\mathrm{int}(\Theta)$.
\end{proof}

\subsection{Straight and spaced sequences}

To prove Theorem~\ref{Thm:Anosov} we apply  
\cite[Theorem~3.18]{KLP25}. Before stating it, 
we introduce 
some further definitions.

Let $p,q,q'\in X$ be three different points such that $pq$ and $pq'$ are regular
(Definition~\ref{Definition:regular}).
There exist   Weyl sectors $V(p,\sigma)$ and 
$V(p,\sigma')$ such that
$q\in V(p,\sigma)$ and $q'\in V(p,\sigma')$, 
that are unique by regularity.
Choose $\zeta\in \mathrm{int}(\sigma_{\mathrm{mod}})$
and let $r$ and $r'$ be rays starting at $p$ with image in 
$V(p,\sigma)$ and $V(p,\sigma')$ respectively, and with 
$\mathrm{type}(r(+\infty))=\mathrm{type}(r'(+\infty))=\zeta$,
Figure~\ref{Figure:SectorsAngle}.

\begin{Definition} For  $\zeta\in \mathrm{int}(\sigma_{\mathrm{mod}})$, 
    the $\zeta$-angle at $p$ of the segments $pq$ and $pq'$ is
    $$
    \angle^\zeta_p(q,q')=\angle_p(r,r') .
    $$
\end{Definition}

This must be understood as a way to measure angles between Euclidean 
Weyl sectors (see Figure~\ref{Figure:SectorsAngle}).

\begin{figure}[ht]
\begin{center}
  \begin{tikzpicture}[line join = round, line cap = round, scale=.8]
\begin{scope}[yscale=0.7]
\draw[white, fill=gray!50!white, opacity=0.3] (3,2)--(0,0)--(4,1/2)  to[out=105,in=-30] (3,2);
\draw (0,0)--(4,1/2);   
\draw (0,0)--(3,2);   
\draw[ thick] (0,0)--(3.5*1.10,5/4*1.10);  
\end{scope}
\draw[ultra thick] (1.4,0.35) to[out=105, in=-45] (1,1.1);
% \draw (4,1/2)  to[out=105,in=-30] (3,2);
\begin{scope}[scale=.8, rotate=10]
    \draw[white, fill=white] (3.5,1.75)--(0,0)--(2.5,3) to[out=-30,in= 105] (3.5,1.75); 
    \draw[white, fill=gray!50!white, opacity=0.3] (3.5,1.75)--(0,0)--(2.5,3) to[out=-30,in= 105] (3.5,1.75); 
    \draw (0,0)--(3.5,1.75);   
    \draw (0,0)--(2.5,3);  
    \draw[ thick] (0,0)--(3*1.10,4.75/2*1.10); 
\end{scope}
% \draw[thick, opacity=.2] (0,0)--(3,2);
\draw[very thick, dotted] (1.4,0.35) to[out=105, in=-45] (1,1.1);
\fill[black]  (0,0) circle[radius=.07];
\fill[black]  (3,.45) circle[radius=.07];
\fill[black]  (1.4,2.1) circle[radius=.07];
\node at (-.4,-.1) {$p$};
\node at (3,0) {$q$};
\node at (1.0,2.2) {$q'$};
\node at (4.1, 1.05) {$r$};
\node at (2.5, 2.8) {$r'$};
\node at (5,0.5) {$V(p,\sigma)$};
\node at (1.3,3.1) {$V(p,\sigma')$};
  \end{tikzpicture}
  \end{center}
  \caption{The $\zeta$-angle $\angle^\zeta_p(q,q')=\angle_p(r,r')$.}
  \label{Figure:SectorsAngle}
  \end{figure}

% We require some notation for the next definition. Let  $\zeta=\mathrm{center}(\sigma_{\mathrm{mod}})$
% denote the center of the spherical Weyl chamber, in particular $\iota(\zeta)=\zeta$. 
% Let $\Theta\subset\mathrm{int}(\sigma_{\mathrm{mod}})$ be a relatively compact 
% interval with $\iota(\Theta)=\Theta$ and $\mathrm{int}(\Theta)\neq \emptyset$.  

\begin{Definition}[\cite{KLP25}]
\label{Def:straight}
Set $\zeta=\mathrm{center}(\sigma_{\mathrm{mod}})$ 
and  let $\Theta\subset\mathrm{int}(\sigma_{\mathrm{mod}})$ be a relatively compact 
interval with $\iota(\Theta)=\Theta$ and $\mathrm{int}(\Theta)\neq \emptyset$.  
For $\varepsilon> 0$ and 
$l>0$,
a sequence $(x_n)_{n\in\mathbb Z}$ in $X$ is:
\begin{itemize}
    \item  $(\Theta,\varepsilon)$-\emph{straight}
if for every $n\in\mathbb Z$ the segments $x_nx_{n+1}$ are $\Theta$-regular and 
$$
\angle^\zeta_{x_n}(x_{n-1}, x_{n+1})\geq \pi-\varepsilon.
$$
    \item $l$-\emph{spaced} if $d(x_n, x_{n+1})\geq l$ for every $n\in\mathbb Z$.
\end{itemize}
\end{Definition}

The following is Theorem~3.18 from  \cite{KLP25} and it is the key tool to prove
Theorem~\ref{Thm:Anosov}. We use the notation 
$$
 V(\bar x_n,\sigma_{\pm}^{\Theta'})=\{\bar x_n\}\cup 
 \{y\in V(\bar x_n,\sigma_{\pm})\mid y\neq \bar x_n\textrm{ and }
 \mathrm{type}(\bar x_ny)\in \Theta' \}.
 $$ 
 Also $F(\sigma_-,\sigma_+)$ denotes the unique maximal flat asymptotic
to two opposite chambers $\sigma_-$ and $\sigma_+$.
Finally, 
$\Theta,\Theta'\subset \mathrm{int} (\sigma_{\mathrm{mod}})$ 
denote
$\iota$-symmetric,
relatively compact intervals, with nonempty interior.

\begin{Theorem}[\cite{KLP25}] 
\label{Thm:LocalMorse}
For  $\Theta,\Theta'$ as above and $\delta>0$,
% For $\Theta,\Theta'\subset \mathrm{int} (\sigma_{\mathrm{mod}})$ 
% $\iota$-symmetric,
% relatively compact intervals, and with nonempty interior, 
there are $l>0$ and $\varepsilon>0$ such that the 
following holds.

If $(x_n)_{n\in\mathbb Z}$ is a  
$(\Theta,\varepsilon)$-straight and $l$-spaced
sequence in $X$,  then $x_n$ is $\delta$-close to a maximal 
flat $F(\sigma_-,\sigma_+)$.
Furthermore the sequence moves from $\sigma_-$ to $\sigma_+$ in the sense that
$\bar x_n$, the nearest point 
projection of $x_n$ to $F(\sigma_-,\sigma_+)$, 
satisfies 
$$
\bar x_{n\pm m}\in V(\bar x_n,\sigma_{\pm}^{\Theta'})\qquad \textrm{for every } m, n\in\mathbb Z, 
\ m\geq 1
$$   
\end{Theorem}

\begin{figure}[ht]
\begin{center}
  \begin{tikzpicture}[line join = round, line cap = round, scale=1] 
   \draw[black]  (-3,1)--(-1.5,.7)--(0,.8)--(1.5,.6)--(3,.9);
   \draw[black]  (-3,-.5)--(-1.55,-.1)--(0,0)--(1.5,-.12)--(3,-.5);
  \draw[black, fill=black] (0,0) circle[radius=.03];
  \draw[black, fill=black] (1.5,-0.12) circle[radius=.03];
  \draw[black, fill=black] (-1.55,-0.1) circle[radius=.03];
  \draw[black, fill=black] (3,-0.5) circle[radius=.03];
  \draw[black, fill=black] (-3,-0.5) circle[radius=.03];
  \draw[black, fill=black] (0,.8) circle[radius=.03];
  \draw[black, fill=black] (1.5,.6) circle[radius=.03];
  \draw[black, fill=black] (-1.5,.7) circle[radius=.03];
  \draw[black, fill=black] (3,0.9) circle[radius=.03];
  \draw[black, fill=black] (-3,1) circle[radius=.03];
  \draw[gray, opacity=0.6]  (-3,1/2)--(-4,-1.2)--(3,-1.2)--(4.15,1/2);
  \begin{scope}[scale=1.4]
        \draw[black, fill=gray!30!white, opacity=0.3]  (2.10,-.75)--(0,0)--(2.85,.35);
        \draw[black, fill=gray!30!white, opacity=0.3]  (-2.75,-.75)--(0,0)--(-2,.4);      
  \end{scope}
  \node at (0.1,.2) {$\bar x_n$};  
  \node at (2.0,-.03) {$\bar x_{n+1}$};
  \node at (-1.2,.1) {$\bar x_{n-1}$};
  \node at (0.1,1) {$x_n$};  
  \node at (3.5,-0.55) {$\bar x_{n+2}$};
  \node at (-2.75,-.8) {$\bar x_{n-2}$};
  \node at (1.9,.9) {$ x_{n+1}$};
  \node at (-1.3,.85) {$ x_{n-1}$};
  \node at (3.5,1) {$ x_{n+2}$};
  \node at (-2.7,1.2) {$x_{n-2}$};
  \node at (4.5,-.15) {$ V(\bar x_n,\sigma_{+}^{\Theta'})$};
  \node at (-4.15,0) {$ V(\bar x_n,\sigma_{-}^{\Theta'})$};
  \node at (6,.75) {$ F(\sigma_-,\sigma_+) $};
  %
    % \draw[white, fill=gray!30!white, opacity=0.3]  (1.5,1)--(-1,0)--(1.5,-1);
    % \draw[white, fill=gray!30!white, opacity=0.3]  (1.5,1) to[out=-30, in=30](1.5,-1);
    % \draw[white, fill=gray!30!white, opacity=0.3]  (-1.5,1)--(1,0)--(-1.5,-1);
    % \draw[white, fill=gray!30!white, opacity=0.3]  (-1.5,1) to[out=240, in=120] (-1.5,-1);
    % \draw[ ultra thick, black, fill=gray!50!white, opacity=0.6] (-1,0)--(0,0.4)--(1,0)--(0,-.4)--(-1,0);
    % \draw[thick, black]  (-1,0)--(1.5,1);
    % \draw[thick, black]  (-1,0)--(1.5,-1);
    % \draw[thick, black]  (1,0)--(-1.5,1);
    % \draw[thick, black]  (1,0)--(-1.5,-1);
    %  \node at (1.3,0) {$x_+$};
    %  \node at (-1.3,0) {$x_-$};
    %  \node at (2.5,.5) {$V(x_-,\sigma_+)$};
    %    \node at (-2.5,.5) {$V(x_+,\sigma_-)$};
  \end{tikzpicture}
  \end{center}
  \caption{The maximal flat $F(\sigma_-,\sigma_+)$ as in Theorem~\ref{Thm:LocalMorse}}
  \label{Figure:Diamonds}
  \end{figure}

This theorem is a local version of the Morse lemma 
(see \cite{KLPMorseLemma}). It requires
uniform regularity and a control on the angles $\angle^{\zeta}$
between sectors (because of maximal flats, as 
there is no Morse lemma in the Euclidean plane).

\subsection{The Cayley graph of subgroups of $\Gamma$}

    Recall that $x\in \mathrm{Fix}(\rho(a))$, $y=\rho(b) x$ and $z=\rho(b)y$. Using the isomorphism $\mathrm{PSL}_2(\mathbb Z)\cong \mathbb Z_2*\mathbb Z_3$, consider the finite-index subgroups
    \begin{align*}
        \Gamma_2& = \ker(\mathrm{PSL}_2(\mathbb Z) \twoheadrightarrow \mathbb Z_2 ),\\
        \Gamma_3& = \ker(\mathrm{PSL}_2(\mathbb Z) \twoheadrightarrow \mathbb Z_3 ),\\
        \Gamma_6&=\Gamma_2\cap \Gamma_3\cong F_2.
    \end{align*}
We first describe the Cayley graph of $\Gamma_6$ embedded in $\mathbb{H}^2$.

\begin{figure}[ht]
\begin{center}
  \begin{tikzpicture}[line join = round, line cap = round, scale=.7] 
       \begin{scope}[opacity=.7]
            \draw [dgray] (0,-6) -- (0,6);
        \draw [dgray] (0,0) circle (6);
        \foreach \i in {0,1,2,3} {%
            \draw [ dgray] (90*\i:6) arc (270+90*\i:180+90*\i:6);}
         \foreach \i in {0,1,...,7} {%
             \draw [dgray] (45*\i:6) arc (270+45*\i:135+45*\i:3.3*.75);}
        \foreach \i in {0,1,...,15} {%
             \draw [dgray] (22.5*\i:6) arc (270+22.5*\i:112.5+22.5*\i:1.6*.75);}
         \foreach \i in {0,1,...,31} {%
             \draw [thin] (11.25*\i:6) arc (270+11.25*\i:101.25+11.25*\i:0.8*.75);}
        % \foreach \i in {0,1,...,63} {%
        %     \draw [ultra thin] (5.625*\i:8) arc (270+5.625*\i:95.625+5.625*\i:0.4);}
        \end{scope}
  \node[fill = white]  at (.5,0) {$\underline x$};
  \node[fill = white]  at (2.7,-1.4) {$\underline y=b\underline x$};
  \node[fill = white]  at (2.7,1.4) {$\underline z=b^2\underline x$};
  \node[fill = white]  at (-0.4,2.1) {$a\underline y=aba\underline x$};
  \node[fill = white]  at (-0.4,-2) {$a\underline z=ab^2a\underline x$};
  \node[fill = white] at (3.2,3.65) {$ b^2ab^2a\underline x$};
  \node[fill = white]  at (4.6,2.1) {$ b^2aba\underline x$};
  \node[fill = white]  at (4.6,-1.9) {$ bab^2a\underline x$};
  \node[fill = white]  at (2.6,-4.05) {$ baba\underline x$};
  \node[fill = white]  at (-3,3.7) {$ abab\underline x$};
  \node[fill = white]  at (-4.6,2.1) {$ abab^2\underline x$};
  \node[fill = white]  at (-4.6,-1.9) {$ ab^2ab\underline x$};
  \node[fill = white!50]  at (-2.8,-4.05) {$ ab^2ab^2\underline x$};
 \fill[black]  (0,0) circle[radius=.15];
     \fill[black]  (1.76,1.76) circle[radius=.15];
     \fill[black]  (1.76,-1.76) circle[radius=.15];
     \fill[black]  (-1.76,1.76) circle[radius=.15];
     \fill[black]  (-1.76,-1.76) circle[radius=.15];
     \draw [very thick, black] (-1.76,-1.76)  to[out=45+15, in=-135-15] (0,0) to[out=45-15, in=-135+15](1.76,1.76);
     \draw [very thick, black] (1.76,-1.76) to[out=135-15, in=-45+15] (0,0) to[out=135+15, in=-45-15] (-1.76,1.76);
     \draw [very thick, black] (1.76,1.76) arc[start angle=180-30, end angle= 180+30, radius= 3.5];
     \draw [very thick, black] (-1.76,-1.76) arc[start angle=-30, end angle= +30, radius= 3.5];
    \draw [very thick, black] (1.76,1.76) arc[start angle=-30, end angle= 15, radius= 2.4]
     arc[start angle=180+15, end angle= 270-15, radius= 2.2];
    \draw [very thick, black] (1.76,1.76) arc[start angle=120, end angle= 75, radius= 2.4];
     \fill[black]  (1.76,1.76) circle[radius=.15];
     \fill[black]  (2.01,3.57) circle[radius=.15];
     \fill[black]  (3.57,2.01) circle[radius=.15];
     \begin{scope}[rotate=90]
        \draw [very thick, black] (1.76,1.76) arc[start angle=-30, end angle= 15, radius= 2.4]
     arc[start angle=180+15, end angle= 270-15, radius= 2.2];
    \draw [very thick, black] (1.76,1.76) arc[start angle=120, end angle= 75, radius= 2.4];
        \fill[black]  (1.76,1.76) circle[radius=.15];
        \fill[black]  (2.01,3.57) circle[radius=.15];
        \fill[black]  (3.57,2.01) circle[radius=.15];
     \end{scope}
     \begin{scope}[rotate=180]
        \draw [very thick, black] (1.76,1.76) arc[start angle=-30, end angle= 15, radius= 2.4]
     arc[start angle=180+15, end angle= 270-15, radius= 2.2];
    \draw [very thick, black] (1.76,1.76) arc[start angle=120, end angle= 75, radius= 2.4];
        \fill[black]  (1.76,1.76) circle[radius=.15];
        \fill[black]  (2.01,3.57) circle[radius=.15];
        \fill[black]  (3.57,2.01) circle[radius=.15];
     \end{scope}
     \begin{scope}[rotate=270]
        \draw [very thick, black] (1.76,1.76) arc[start angle=-30, end angle= 15, radius= 2.4]
     arc[start angle=180+15, end angle= 270-15, radius= 2.2];
    \draw [very thick, black] (1.76,1.76) arc[start angle=120, end angle= 75, radius= 2.4];
        \fill[black]  (1.76,1.76) circle[radius=.15];
        \fill[black]  (2.01,3.57) circle[radius=.15];
        \fill[black]  (3.57,2.01) circle[radius=.15];
     \end{scope}
       \begin{scope}[opacity=.3]
            \draw [dgray] (0,-6) -- (0,6);
        \draw [dgray] (0,0) circle (6);
        \foreach \i in {0,1,2,3} {%
            \draw [ dgray] (90*\i:6) arc (270+90*\i:180+90*\i:6);}
         \foreach \i in {0,1,...,7} {%
             \draw [dgray] (45*\i:6) arc (270+45*\i:135+45*\i:3.3*.75);}
        \foreach \i in {0,1,...,15} {%
             \draw [dgray] (22.5*\i:6) arc (270+22.5*\i:112.5+22.5*\i:1.6*.75);}
         \foreach \i in {0,1,...,31} {%
             \draw [thin] (11.25*\i:6) arc (270+11.25*\i:101.25+11.25*\i:0.8*.75);}
        % \foreach \i in {0,1,...,63} {%
        %     \draw [ultra thin] (5.625*\i:8) arc (270+5.625*\i:95.625+5.625*\i:0.4);}
        \end{scope}
  \end{tikzpicture}
  \end{center}
  \caption{The Cayley graph of $\Gamma_2$ for the generating set $\{b, b^2, aba, ab^2a\}$, with vertices the $\mathrm{PSL}_2(\mathbb Z)$-orbit of $\underline x$ (which equals its $\Gamma_2$-orbit).}
  \label{Figure:CayleyGamma2}
  \end{figure}

\begin{figure}[ht]
\begin{center}
  \begin{tikzpicture}[line join = round, line cap = round, scale=.9] 
     \fill[black]  (0,0) circle[radius=.15];
     \draw [ultra thick, black] (-1.76,-1.76)  to[out=45+15, in=-135-15] (0,0) to[out=45-15, in=-135+15] (1.76,1.76);
     \draw [ultra thick, black] (1.76,-1.76) to[out=135-15, in=-45+15] (0,0) to[out=135+15, in=-45-15] (-1.76,1.76);
     \draw [dashed, black, opacity=.7] (1.76,1.76) arc[start angle=180-30, end angle= 180+30, radius= 3.5];
     \draw [dashed, black, opacity=.7] (-1.76,-1.76) arc[start angle=-30, end angle= +30, radius= 3.5];
    \draw [dashed, black, opacity=.7] (1.76,1.76) arc[start angle=-30, end angle= 15, radius= 2.4];
    \draw [ultra thick, black] (1.76,1.76) arc[start angle=120, end angle= 75, radius= 2.4]
    arc[start angle=270-15, end angle= 180+15, radius= 2.2];
     \fill[white]  (1.76,1.76) circle[radius=.08];
     \fill[black, opacity=.4]  (1.76,1.76) circle[radius=.08];
     \fill[white]  (2.01,3.57) circle[radius=.08];
     \fill[black, opacity=.4]  (2.01,3.57) circle[radius=.08];
     \fill[black]  (3.57,2.01) circle[radius=.15];
     \begin{scope}[rotate=90]
        \draw [ultra thick, black] (1.76,1.76) arc[start angle=-30, end angle= 15, radius= 2.4]
     arc[start angle=180+15, end angle= 270-15, radius= 2.2];
    \draw [dashed, black, opacity=.7] (1.76,1.76) arc[start angle=120, end angle= 75, radius= 2.4];
        \fill[white]  (1.76,1.76) circle[radius=.08];
        \fill[black, opacity=.4]  (1.76,1.76) circle[radius=.08];
        \fill[black]  (2.01,3.57) circle[radius=.15];
        \fill[white]  (3.57,2.01) circle[radius=.08];
        \fill[black, opacity=.4]  (3.57,2.01) circle[radius=.08];
     \end{scope}
     \begin{scope}[rotate=180]
        \draw [dashed, black, opacity=.7] (1.76,1.76) arc[start angle=-30, end angle= 15, radius= 2.4];
    \draw [ultra thick, black] (1.76,1.76) arc[start angle=120, end angle= 75, radius= 2.4]
         arc[start angle=270-15, end angle= 180+15, radius= 2.2];
        \fill[white]  (1.76,1.76) circle[radius=.08];
        \fill[black, opacity=.4]  (1.76,1.76) circle[radius=.08];
        \fill[white]  (2.01,3.57) circle[radius=.08];
        \fill[black, opacity=.4]  (2.01,3.57) circle[radius=.08];
        \fill[black]  (3.57,2.01) circle[radius=.15];
     \end{scope}
     \begin{scope}[rotate=270]
        \draw [ultra thick, black] (1.76,1.76) arc[start angle=-30, end angle= 15, radius= 2.4]
     arc[start angle=180+15, end angle= 270-15, radius= 2.2];
    \draw [dashed, black, opacity=.7] (1.76,1.76) arc[start angle=120, end angle= 75, radius= 2.4];
        \fill[white]  (1.76,1.76) circle[radius=.08];
        \fill[black, opacity=.4]  (1.76,1.76) circle[radius=.08];
        \fill[black]  (2.01,3.57) circle[radius=.15];
        \fill[white]  (3.57,2.01) circle[radius=.08];
        \fill[black, opacity=.4]  (3.57,2.01) circle[radius=.08];
     \end{scope}
  \node at (.5,0) {$\underline x$};
  \node[opacity=.7] at (2.5,-1.5) {$\underline y=b\underline x$};
  \node[opacity=.7] at (2.5,1.6) {$\underline z=b^2\underline x$};
  \node[opacity=.7] at (-0.7,1.9) {$a\underline y=aba\underline x$};
  \node[opacity=.7] at (-0.6,-1.8) {$a\underline z=ab^2a\underline x$};
  \node at (4.5, -2) {$bab^2a\underline x$};
  \node at (4.5, 2.1) {$b^2aba\underline x$};
  \node at (-4.5, 2.1) {$abab^2\underline x$};
  \node at (-4.5, -2) {$ab^2ab\underline x$};
  \node[opacity=.7] at (2.8,3.7) {$b^2ab \underline y$};
  \node[opacity=.7] at (2.8,-3.5) {$bab^2 \underline z$};
  \node[opacity=.7] at (-2.8,3.7) {$abab^2 \underline z$};
  \node[opacity=.7] at (-2.8,-3.5) {$ab^2ab \underline y$};
  \end{tikzpicture}
  \end{center}
  \caption{The Cayley graph of $\Gamma_6$ for the generating set
  $\{b^2aba, bab^2a\}$. The figure shows 4 edges, 5 vertices, 4 half edges and 
  8 midpoints.}
  \label{Figure:CayleyGamma6}
  \end{figure}

We consider the action of 
$\Gamma= \mathrm{PSL}_2(\mathbb Z)
<
\mathrm{PSL}_2(\mathbb R)
$ on 
$\mathbb H^2$.
In analogy with the action on $X$, we denote
$$
\underline x=\mathrm{Fix}(a)\in\mathbb H^2, \qquad
\underline y=b \underline x\in\mathbb H^2, 
\qquad \underline z=b^2\underline x\in\mathbb H^2.
$$
Each point $\underline x$, $\underline y$, and $\underline z$ lies
on an edge of the Farey tessellation.

Consider the whole $\mathrm{PSL}_2(\mathbb Z)$-orbit of the 
triangle $\underline x\underline y\underline z\in\mathbb H^2$;
it can be seen as
 the Cayley graph of $\Gamma_2$.
Figure~\ref{Figure:CayleyGamma2} shows the Cayley graph of $\Gamma_2$, 
with one vertex precisely at each 
$\mathrm{PSL}_2(\mathbb Z)$-orbit
of $\underline x$
(the $\Gamma_2$-orbit of $\underline x$ equals the $\mathrm{PSL}_2(\mathbb Z)$-orbit
of $\underline x$). 

The Cayley graph of $\Gamma_6\cong F_2$ embeds in $\mathbb H^2$ with 
vertices the $\Gamma_6$-orbit of $\underline x$ and edges 
the union of two segments, in the $\Gamma_3$-orbit 
of the segments $\underline x\underline y$ and $\underline x\underline z$ 
(the $\Gamma_3$-orbit of $\underline y\underline z$ does not occur). 
Points in the $\Gamma_3$-orbits of $\underline y$ and $\underline z$  
are precisely the midpoints of the embedded Cayley graph of 
$\Gamma_6\cong F_2$, Figure~\ref{Figure:CayleyGamma6}. We have:

\begin{Remark}
    \label{Rmk:midpoints}
    The midpoint of each  edge of the Cayley graph of $\Gamma_6$ embedded in $\mathbb{H}^2$ as above
    is a vertex of the  Cayley graph of $\Gamma_2$.
\end{Remark}

\subsection{Orbit map of the Cayley graph}

The orbit map yields a local embedding of the Cayley graph of $\mathrm{PSL}_2(\mathbb Z)$
in $X$. 
We shall use Theorem~\ref{Thm:LocalMorse} 
to show that this local embedding maps geodesics in the
Cayley graph to regular quasi-geodesics in $X$.
We shall not apply Theorem~\ref{Thm:LocalMorse} to the vertices
of this graph but to the  ``midpoints'' (in the sense of 
Remark~\ref{Rmk:midpoints})
of the Cayley graph of
$$
\Gamma_6=\ker(\mathrm{PSL}_2(\mathbb Z)\cong \mathbb Z_2*\mathbb Z_3
\twoheadrightarrow \mathbb Z_2\times\mathbb Z_3)\cong F_2,
$$
which is the group freely generated by 
$ bab^2a$ and  $b^2aba$.

Recall that $x\in \mathrm{Fix}(\rho(a))$ and consider the orbit map
$$
\begin{array}{rcl}
     \Gamma_6 & \to & X  \\
     g & \mapsto & \rho(g) x
\end{array}.
$$
Given a discrete geodesic $\{\dotsc,g_{n-1},g_n,g_{n+1},\dotsc\}$ 
in $\Gamma_6$, consider 
the midpoints in  $\Gamma_2$ obtained by 
taking midpoints in the Cayley graph as in 
Remark~\ref{Rmk:midpoints}, 
$\mathrm{mid}(g_n, g_{n+1})\in\Gamma_2$, and consider their 
orbits
$$
m_n=\rho(\mathrm{mid}(g_n, g_{n+1}))(x).
$$
(Recall that $\mathrm{mid}(g_n, g_{n+1})\in \Gamma_2$
and   $\mathrm{mid}(g_n, g_{n+1}) \underline x\in\mathbb H^2$
is the midpoint between $g_n\underline x$ and $ g_{n+1}\underline x$
in the Cayley graph of $\Gamma_2$ embedded in 
$\mathbb{H}^2$).  

\begin{figure}[ht]
\begin{center}
  \begin{tikzpicture}[line join = round, line cap = round, scale=1.4] 
        \draw[thick] (0,0)  to[out=-25,in=180-2] (1,-.25) to[out=+30,in=180+5] (2,0) to[out=-30,in=180-5] (3,-.5) to[out=+45,in=180+5] (4,0)
        to[out=-10,in=180-30] (5,-.3) 
        to[out=+30,in=180+15] (6,0);
        \fill[black]  (0,0) circle[radius=.04];
        \fill[black]  (2,0) circle[radius=.04];
        \fill[black]  (4,0) circle[radius=.04];
        \fill[black]  (6,0) circle[radius=.04];
        \fill[black]  (1,-.25) circle[radius=.04];
        \fill[black]  (3,-.5) circle[radius=.04];
        \fill[black]  (5,-.3) circle[radius=.04];
         \draw[ultra thick, opacity=.3] (1,-.25)  -- (3,-.5)-- (5,-.3) ;
         \node at (0,.25) {$\rho(g_{n+2})x$};
         \node at (2,.25) {$\rho(g_{n+1})x$};
         \node at (4,.25) {$\rho(g_{n})x$};
         \node at (6,.25) {$\rho(g_{n-1})x$};
         \node at (1,-.45) {$m_{n+1}$};
         \node at (3,-.7) {$m_{n}$};
         \node at (5,-.5) {$m_{n-1}$};
  %    \fill[black]  (0,0) circle[radius=.08]
  %    \draw[black,thick] (0,0) circle[radius=.08];
  %    \draw [thick, black] (-1.76,-1.76) to[out=45+15, in=-135-15] (0,0) to[out=45-15, in=-135+15] (1.76,1.76);
  %    \draw [thick, black] (1.76,-1.76) to[out=135-15, in=-45+15] (0,0) to[out=135+15, in=-45-15](-1.76,1.76);
  %    \draw [dashed, black, opacity=.7] (1.76,1.76) arc[start angle=180-30, end angle= 180+30, radius= 3.5];
  %    \draw [dashed, black, opacity=.7] (-1.76,-1.76) arc[start angle=-30, end angle= +30, radius= 3.5];
  %    \fill[white]  (1.76,1.76) circle[radius=.08];
  %    \fill[black]  (1.76,-1.76) circle[radius=.08];
  %    \fill[white]  (-1.76,1.76) circle[radius=.08];
  %    \fill[white]  (-1.76,-1.76) circle[radius=.08];
  %    \draw[black,thick]  (1.76,1.76) circle[radius=.08];
  %    \draw[black,thick]  (1.76,-1.76) circle[radius=.08];
  %    \draw[black,thick]  (-1.76,1.76) circle[radius=.08];
  %    \draw[black,thick]  (-1.76,-1.76) circle[radius=.08];
  % % 
  % \fill[white] (1.3,0) circle[radius=.3];
  % \node at (1.2,0) {$x_{n+1}=x$};
  % \node[opacity=1] at (2.7,-1.8) {$m_{n}=y$};
  % \node[opacity=1] at (2.9,1.8) {$m_{n+1}=z$};
  % \node[opacity=1] at (-3.3,1.8) {$m_{n+1}=\rho(a)y$};
  % \node[opacity=1] at (-3.3,-1.8) {$m_{n+1}=\rho(a)z$};
  \end{tikzpicture}
  \end{center}
  \caption{The path of orbits 
  $\rho(g_n)x$ 
  and
  the path orbits of midpoints 
  $m_n=\rho(\mathrm{mid}(g_n,g_{n+1} ))(x)$.}
  \label{Figure:midpointspath}
  \end{figure}

\begin{Proposition}
\label{Prop:StraightlSpaced}
    Given $C,l,\varepsilon>0$ (with $s\leq C$) and $\Theta\subset \mathrm{int}(\sigma_{\mathrm{mod}})$ (a $\iota$-symmetric, compact interval with
    nonempty interior),
    for $t$ sufficiently large the sequence $(m_n)_{n\in\mathbb Z}$ is
   $(\Theta, \varepsilon)$-straight and $l$-spaced (Definition~\ref{Def:straight}). 
\end{Proposition}

\begin{proof}
% \smallskip
% {\emph{The Cayley graph of $\Gamma_6$ in $X$.}} 
To prove Proposition~\ref{Prop:StraightlSpaced}, 
we use  the 
local embedding in $X$ of the Cayley graphs of $\Gamma_6$
and $\Gamma_2$
via the $\rho(\Gamma_2)$-orbit of  $x$.
 Consider 
the discrete geodesic 
$$\dotsc,g_{n-1},g_n,g_{n+1},\dotsc$$ 
in $\Gamma_6$
and the corresponding sequence of orbits
$$\dotsc,x_{n-1},x_n,x_{n+1},\dotsc$$ 
with $x_j=\rho(g_j)x$, and 
orbits of
midpoints in the sense of Remark~\ref{Rmk:midpoints}
$m_j=\rho(\mathrm{mid}(g_j,g_{j+1}))(x)$,
see Figure~\ref{Figure:midpointspath}. 
We prove first the following estimate:

\begin{Lemma}
    \label{Lem:angleestimate} 
   For sufficiently large $t$, $\angle_{m_{n}}(x_{n+1},m_{n+1})$
is arbitrarily small.
\end{Lemma}

\begin{figure}[ht]
\begin{center}
  \begin{tikzpicture}[line join = round, line cap = round, scale=.9] 
     \fill[black]  (0,0) circle[radius=.08];
     \draw[black,thick] (0,0) circle[radius=.08];
     \draw [thick, black] (-1.76,-1.76) to[out=45+15, in=-135-15] (0,0) to[out=45-15, in=-135+15] (1.76,1.76);
     \draw [thick, black] (1.76,-1.76) to[out=135-15, in=-45+15] (0,0) to[out=135+15, in=-45-15](-1.76,1.76);
     \draw [dashed, black, opacity=.7] (1.76,1.76) arc[start angle=180-30, end angle= 180+30, radius= 3.5];
     \draw [dashed, black, opacity=.7] (-1.76,-1.76) arc[start angle=-30, end angle= +30, radius= 3.5];
     \fill[white]  (1.76,1.76) circle[radius=.08];
     \fill[black]  (1.76,-1.76) circle[radius=.08];
     \fill[white]  (-1.76,1.76) circle[radius=.08];
     \fill[white]  (-1.76,-1.76) circle[radius=.08];
     \draw[black,thick]  (1.76,1.76) circle[radius=.08];
     \draw[black,thick]  (1.76,-1.76) circle[radius=.08];
     \draw[black,thick]  (-1.76,1.76) circle[radius=.08];
     \draw[black,thick]  (-1.76,-1.76) circle[radius=.08];
  \fill[white] (1.3,0) circle[radius=.3];
  \node at (1.2,0) {$x_{n+1}=x$};
  \node[opacity=1] at (2.7,-1.8) {$m_{n}=y$};
  \node[opacity=1] at (2.9,1.8) {$m_{n+1}=z$};
  \node[opacity=1] at (-3.3,1.8) {$m_{n+1}=\rho(a)y$};
  \node[opacity=1] at (-3.3,-1.8) {$m_{n+1}=\rho(a)z$};
  \end{tikzpicture}
  \end{center}
  \caption{The points $m_{n}=y$, $x_{n+1}=x$ and the three possibilities for $m_{n+1}$: $m_{n+1}=z$, $m_{n+1}=\rho(a)y$, or $m_{n+1}=\rho(a)z$.}
  \label{Figure:midpoint}
  \end{figure}

 \begin{proof}
     By equivariance, we assume that $x_{n+1}=x$ and $m_{n}=y$. 
There are three possibilities for 
$m_{n+1}$: $m_{n+1}=z$, $m_{n+1}=\rho(a)y$, or $m_{n+1}=\rho(a)z$ (see Figure~\ref{Figure:midpoint}).
\begin{itemize}
    \item If $m_{n+1}=z$, then $\angle_{y}(x,z)$ is arbitrarily small 
    by Lemma~\ref{lem:smallangles}.
    \item If $m_{n+1}=\rho(a)y$, then
    by 
    Lemma~\ref{Lemma:almostperp}
    $\angle_x (y, \rho(a)y)\approx\pi$ 
     and, as the angles add at most $\pi$, $\angle_{y}(x,\rho(a)y)\approx 0$.
    \item If $m_{n+1}=\rho(a)z$, then the argument is as follows. 
    By Lemma~\ref{lem:smallangles}
    $$
       \angle_{x}(\rho(a)(y),\rho(a)(z))= 
       \angle_{x}(y,z)\approx 0.
    $$
    Thus
    $$
    \angle_{x}(y,\rho(a)z)\geq \angle_{x}(y,\rho(a)y)-
    \angle_{x}(\rho(a)y,\rho(a)z)\approx\pi
    $$
    because $\angle_{x}(y,\rho(a)y)\approx \pi$,
    Corollary~\ref{Corollary:xyay}. 
    Hence, as the addition of angles of a triangle is at most  $\pi$, we have
    $\angle_{y}(x,\rho(a)z)\approx 0$.
\end{itemize}
This proves Lemma~\ref{Lem:angleestimate}. 
 \end{proof}

% Since both $\angle_{m_{n}}(x_{n+1},m_{n+1})$ and 
% $\angle_{m_{n}}(x_{n},m_{n-1})$ are arbitrarily small, and using that  
%  $m_{n}=\mathrm{mid}(x_n, x_{n+1}) $,

We continue the proof of Proposition~\ref{Prop:StraightlSpaced}.
By Lemma~\ref{Lem:angleestimate},  
$\angle_{m_{n}}(x_{n+1},m_{n+1})\approx 0$,
 and the same argument yields  
$\angle_{m_{n}}(x_{n},m_{n-1})\approx 0$. 
%, Figure~\ref{Figure:midpointspath}. 
Furthermore, 
$\angle_{m_n}(x_n, x_{n+1})\approx\pi$,
by Corollary~\ref{Corollary:xyay}. 
It follows from these estimates that 
$$
\angle_{m_n}(m_{n-1}, m_{n+1})\geq \angle_{m_n}(x_n, x_{n+1})
-\angle_{m_n}(x_n, m_{n-1})-\angle_{m_n}(x_{n+1},m_{n+1})
\approx \pi.
$$  
% because  $\angle_{m_n}(x_n, x_{n+1})=\pi$, 
% as $m_{n}=\mathrm{mid}(x_n, x_{n+1}) $.

On the other hand, the type of the segment 
$yx$
 is arbitrarily 
close to the center of the Weyl chamber $\zeta=\mathrm{Fix}(\iota)=\mathrm{center}(\sigma_\mathrm{mod})$, 
by Lemma~\ref{Lem:regular}.
In our sequence, this yields that
$\mathrm{type}(m_n x_{n+1})\approx \zeta$.
(Notice that the type map is 
$\mathrm{SL}_3(\mathbb R)$-invariant, 
so we can work equivariantly to get uniform estimates
on the type map.)
By Lemma~\ref{Lem:angleestimate}
$\angle_{m_{n}}(x_{n+1},m_{n+1})\approx 0$,
therefore $\mathrm{type}(m_n m_{n+1})\approx \zeta$.
In particular,   $ m_{n} m_{n+1} $ is $\Theta$-regular. 
% The argument in the proof of Lemma~\ref{Lem:regular}
% shows that $\angle_x(y, \xi)$ is arbitrarily small, 
% where $\xi$ is the limit of the segment 
% starting at $\bar x$ and containing $\bar y$ 
% (and has type $\zeta=\mathrm{Fix}(\iota)$). 
As $\mathrm{type}(m_n m_{n+1})\approx \zeta$
and similarly $\mathrm{type}(m_n m_{n-1})\approx \zeta$,
the estimate on the angle
$\angle_{m_n}(m_{n-1},m_{n+1})\cong \pi$
yields the assertion on the 
$\zeta$-angles: 
$$\angle^\zeta_{m_n}(m_{n-1},m_{n+1})\cong \pi .
$$
This proves straightness.

To establish that the sequence is $l$-spaced, 
as in the proof of Lemma~\ref{Lem:angleestimate}
we may
assume $m_n=y$, $x_{n+1}=x$ and
distinguish the three possible positions of $m_{n+1}$ 
(either $z$, $\rho(a)y$, or $\rho(a)z$), 
% as in the proof of Lemma~\ref{Lem:angleestimate}, 
cf. Figure~\ref{Figure:midpoint}. The claim is immediate when 
  $m_{n+1}=z$ since  
$d(m_n,m_{n+1}) $ is 
the length of an edge of the triangle $xyz$, which
goes to $+\infty$ as $t\to+\infty$.

It remains to handle the cases $m_{n+1}= \rho(a)z$
and $m_{n+1}=\rho(a)y$.  In both cases
we have checked that 
$\angle_x(y,m_{n+1})\approx \pi$. As
$X$ is a 
Cartan-Hadamard manifold, the \emph{first cosine inequality} in \cite[Proposition I.5.1]{Ballmann} applied to the
triangle
$x,y,m_{n+1}$ yields that $d(y,m_{n+1})>d(y,x)$,
%
%
% We claim that the distance to $y$ is strictly increasing along the oriented segment from $x$ to $\rho(a)z$. First, because $\angle_x(y,\rho(a)z)\approx \pi$, the first variation formula implies that the distance to $y$ has positive derivative at the initial point $x$.
% Next, for any point  $w$ in the interior of the segment from $x$ to $\rho(a)z$
% the angle is even larger,
% \[
% \angle_w(y,\rho(a)z) > \angle_x(y,\rho(a)z),
% \]
% because the sum of angles of the triangle $yxw$ is at most $\pi$, and both 
% \( 
% \angle_x(y,w) = \angle_x(y,\rho(a)z)\)
% and
% % \quad\text{and}\quad
% \(\angle_w(x,y) = \pi - \angle_w(y,\rho(a)z)
% \)
% are angles of this triangle. Thus by the first variation formula, the distance to $y$ is strictly increasing at every interior point of the segment from $x$ to $\rho(a)z$. This  yields $d(y,\rho(a)z)>d(y,x)$, 
which goes to $+\infty$ with $t$.

This completes the proof of Proposition~\ref{Prop:StraightlSpaced}.
\end{proof}

To conclude the proof of Theorem~\ref{Thm:Anosov}, we use the
characterization of Anosov groups 
as Morse groups.

\begin{Definition}
    Given two points $x_\pm\in X$ so that the segment $x_-x_+$ is regular, the
    \emph{diamond}  with tips $x_-$ and $x_+$ is
$$
\Diamond(x_-,x_+)= V(x_-,\sigma_+)\cap V(x_+,\sigma_-),
$$
where $V(x_-,\sigma_+)$ is the Euclidean Weyl sector with tip $x_-$ that contains $x_+$, and analogously for $V(x_+,\sigma_-)$.
\end{Definition}

\begin{figure}[ht]
\begin{center}
  \begin{tikzpicture}[line join = round, line cap = round, scale=1] 
    \draw[white, fill=gray!30!white, opacity=0.3]  (1.5,1)--(-1,0)--(1.5,-1);
    \draw[white, fill=gray!30!white, opacity=0.3]  (1.5,1) to[out=-30, in=30](1.5,-1);
    \draw[white, fill=gray!30!white, opacity=0.3]  (-1.5,1)--(1,0)--(-1.5,-1);
    \draw[white, fill=gray!30!white, opacity=0.3]  (-1.5,1) to[out=240, in=120] (-1.5,-1);
    \draw[ ultra thick, black, fill=gray!50!white, opacity=0.6] (-1,0)--(0,0.4)--(1,0)--(0,-.4)--(-1,0);
    \draw[thick, black]  (-1,0)--(1.5,1);
    \draw[thick, black]  (-1,0)--(1.5,-1);
    \draw[thick, black]  (1,0)--(-1.5,1);
    \draw[thick, black]  (1,0)--(-1.5,-1);
     \node at (1.3,0) {$x_+$};
     \node at (-1.3,0) {$x_-$};
     \node at (2.5,.5) {$V(x_-,\sigma_+)$};
       \node at (-2.5,.5) {$V(x_+,\sigma_-)$};
  %    \fill[black]  (0,0) circle[radius=.08];
  %    \draw[black,thick] (0,0) circle[radius=.08];
  %    \draw [thick, black] (-1.76,-1.76) to[out=45+15, in=-135-15] (0,0) to[out=45-15, in=-135+15] (1.76,1.76);
  %    \draw [thick, black] (1.76,-1.76) to[out=135-15, in=-45+15] (0,0) to[out=135+15, in=-45-15](-1.76,1.76);
  %    \draw [dashed, black, opacity=.7] (1.76,1.76) arc[start angle=180-30, end angle= 180+30, radius= 3.5];
  %    \draw [dashed, black, opacity=.7] (-1.76,-1.76) arc[start angle=-30, end angle= +30, radius= 3.5];
  %    \fill[white]  (1.76,1.76) circle[radius=.08];
  %    \fill[black]  (1.76,-1.76) circle[radius=.08];
  %    \fill[white]  (-1.76,1.76) circle[radius=.08];
  %    \fill[white]  (-1.76,-1.76) circle[radius=.08];
  %    \draw[black,thick]  (1.76,1.76) circle[radius=.08];
  %    \draw[black,thick]  (1.76,-1.76) circle[radius=.08];
  %    \draw[black,thick]  (-1.76,1.76) circle[radius=.08];
  %    \draw[black,thick]  (-1.76,-1.76) circle[radius=.08];
  % % 
  % \fill[white] (1.3,0) circle[radius=.3];
  % \node at (1.2,0) {$x_{n+1}=x$};
  % \node[opacity=1] at (2.7,-1.8) {$m_{n}=y$};
  % \node[opacity=1] at (2.9,1.8) {$m_{n+1}=z$};
  % \node[opacity=1] at (-3.3,1.8) {$m_{n+1}=\rho(a)y$};
  % \node[opacity=1] at (-3.3,-1.8) {$m_{n+1}=\rho(a)z$};
  \end{tikzpicture}
  \end{center}
  \caption{The diamond $\Diamond(x_-,x_+)= V(x_-,\sigma_+)\cap V(x_+,\sigma_-)$}
  \label{Figure:Diamond}
  \end{figure}

\begin{Definition}
    A discrete subgroup  
    $  H< \mathrm{Isom}(X)$ is 
    $\sigma_{\mathrm{mod}}$-Morse if it is 
    regular, word
hyperbolic, and satisfies the following property:

For every discrete geodesic segment 
$s\colon [n_-, n_+] \cap \mathbb Z\to H$, the path of orbits 
$sx$ is contained in a tubular neighborhood of uniform radius
$R$ (depending on $H$ and $x$) of a diamond
$\Diamond(x_-,x_+)$ with tips at distance $d(x_{\pm }, s(n_{\pm}) x)\leq R$ from 
the end-points.
\end{Definition}

\begin{proof}[Proof of Theorem~\ref{Thm:Anosov}]
By  Proposition~\ref{Prop:StraightlSpaced}, 
Theorem~\ref{Thm:LocalMorse} applies to the sequence of midpoints
of orbits of any discrete geodesic sequence in $\Gamma_6\cong F_2$.

Theorem~\ref{Thm:LocalMorse}
% This implies that the orbit of a discrete geodesic sequence 
% in $\Gamma_6\cong F_2$ is a quasi-geodesic in $X$.  Thus $\rho(\Gamma_6)$ is a 
% discrete subgroup, and
% by Proposition~\ref{Prop:StraightlSpaced} and 
% Theorem~\ref{Thm:LocalMorse} again, $\rho(\Gamma_6)$ 
% is a Morse group.
gives Morse control for midpoint sequences; since vertices stay 
uniformly close to their adjacent midpoints, the orbit map of 
\(\Gamma_6\) is Morse; hence \(\rho(\Gamma_6)\) is discrete and 
Morse 
and, since $\Gamma_6$ has finite index in $\Gamma$,  
the conclusion follows for
\(\rho(\Gamma)\).
Morse is one of the characterizations 
of Anosov subgroups~\cite{KLPEJM}.
\end{proof}

\bibliographystyle{abbrv}
\bibliography{Modular}

\begin{small}
\noindent \textsc{Departament de Matem\`atiques,  Universitat Aut\`onoma de Barcelona,
and Centre de Recerca Matem\`atica (UAB-CRM)\\
08193 Cerdanyola del Vall\`es, Spain }

\noindent \textsf{joan.porti@uab.cat}\\
\noindent \textsf{joanporti@gmail.com}
\end{small}

\end{document}